\documentclass[11pt]{report}
\usepackage[margin=1.5in]{geometry}
\usepackage{amssymb,amsfonts,amsmath,amsthm,amscd,dsfont,mathrsfs,bbold,pifont}
\usepackage{blkarray}
\usepackage{graphicx,float,psfrag,epsfig,color}
\usepackage{comment}
\usepackage{subcaption}

\usepackage{algorithm}
\usepackage[noend]{algpseudocode}

\makeatletter
\def\BState{\State\hskip-\ALG@thistlm}
\makeatother

	\usepackage{microtype}
	\usepackage[pdftex,pagebackref=false,colorlinks]{hyperref}
	\hypersetup{linkcolor=[rgb]{.7,0,0}}
	\hypersetup{citecolor=[rgb]{0,.7,0}}
	\hypersetup{urlcolor=[rgb]{.7,0,.7}}

\newcommand\independent{\protect\mathpalette{\protect\independent}{\perp}} 
\def\independent#1#2{\mathrel{\rlap{$#1#2$}\mkern2mu{#1#2}}}

\def\mmse{{\sf mmse}}
\def\MMSE{{\sf MMSE}}

\def\overlap{{\sf Overlap}}
\def\sT{{\sf T}}
\def\hx{\widehat{x}}
\def\cG{\mathcal{G}}
\def\reals{\mathbb{R}}
\def\Unif{{\sf Unif}}
\def\normal{{\sf N}}
\def\Info{{\sf I}}
\def\de{{\rm d}}

\def\hbs{{\boldsymbol{\widehat{s}}}}
\def\<{\langle}
\def\>{\rangle}

\newcommand{\snr}{\mathrm{SNR}}
\newcommand{\sbm}{\mathrm{SBM}}
\newcommand{\ssbm}{\mathrm{SSBM}}
\newcommand{\map}{\mathrm{map}}

\newcommand{\pin}{q_\mathrm{in}}
\newcommand{\pou}{q_\mathrm{out}}
\newcommand{\sgn}{\mathrm{sgn}}

\newcommand{\bin}{\mathrm{Bin}}

\newcommand{\argmax}{\mathrm{argmax}}

\newcommand{\mR}{\mathbb{R}} 

\newcommand{\mZ}{\mathbb{Z}}

\newcommand{\tr}{\mathrm{tr}}

\newcommand{\pp}{\mathbb{P}}
\newcommand{\E}{\mathbb{E}}
\newcommand{\Var}{\mathbb{V}\mathrm{ar}}
\newcommand{\e}{\varepsilon}

\newcommand{\rank}{\mathrm{rank}}

\DeclareMathOperator{\diag}{diag}

\newcommand{\X}{\mathcal{X}}
\newcommand{\Y}{\mathcal{Y}}

\newcommand{\1}{\mathbb{1}}

\newtheorem{theorem}{Theorem}
\newtheorem{lemma}{Lemma}
\newtheorem{corollary}{Corollary}
\newtheorem{definition}{Definition}
\newtheorem{remark}{Remark}

\newtheorem{conjecture}{Conjecture}

\begin{document}

\title{Community Detection and Stochastic Block Models
}
%\author{Emmanuel Abbe (Princeton University) and Martin Wainwright (UC Berkeley)}

\author{Emmanuel Abbe\thanks{Program in Applied and Computational Mathematics, and Department of Electrical Engineering, Princeton University, Princeton, USA. Email: \texttt{eabbe@princeton.edu}, URL: {\tt www.princeton.edu/$\sim$eabbe}.  
This research was partly supported by the Bell Labs Prize, the NSF CAREER Award CCF-1552131, ARO grant W911NF-16-1-0051, NSF Center for the Science of Information CCF-0939370, and the Google Faculty Research Award. 
} \\
Princeton University
%\and 
%Martin Wainwright\thanks{Departments of Electrical Engineering and Computer Science, and Department of Statistics, University of California at Berkeley, Berkeley, USA, \texttt{wainwrig@berkeley.edu}, {\tt http://www.cs.berkeley.edu/$\sim$wainwrig}.
%}
}

\date{}
%\date{February 28, 2015}
\maketitle

\begin{abstract}
The stochastic block model (SBM) is a random graph model with different group of vertices connecting differently. It is widely employed as a canonical model to study clustering and community detection, and provides a fertile ground to study the information-theoretic and computational tradeoffs that arise in combinatorial statistics and more generally data science. 
 
This monograph surveys the recent developments that establish the fundamental limits for community detection in the SBM, both with respect to information-theoretic and computational tradeoffs, and for various recovery requirements such as exact, partial and weak recovery. The main results discussed are the phase transitions for exact recovery at the Chernoff-Hellinger threshold, the phase transition for weak recovery at the Kesten-Stigum threshold, the optimal SNR-mutual information tradeoff for partial recovery, and the gap between information-theoretic and computational thresholds.

The monograph gives a principled derivation of the main algorithms developed in the quest of achieving the limits, in particular two-round algorithms via graph-splitting, semi-definite programming, (linearized) belief propagation, classical/nonbacktracking spectral methods and graph powering. Extensions to other block models, such as geometric block models, and a few open problems are also discussed.
\end{abstract}

%\chapter{The Distribution and Installation}
%\label{c-intro} % a label for the chapter, to refer to it later
%
%%This document explains to typesetters how to install and use the \LaTeX\ style files for Foundations and Trends journal articles from Now Publishers. 
%
%\chapter{Pre-requisites}
%You will need a working LaTeX installation. We recomend using pdflatex to process the files. You will also need
%biber.exe installed. This is distributed as part of the latest versions of LiveTex and MikTex. If you have
%problems, please let us know.
%
%\chapter{The Distribution}
%The distribution contains 2 folders: \texttt{nowfnt}  and \texttt{nowfnttexmf}. 
%
%\chapter{Folder \texttt{nowfnt}}
%This folder contains the following files using a flat stucture required to compile a FnT issue:
%
%This folder contains the class files both as a texmf structure in the folder \texttt{nowtexmf} and as a flat files.
%The following files are required to compile a FnT issue (comprising an article, a book, and an e-book version):
%
%\begin{itemize}
%\setlength{\parsep}{0pt}
%\setlength{\itemsep}{0pt}
%\item essence\_logo.eps
%\item essence\_logo.pdf
%\item now\_logo.eps
%\item now\_logo.pdf
%\item nowfnt.cls
%\item nowfnt-biblatex.sty
%%\item <jrnlcode>-editorialboard.tex
%%\item <jrnlcode>-journaldata.tex
%\item NOWFnT-data.tex 
%\end{itemize}

\tableofcontents

\chapter{Introduction}
%totodo:
% check proof I2

\section{Community detection, clustering and block models}
The most basic task of {\it community detection}, or {\it graph clustering},
%\footnote{In this note, the terms communities and clusters are used exchangeably.}
%\footnote{Terminology summary: In this note, we view community detection as the general task of inferring similarity classes in labelled graphs. Graph clustering corresponds to a special case of community detection where the communities are assortative, i.e., with denser intra-connectivity. The general clustering problem extends the graph clustering problem by providing not only presence of absence of edges but also labels on the edges (such as similarity or dissimilarity functions) and thus clustering is viewed as a special case of community detection.} 
consists in partitioning the vertices of a graph into clusters that are more densely connected. From a more general point of view, community structures may also refer to groups of vertices that connect similarly to the rest of the graph without having necessarily a higher inner density, such as disassortative communities that have higher external connectivity. Note that the terminology of `community' is sometimes used only for assortative clusters in the literature, but we adopt here the more general definition. Community detection may also be performed on graphs where edges have labels or intensities, and if these labels represent similarities among data points, the problem may be called {\it data clustering.} In this monograph, we will use the terms communities and clusters exchangeably.
Further, one may also have access to interactions that go beyond pairs of vertices, such as in hypergraphs, and communities may not always be well separated due to overlaps. In the most general context, community detection refers to the problem of inferring similarity classes of vertices in a network by having access to  measurements of local interactions.

\begin{figure}[h]
\centering
\begin{subfigure}{.5\textwidth}
  \centering
  \includegraphics[width=1\linewidth]{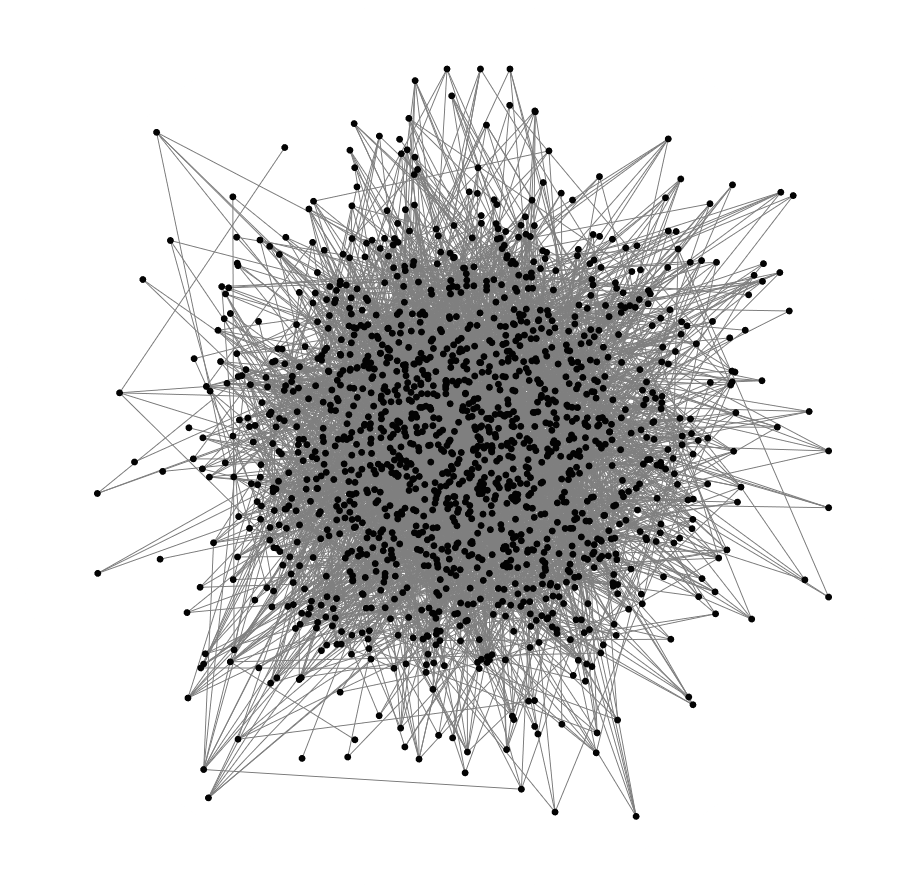}
%  \caption{Random configuration.}
  \label{fig:sub1}
\end{subfigure}%
\begin{subfigure}{.5\textwidth}
  \centering
  \includegraphics[width=1\linewidth]{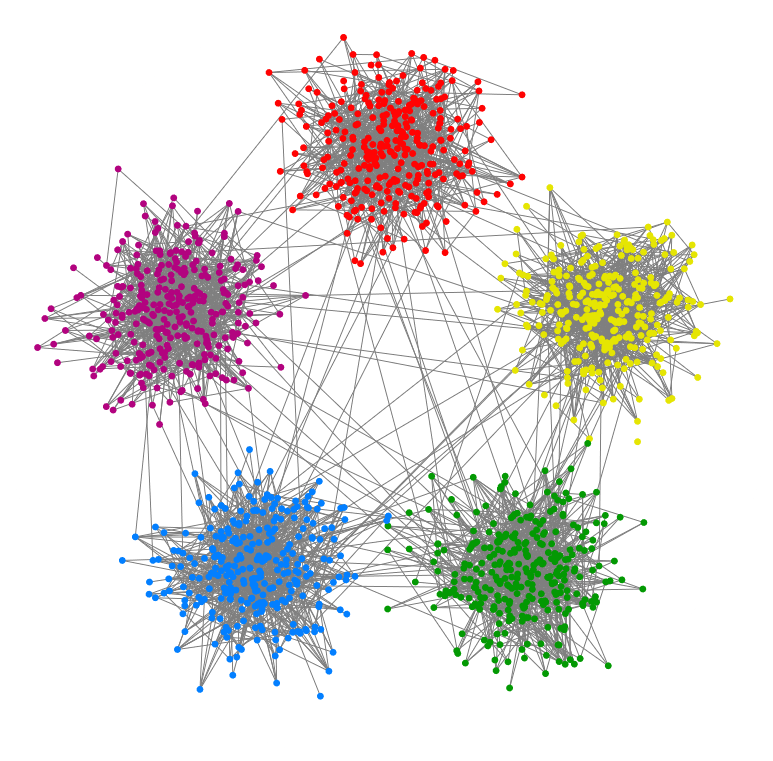}
 % \caption{Clustered configuration.}
  \label{fig:sub2}
\end{subfigure}
\caption{The above two graphs are the same graph re-organized and drawn from the SBM model with 1000 vertices, 5 balanced communities, within-cluster probability of $1/50$ and across-cluster probability of $1/1000$. The goal of community detection in this case is to obtain the right graph (with five communities) from the left graph (scrambled) up to some level of accuracy. In such a context, community detection may be called graph clustering. In general, communities may not only refer to denser clusters but more generally to groups of vertices that behave similarly.}
\label{fig:test1}
\end{figure}

Community detection and clustering are central problems in machine learning and data science. A large number of data sets can be represented as a network of interacting items, and one of the first features of interest in such networks is to understand which items are ``alike,'' as an end or as a preliminary step towards other learning tasks. Community detection is used in particular to understand sociological behavior \cite{airoldi,fortunato,social1}, protein to protein interactions \cite{ppi2,marcotte}, gene expressions \cite{gene-expr2,gene-survey}, recommendation systems \cite{amazon,reco1,srikant-rec}, medical prognosis \cite{tumor}, DNA 3D folding \cite{irineo}, image segmentation \cite{image1}, natural language processing \cite{ball}, product-customer segmentation \cite{clauset}, webpage sorting \cite{cyber_comm}, and more. 

The field of community detection has been expanding greatly since the 1980's, with a remarkable diversity of models and algorithms developed in different communities such as machine learning, computer science, network science, social science and statistical physics. These rely on various benchmarks for finding clusters, in particular, cost functions based on cuts or modularities \cite{newman-girvan}. We refer to \cite{newman-book,fortunato,airoldi,social1} for an overview of these developments. 

Nonetheless, various fundamental questions remain unsettled, such as:
\begin{itemize}
\item When are there really communities?  Algorithms may output community structures, but are these meaningful or artefacts?
\item Can we always extract the communities when they are present; fully, partially?
\item What is a good benchmark to measure the performance of algorithms, and how good are the current algorithms?
\end{itemize}
The goal of this monograph is to describe recent developments aimed at answering these questions in the context of block models. Block models are a family of random graphs with planted clusters. The ``mother model'' is the {\it stochastic block model (SBM)}, which has been widely employed as a canonical model for community detection. It is arguably the simplest model of a graph with communities (see definitions in the next section). Since the SBM is a generative model, it benefits from a ground truth for the communities, which allows us to consider the previous questions in a formal context. 
Like any model, it is not necessarily realistic, but it is insightful - judging for example from the powerful algorithms that have emerged from its study. 

In a sense, the SBM plays a similar role to the discrete memoryless channel (DMC) in information theory. While the task of modelling external noise may be more amenable to simplifications than real data sets, the SBM captures some of the key bottleneck phenomena for community detection and admits many possible refinements that improve its fit to real data. Our focus here will be on the fundamental understanding of the ``canonical SBM,'' without diving too much into the refined extensions.  

The SBM is defined as follows. For positive integers $k, n$, a probability vector $p$ of dimension $k$, and a symmetric matrix $W$ of dimension $k \times k$ with entries in $[0,1]$, the model $\sbm(n,p,W)$ defines an $n$-vertex random graph with vertices split in $k$ communities, where each vertex is assigned a community label in $\{1,\dots,k\}$ independently under the community prior $p$, and pairs of vertices with labels $i$ and $j$ connect independently with probability $W_{i,j}$. 

Further generalizations allow for labelled edges and continuous vertex labels, connecting to low-rank approximation models and graphons (using the latter terminology as adapted in the statistics literature). For example, a spiked Wigner model with observation $Y=XX^T + Z$, where $X$ is an unknown vector and $Z$ is Wigner, can be viewed as a labeled graph where edge $(i,j)$'s label is given by $Y_{ij}=X_iX_j + Z_{ij}$. If the $X_i$'s take discrete values, e.g., $\{1,-1\}$, this is closely related to the stochastic block model---see \cite{yash_sbm} for a precise connection. Continuous labels can also model Euclidean connectivity kernels, an important setting for data clustering. In general, models where a collection of variables $\{X_i\}$ have to be recovered from noisy observations $\{Y_{ij}\}$ that are stochastic functions of $X_i,X_j$, or more generally that depend on local interactions of some of the $X_i$'s, can be viewed as inverse problems on graphs or hypergraphs that bear similarities with the basic community detection problems discussed here. This concerns in particular topic modelling, ranking, synchronization problems and other unsupervised learning problems. We refer to Section \ref{others} for further discussion on these.  The specificity of the stochastic block model is that the input variables are discrete.

%A general abstraction of these problems can also be obtained from an information theoretic point of view, with graphical channels \cite{abbetoc}, themselves a special case of conditional random fields \cite{lafferty2001conditional}, which model conditional distributions between a collection of vertex variables $X^V$ and a collection of edge variables $Y^E$ on a hyper-graph $G=(V,E)$, where the conditional probability distributions factors over each edge with a local kernel $Q$:
%\begin{align*}
%P(y^E|x^V) =\prod_{I \in E} Q_I(y_I|x[I]),
%\end{align*}
%where $y_I$ is the realization of $Y$ on the hyperedge $I$ and $x[I]$ is the realization of $X^V$ over the vertices incident to the hyperedge $I$. Our goal in this note is to discuss tools and methods for the SBM that are likely to extend to the analysis of such general models.
%

A first hint at the centrality of the SBM comes from the fact that the model appeared independently in numerous scientific communities. It appeared under the SBM terminology in the context of social networks in the machine learning and statistics literature \cite{holland}, while the model is typically called the planted partition model in theoretical computer science \cite{bui,dyer,boppana}, and the inhomogeneous random graph in the mathematics literature \cite{bollo_inhomo}. The model takes also different interpretations, such as a planted spin-glass model \cite{decelle}, a sparse-graph code \cite{abh,colin1focs} or a low-rank (spiked) random matrix model \cite{mcsherry,Vu-arxiv,yash_sbm} among others.

In addition, the SBM has recently turned into more than a model for community detection. It provides a fertile ground for studying various central questions in machine learning, computer science and statistics: It is rich in phase transitions \cite{decelle,massoulie-STOC,Mossel_SBM2,abh,colin1focs}, allowing us to study the interplay between statistical and computational barriers \cite{chen-xu,colin2nips,banks2,colin3cpam}, as well as the discrepancies between probabilstic and adversarial models \cite{ankur_SBM}, and it serves as a test bed for algorithms, such as SDPs \cite{abh,new-xu,new-xu2,sbm-groth,levina,afonso3,montanari_sen,wein_sdp}, spectral methods \cite{Vu-arxiv,jiaming,massoulie-STOC,redemption,bordenave,prout}, and belief propagation \cite{decelle,colin3}.  

%\section{Inference on graphs}\label{graphical_channel}

\section{Fundamental limits: information and computation}
This monograph focuses on the {\it fundamental limits} of community detection. The term `fundamental limit' is used to emphasize the fact that we seek conditions for recovering the communities that are {\it necessary and sufficient}. In the information-theoretic sense, this means finding conditions under which a given task can or cannot be resolved irrespective of complexity or algorithmic considerations, whereas in the computational sense, this further constrains the algorithms to run in polynomial time in the number of vertices. As we shall see in this monograph, such fundamental limits are often expressed through {\it phase transitions}, which provide sharp transitions in the relevant regimes between phases where the given task can or cannot be resolved.

Fundamental limits have proved to be instrumental in the developments of algorithms. A prominent example is Shannon's coding theorem \cite{shannon48}, which gives a sharp threshold for coding algorithms at the channel capacity, and which has led the development of coding algorithms for more than 60 years (e.g., LDPC, turbo or polar codes) at both the theoretical and practical level \cite{RU01}. 
Similarly, the SAT threshold \cite{AchlioetalNature} has driven the developments of a variety of satisfiability algorithms such as survey propagation \cite{MezParZec_science}.

In the area of clustering and community detection, where establishing rigorous benchmarks is a long standing challenge, the quest of fundamental limits and phase transitions is also impacting the development of algorithms. As discussed in this monograph, this has already lead to developments of algorithms such as sphere-comparisons, linearized belief propagation, nonbacktracking spectral methods. Fundamental limits also shed light on the limitations of the model versus those of the algorithms used; see Section \ref{datab}. However, unlike in the data transmission context of Shannon, information-theoretic limits may not always be efficiently achievable in community detection, with information-computation gaps that may emerge as discussed in Section \ref{info}.

\section{An example on real data}\label{datab}
This monograph focuses on the fundamentals of community detection, but we want to give an application example here. We use the blogosphere data set from the 2004 US political elections \cite{blogs} as an archetype example.

Consider the problem where one is interested in extracting features about a collection of items, in our case $n=1,222$ individuals writing about US politics, observing only some of their interactions. In our example, we have access to which blogs refers to which (via hyperlinks), but nothing else about the content of the blogs. The hope is to extract knowledge about the individual features from these simple interactions.   

To proceed, build a {\it graph of interaction} among the $n$ individuals, connecting two individuals if one refers to the other, ignoring the direction of the hyperlink for simplicity. Assume next that the data set is generated from a stochastic block model; assuming two communities is an educated guess here, but one can also estimate the number of communities (e.g., as in \cite{colin2nips}). The type of algorithms developed in Sections \ref{weak} and \ref{exact} can then be run on this data set, and two assortative communities are obtained. In the paper \cite{blogs}, Adamic and Glance recorded which blogs are right or left leaning, so that we can check how much agreement the algorithms give with the true partition of the blogs. The results give about 95\% agreement on the blogs' political inclinations (which is roughly the state-of-the-art \cite{newman1,score,harrison}).

\begin{figure}[h]
\centering
\begin{subfigure}{.45\textwidth}
  \centering
  \includegraphics[width=.82\linewidth]{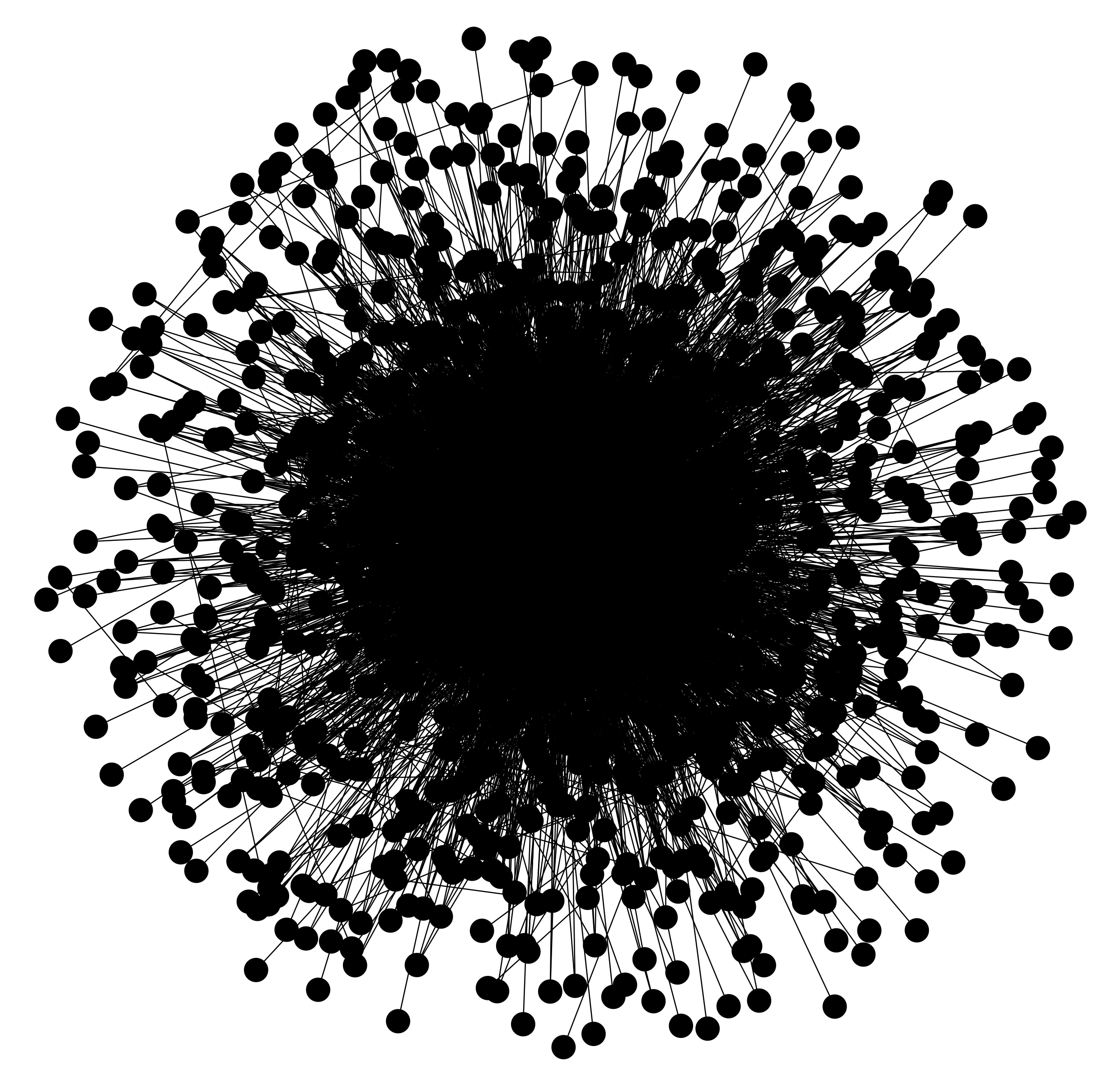}
%  \caption{Random configuration.}
  \label{blog1}
\end{subfigure}%
\begin{subfigure}{.5\textwidth}
  \centering
  \includegraphics[width=.95\linewidth]{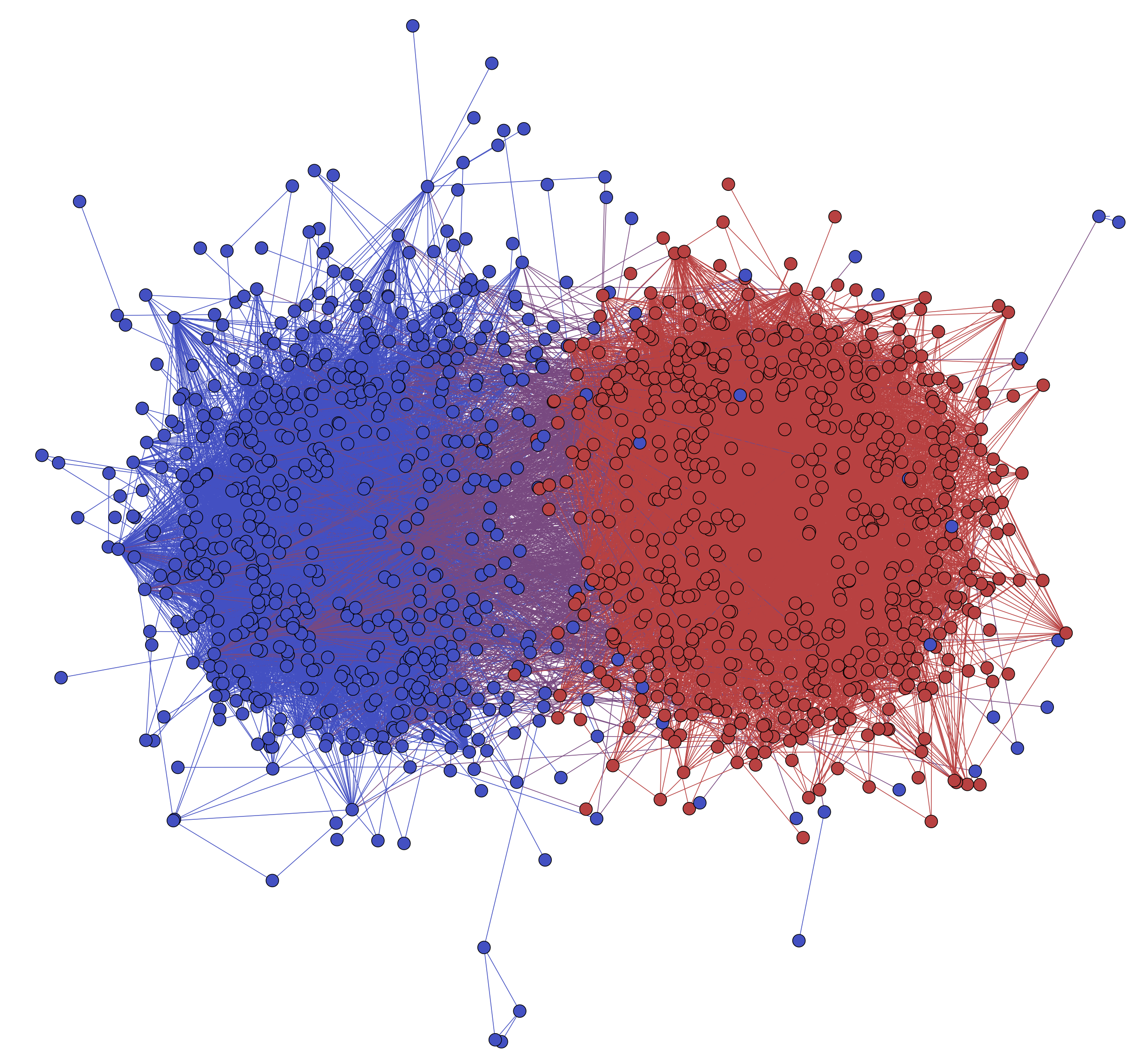}
 % \caption{Clustered configuration.}
  \label{blog2}
\end{subfigure}
\caption{The above graphs represent the real data set of the political blogs from \cite{blogs}. Each vertex represents a blog and each edge represents the fact that one of the blogs refers to the other. The left graph is plotted with a random arrangement of the vertices, and the right graph is the output of the ABP algorithm described in Section \ref{weak}, which gives 95\% accuracy on the reconstruction of the political inclination of the blogs (blue and red colors correspond to left and right leaning blogs).   
}
\label{fig:test}
\end{figure}

Despite the fact that the blog data set is particularly `well behaved'--there are two dominant clusters that are well balanced and well separated--the above approach can be applied to a broad collection of data sets to extract knowledge about the data from graphs of similarities or interactions. In some applications, the graph is obvious (such as in social networks with friendships), while in others, it is engineered from the data set based on metrics of similarity/interactions that need to be chosen properly (e.g, similarity of pixels in image segmentation). The goal is to apply such approaches to problems where the ground truth is unknown, such as to understand biological functionality of protein complexes; to find genetically related sub-populations; to make accurate recommendations; medical diagnosis; image classification; segmentation; page sorting; and more (see references in the introduction).

In such cases where the ground truth is not available, a key question is to understand how reliable the algorithms' outputs may be. 
We now discuss how the results presented in this monograph add to this question. 
%On this matter, the theory discussed in this note gives a new perspective as follows. 
% \cite{newman-girvan,social1,ppi2,marcotte,gene-expr,gene-expr2,genetics,gene-survey,palla,prem,citation,comm-net,comm2,dixon,image1,image2,newman1,fortunato,tumor,Geography,KleinbergGroup,sparse-network1,reco2,comparatif}. 
Following the definitions from Sections \ref{weak} and \ref{exact}, 
the parameters estimated by fitting an SBM on this data set in the 
constant degree regime are
\begin{align}
p_1=0.48, \,\,\, p_2=0.52, \quad 
Q=\begin{pmatrix}
52.06 & 5.16\\
5.16 & 47.43
\end{pmatrix}.
\end{align}
and in the logarithmic degree regime
\begin{align}
p_1=0.48, \,\,\, p_2=0.52, \quad 
Q=\begin{pmatrix}
7.31 & 0.73\\
0.73 & 6.66
\end{pmatrix}.
\end{align}
Following the definitions of Theorem \ref{main} from Section \ref{weak}, we can now compute the SNR for these parameters in the constant-degree regime, obtaining $\lambda_2^2/\lambda_1 \approx 18 $ which is much greater than 1. Thus, under an SBM model, the data is largely in a regime where communities can be detected, i.e., above the weak recovery threshold. Following the definitions of Theorem \ref{exact_thm} from Section \ref{exact}, we can also compute the CH-divergence for these parameters in the logarithmic-degree regime, obtaining $J(p,Q) \approx 2$ which is also greater than 1. Thus, under an SBM and with an asymptotic approximation, the data is in a regime where the graph communities can in fact be recovered entirely, i.e, above the exact recovery threshold. This does not answer whether the SBM is a good or a bad model, but it gives that under this model, the data appears to be in a strong `clusterable regime.'

Note also that such a conclusion may not appear using a specific algorithm, e.g., one that is sensitive to the degree variations and that may split the vertices into high vs.\ low-degree vertices. This prompted for example the development of degree-corrected SBMs in \cite{newman2}, as the algorithm used in \cite{newman2} for the blog data set with the fitting of an SBM failed for such reasons. However, how do we know whether the failure is due to the model or the algorithm? By establishing the fundamental limits on the SBM, we will find algorithms that are `maximally' robust by succeeding in the most challenging regimes, i.e., down to the fundamental limits, which achieve in particular the positive accuracy for the blog data set described in Figure \ref{fig:test}. We also refer to Section \ref{robust} for discussions on the robustness of algorithms to degree variations.

\section{Historical overview of the recent developments}\label{history}
This section provides a brief historical overview of the recent developments discussed in this monograph. The resurgent interest in the SBM and its `modern study' have been initiated in part due to the paper of Decelle, Krzakala, Moore and Zdeborov\'a \cite{decelle}, which conjectured\footnote{The conjecture of the Kesten-Stigum threshold in \cite{decelle} was formulated with what we call in this note the max-detection criteria, asking for an algorithm to output a reconstruction of the communities that strictly improves on the trivial performance achieved by putting all the vertices in the largest community. This conjecture is formally incorrect for general SBMs, see \cite{colin3cpam} for a counter-example, as the notion of max-detection is too strong in some cases. The conjecture is believed to hold for symmetric SBMs, as re-stated in \cite{Mossel_SBM1}, but it requires a different notion of detection to hold for general SBMs; see definitions from \cite{colin3cpam} discussed in Section \ref{weak}.} phase transition phenomena for the weak recovery (a.k.a.\ detection) problem at the Kesten-Stigum threshold and the information-computation gap at 4 symmetric communities in the symmetric case. These conjectures are backed in \cite{decelle} with insights from statistical physics, based on the cavity method (belief propagation), and provide a detailed picture of the weak recovery problem, both for the algorithmic and information-theoretic behavior. With such insights, a new research program started driven by the phase transition phenomena. %and the development of algorithms that can achieve these. 

One of the first papers that obtained a non-trivial algorithmic result for the weak recovery problem is \cite{coja-sbm} from 2010, which appeared before the conjecture (and does not achieve the threshold by a logarithmic degree factor). The first paper that made progress on the conjecture is \cite{Mossel_SBM1} from 2012, which proved the impossibility part of the conjecture for two symmetric communities, introducing various key concepts in the analysis of block models. In 2013, \cite{mossel2} also obtained a result on the partial recovery of the communities, expressing the optimal fraction of mislabelled vertices when the signal-to-noise ratio is large enough in terms of the broadcasting problem on trees \cite{ks1,evans}.

The positive part of the conjecture for efficient algorithm and two communities was first proved in 2014 with \cite{massoulie-STOC} and \cite{Mossel_SBM2}, using respectively a spectral method from the matrix of self-avoiding walks and weighted non-backtracking walks between vertices. 

In 2014, \cite{abh_arxiv,abh} and \cite{mossel-consist} found that the exact recovery problem for two symmetric communities has also a phase transition, in the logarithmic rather than constant degree regime, shown to be also efficiently achievable. This relates to a large body of work from the first decades of research on the SBM \cite{bui,dyer,boppana,snij,condon,mcsherry,bickel,choi,Vu-arxiv,chen-xu}, driven by the exact or almost exact recovery problems without sharp thresholds.

In 2015, the phase transition for exact recovery was obtained for the general SBM \cite{colin1focs,colin2nips}, and shown to be efficiently achievable irrespective of the number of communities. For the weak recovery problem, \cite{bordenave} showed that the Kesten-Stigum threshold can be achieved with a spectral method based on the nonbacktracking (edge) operator in a fairly general setting (covering SBMs that are not necessarily symmetric), but felt short of settling the conjecture for more than two communities in the symmetric case due to technical reasons. The approach of \cite{bordenave} is based on the `spectral redemption' conjecture made in 2013 in \cite{redemption}, which introduced the use of the nonbacktracking operator as a linearization of belief propagation. This is one of the most elegant approaches to the weak recovery problem, except perhaps for the fact that the matrix is not symmetric (note that the first proof of \cite{massoulie-STOC} does provide a solution with a symmetric matrix via the count of self-avoiding walks, albeit less direct to construct). The general conjecture for arbitrary many symmetric or asymmetric communities is settled later in 2015 with \cite{colin3,colin3cpam}, relying on a higher-order nonbacktracking matrix and a message passing implementation. It was further shown in \cite{colin3,colin3cpam} that it is possible to cross information-theoretically the Kesten-Stigum threshold in the symmetric case at 4 communities, settling both positive parts of the conjectures from \cite{decelle}. Crossing at 5 rather than 4 communities is also obtained in \cite{banks,banks2}, which further obtains the scaling of the information-theoretic threshold for a growing number of communities. 

In 2016, a tight expression was obtained for partial recovery with two communities in the regime of finite SNR with diverging degrees in \cite{yash_sbm} and \cite{mossel-xu} for a different distortion measure. This also gives the threshold for weak recovery in the regime where the SNR is finite while the degrees are diverging.

Other major lines of work on the SBM have been concerned with the performance of SDPs, with a precise picture obtained in \cite{sbm-groth,montanari_sen,adel_sbm} for the weak recovery problem and in \cite{abh,new-xu,levina,afonso_single,afonso3,wein_sdp} for the (almost) exact recovery problem; as well as with spectral methods on classical operators \cite{mcsherry,coja-sbm,new-vu,jiaming,Vu-arxiv,prout,prout2}. A detailed picture has also been developed for the problem of a single planted community in \cite{am_1comm,hajek_1comm,1comm_info,lelarge_ks}. Recently, attention has been paid to graphs that have a larger number of short loops \cite{our_grid,abishek2,arya,powering1}. There is a much broader list of works on the SBMs that is not covered in this monograph, especially before the `recent developments' discussed above but also after. It is particularly challenging to track the vast literature on this subject as it is split between different communities of statistics, machine learning, mathematics, computer science, information theory, social sciences and statistical physics. This monograph mainly covers developments until 2016, with some references from 2017 There a few additional surveys available; community detection and statistical network models are discussed in  \cite{newman-book,fortunato,airoldi}, and C.~Moore has a recent overview paper \cite{moore_review} that focuses on the weak recovery problem with emphasis on the cavity method.

In the table below, we summarize the main thresholds proved for weak and exact recovery, covered in several chapters of this monograph:

{\footnotesize
\begin{center}
  \begin{tabular}{| l |  c | c  | }
   \hline
            & Exact recovery & Weak recovery (detection)  \\
           & (logarithmic degrees) & (constant degrees)   \\
   \hline        
    2-SSBM & $| \sqrt{a}-\sqrt{b}|  >\sqrt{2}$  \cite{abh_arxiv,mossel-consist} & $(a-b)^2 >2(a+b)$   \cite{massoulie-STOC,Mossel_SBM2}  \\
    %& (Massouli\'e '13, Mossel et al.\ '13) & (This paper)  \\
   % &\cite{massoulie-STOC,Mossel_SBM2} & \cite{abh,mossel-consist} \\
    \hline
               General SBM  & $\min\limits_{i<j} D_+((PQ)_i , (PQ)_j) > 1$  \cite{colin1focs}  &  $\lambda_2^2(PQ) > \lambda_1(PQ)$ \cite{bordenave,colin3}   \\
    %& (Massouli\'e '13, Mossel et al.\ '13) & (This paper)  \\
%    &\cite{bordenave,colin3} &  \cite{colin1focs} \\
    \hline
  \end{tabular}
\end{center}
}

\section{Outline}
In the next section, we formally define the SBM and various recovery requirements for community detection, namely exact, weak, and partial recovery. 
We then start with a quick overview of the key approaches for these recovery requirements in Section \ref{tackle}, introducing the key new concepts obtained in the recent developments. 
We then treat each of these three recovery requirements separately for the two community SBM in Sections \ref{exact}, \ref{weak} and \ref{partial} respectively, discussing both fundamental limits and efficient algorithms. 
We give complete (and revised) proofs for exact recovery and partial proofs for weak and partial recovery. 
We then move to the results for the general SBM in Section \ref{general}.
In Section \ref{others} we discuss other block models, such as geometric block models, and in Section \ref{open} we give concluding remarks and open problems.

\section{Notations}
We use the standard little-o and big-o notations.  
Recall that $a_n=\Omega(b_n)$ means that $b_n = O(a_n)$,
and $a_n=\omega(b_n)$ means that $b_n = o(a_n)$.
In particular, $a_n = o(1)$ means that $a_n$ is vanishing, $a_n = \Omega(1)$ means that $a_n$ is non-vanishing, and $a_n=\omega(1)$ means that $a_n$ is diverging. 
We use 
$a_n \lesssim b_n$ when $a_n=\Omega(b_n)$; 
$a_n \ll b_n$ when $a_n= o(b_n)$ (and $a_n \gg b_n$ when $b_n= o(a_n)$);
$a_n = \Theta(b_n)$, or equivalently $a_n \asymp b_n$, when we simultaneously have $a_n=\Omega(b_n)$ and $a_n = O(b_n)$;
$a_n \sim b_n$ when $a_n=b_n(1 + o(1))$.

We say that an event $E_n$ takes place with high probability if its probability tends to 1 as $n$ diverges, i.e., $\pp\{E_n\}=1-o(1)$.
We also use a.a.e.\ and a.a.s.\ for asymptotically almost everywhere and asymptotically almost surely (respectively). 

We usually use superscripts to specify the dimensions of vectors; in particular, $1^n$ is the all-one vector of dimension $n$, $0^n$ the all-zero vector of dimension $n$, and $x^n=(x_1,\dots,x_n)$.

\chapter{The stochastic block model}\label{model}
The history of the SBM is long, and we omit a comprehensive treatment here. As mentioned earlier, the model appeared independently in multiple scientific communities: the term SBM, which seems to have dominated in the recent years, comes from the machine learning and statistics literature \cite{holland}, while the model is typically called the planted partition model in theoretical computer science \cite{bui,dyer,boppana}, and the inhomogeneous random graphs model in the mathematics literature \cite{bollo_inhomo}. 

\section{The general SBM}
%The SBM is a random graph with planted clusters defined as follows. 
\begin{definition}
Let $n$ be a positive integer (the number of vertices), $k$ be a positive integer (the number of communities), $p=(p_1,\dots,p_k)$ be a probability vector on $[k]:=\{1,\dots,k\}$ (the prior on the $k$ communities) and $W$ be a $k\times k$ symmetric matrix with  entries in $[0,1]$ (the connectivity probabilities). The pair $(X,G)$ is drawn under $\mathrm{SBM}(n,p,W)$ if $X$ is an $n$-dimensional random vector with i.i.d.\ components distributed under $p$, and $G$ is an $n$-vertex simple graph where vertices $i$ and $j$ are connected with probability $W_{X_i,X_j}$, independently of other pairs of vertices. We also define the community sets by $\Omega_i=\Omega_i(X) : = \{v \in [n] : X_v = i\}, i \in [k].$
\end{definition}
\noindent
Thus the distribution of $(X,G)$ where $G=([n],E(G))$ is defined as follows; for $x \in [k]^n$ and $y \in \{0,1\}^{{n \choose 2}}$,  
\begin{align}
\pp\{X=x\} &:= \prod_{u=1}^n p_{x_u}=\prod_{i=1}^k p_{i}^{|\Omega_i(x)|},\\
\pp\{E(G)=y | X=x\} &:= \prod_{1 \le u < v \le n} W_{x_u,x_v}^{y_{uv}} (1-W_{x_u,x_v})^{1-y_{uv}} \\
& = \prod_{1 \leq i \leq j \leq k} W_{i,j}^{N_{ij}(x,y)} (1-W_{i,j})^{N^c_{ij}(x,y)} ,
\end{align}
where
\begin{align}
N_{ij}(x,y)&:=\sum_{u<v \atop x_u=i,x_v=j} \1(y_{uv}=1),\\
N^c_{ij}(x,y)&:=\sum_{u<v \atop x_u=i,x_v=j} \1(y_{uv}=0)=|\Omega_i(x)| |\Omega_j(x)|-N_{ij}(x,y), \,\,\, i\neq j \\
N^c_{ii}(x,y)&:=\sum_{u<v \atop x_u=i,x_v=i} \1(y_{uv}=0)=|\Omega_i(x)| (|\Omega_i(x)|-1)/2-N_{ii}(x,y), 
\end{align}
are the number of edges and non-edges between pairs of communities. We may also talk about $G$ drawn under $\mathrm{SBM}(n,p,W)$ without specifying the underlying community labels $X$.

\begin{remark}
Except for Section \ref{open}, we assume that $p$ does not scale with $n$, whereas $W$ typically does. As a consequence, the number of communities does not scale with $n$ and the communities have linear size. Nonetheless, various results discussed in this monograph should extend (by inspection) to cases where $k$ is growing slowly enough.
\end{remark}
\begin{remark}
Note that by the law of large numbers, almost surely, 
$$ \frac{1}{n} | \Omega_i | \to p_i .$$
Alternative definitions of the SBM require $X$ to be drawn uniformly at random with the constraint that $\frac{1}{n} |\{v \in [n] : X_v = i \}| = p_i +o(1)$, or $\frac{1}{n} |\{v \in [n] : X_v = i \}| = p_i$ for consistent values of $n$ and $p$ (e.g., $n/2$ being an integer for two symmetric communities). 
%In the latter case, the joint distribution of $(X,G)$ is given by 
%\begin{align}
%\pp\{E(G)=y , X=x\} & = \prod_{1 \leq i \leq j \leq k} W_{i,j}^{N_{ij}(x,y)} (1-W_{i,j})^{n^2p_ip_j-N_{ij}(x,y)}  \prod_{i=1}^k \1(|\Omega_i(x)|=np_i).
%\end{align}
For the purpose of this paper, these definitions are essentially equivalent, and we may switch between the models to simplify some proofs. 
\end{remark}
Note also that if all entries of $W$ are the same, then the SBM collapses to the Erd\H{o}s-R\'enyi random graph, and no meaningful reconstruction of the communities is possible.

\section{The symmetric SBM}
The SBM is called symmetric if $p$ is uniform and if $W$ takes the same value on the diagonal and the same value outside the diagonal. 
\begin{definition}
$(X,G)$ is drawn under $\mathrm{SSBM}(n,k,\pin,\pou)$ if $W$ takes value $\pin$ on the diagonal and $\pou$ off the diagonal, and if the community prior is $p=\{1/k\}^k$ in the Bernoulli model, and $X$ is drawn uniformly at random with the constraints $|\{v \in [n] : X_v = i \}| = n/k$, $n$ a multiple of $k$, in the uniform or strictly balanced model.
\end{definition}

\section{Recovery requirements}\label{alldef}
The goal of community detection is to recover the labels $X$ by observing $G$, up to some level of accuracy. We next define the notions of agreement. 
\begin{definition}[Agreement and normalized agreement]
The agreement between two community vectors $x,y \in [k]^n$ is obtained by maximizing the common components between $x$ and any relabelling of $y$, i.e., \begin{align}
A(x,y) &= \max_{\pi \in S_k} \frac{1}{n} \sum_{i=1}^n \1(x_i=\pi(y_i)),\\
\tilde{A}(x,y) &= \max_{\pi \in S_k} \frac{1}{k} \sum_{i=1}^k \frac{\sum_{u \in [n]}\1(x_u=\pi(y_u), x_u=i)}{\sum_{u \in [n]}\1(x_u=i)},
\end{align}
\end{definition}
Note that the relabelling permutation is used to handle symmetric communities such as in SSBM, as it is impossible to recover the actual labels in this case, but we may still hope to recover the {\it partition}. In fact, one can alternatively work with the community partition $\Omega=\Omega(X)$, defined earlier as the unordered collection of the $k$ disjoint unordered subsets $\Omega_1,\dots,\Omega_k$ covering $[n]$ with $\Omega_i=\{u \in [n]: X_u =i \}$. It is however often convenient to work with vertex labels. Further, upon solving the problem of finding the partition, the problem of assigning the labels is often a much simpler task. It cannot be resolved if symmetry makes the community label non identifiable, such as for SSBM, and it is trivial otherwise by using the community sizes and clusters/cuts densities. 

For $(X,G)\sim \mathrm{SBM}(n,p,W)$ one can always attempt to reconstruct $X$ without even taking into account $G$, simply drawing each component of $\hat{X}$ i.i.d.\ under $p$. Then the agreement satisfies almost surely 
\begin{align}
A(X,\hat{X}) \to \| p\|_2^2,
\end{align}
and $\| p\|_2^2=1/k$ in the case of $p$ uniform. Thus an agreement becomes interesting only when it is above this value. 

One can alternatively define a notion of component-wise agreement. Define the overlap between two random variables $X,Y$ on $[k]$ as 
\begin{align}
O(X,Y)= \sum_{z \in [k]} (\pp\{ X=z,Y=z\} - \pp\{ X=z\} \pp\{Y=z\} )  
\end{align}
and $O^*(X,Y)= \max_{\pi \in S_k}  O(X,\pi(Y))$.
In this case, for $X,\hat{X}$ i.i.d.\ under $p$, we have $O^*(X,\hat{X})=0$.

When discussing impossibility for weak recovery in the SSBM (Section \ref{bot-conv}), we also use an alternative definition for which the conditional mutual information between an arbitrary pair of vertices $u \ne v$ vanishes, i.e., $$I(X_u;X_v|G) \to 0$$ as $n \to \infty$. Note that if this mutual information between two arbitrary vertices $u$ and $v$ is vanishing, i.e., the two vertex labels are asymptotically independent conditioned on the graph, then it is not possible to obtain a reconstruction of {\it all} vertices that solves weak recovery. 

All recovery requirements in this note are going to be asymptotic, taking place with high probability as $n$ tends to infinity. We also assume in the following sections---except for Section \ref{learn}---that the parameters of the SBM are known when designing the algorithms. 
\begin{definition}\label{all-defs}
Let $(X,G)\sim \mathrm{SBM}(n,p,W)$. The following recovery requirements are solved if there exists an algorithm that takes $G$ as an input and outputs $\hat{X}=\hat{X}(G)$ such that 
\begin{itemize}
\item {\bf Exact recovery:} $\pp\{A(X,\hat{X})=1\} = 1-o(1)$,
\item {\bf Almost exact recovery:} $\pp\{A(X,\hat{X})=1-o(1)\} = 1-o(1)$,
\item {\bf Partial recovery:} $\pp\{\tilde{A}(X,\hat{X}) \geq \alpha \} = 1-o(1)$, $\alpha \in (1/k,1)$,
\item {\bf Weak recovery:} $\pp\{\tilde{A}(X,\hat{X}) \geq 1/k + \Omega(1)\} = 1-o(1)$.
\end{itemize}
\end{definition}
In other words, exact recovery requires the entire partition to be correctly recovered, almost exact recovery allows for a vanishing fraction of misclassified vertices, partial recovery allows for a constant fraction of misclassified vertices and weak recovery allows for a non-trivial fraction of misclassified vertices. We call $\alpha$ the agreement or accuracy of the algorithm. Note that partial and weak recovery are defined by means of the normalized agreement $\tilde{A}$ rather the agreement $A$. The reason for this is discussed in details below and matters for asymmetric SBMs; for the case of symmetric SBMs, $A$ can be used for all four definitions.

Different terminologies are sometimes used in the literature, with following equivalences:
\begin{itemize}
\item exact recovery $\iff$ strong consistency
\item almost exact recovery $\iff$ weak consistency
\end{itemize}
Sometimes `exact recovery' is also called just `recovery' and `almost exact recovery' is called `strong recovery.'

As mentioned above, values of $\alpha$ that are too small may not be interesting or realizable. In the symmetric SBM with $k$ communities, an algorithm that ignores the graph and simply draws $\hat{X}$ i.i.d.\ under $p$ achieves an accuracy of $1/k$. Thus the problem becomes interesting when $\alpha > 1/k$, leading to the following definition. 
\begin{definition}\label{detect-symm}
Weak recovery (or detection) is solvable in $\mathrm{SSBM}(n,k,$ $\pin,\pou)$ if for $(X,G)\sim \mathrm{SSBM}(n,k,\pin,\pou)$, there exists $\e>0$ and an algorithm that takes $G$ as an input and outputs $\hat{X}$ such that $\pp\{A(X,\hat{X})\geq 1/k+\e\} = 1-o(1)$. 
\end{definition}

Equivalently, $\pp\{O^*(X_V,\hat{X}_V) \geq \e \}=1-o(1)$ where $V$ is uniformly drawn in $[n]$.
Determining the counterpart of weak recovery in the general SBM requires some discussion. Consider an SBM with two communities of relative sizes $(0.8,0.2)$. A random guess under this prior gives an agreement of $0.8^2 + 0.2^2=0.68$, however an algorithm that simply puts every vertex in the first community achieves an agreement of $0.8$. In \cite{decelle}, the latter agreement is used as the one to improve upon in order to detect communities, leading to the following definition:
\begin{definition}
Max-detection is solvable in $\mathrm{SBM}(n,p,W)$ if for $(X,G)$ $\sim \mathrm{SBM}(n,p,W)$, there exists $\e>0$ and an algorithm that takes $G$ as an input and outputs $\hat{X}$ such that $\pp\{A(X,\hat{X})\geq \max_{i \in [k]} p_i +\e\} = 1-o(1)$. 
\end{definition}
As shown in \cite{colin3cpam}, the previous definition is not the right definition to capture the Kesten-Stigum threshold in the general case. In other words, the conjecture that max-detection is always possible above the Kesten-Stigum threshold is not accurate in general SBMs.
%We provide next a different definition of detection for the general SBM, which differs from the above to cope with the following point. 
Back to our example with communities of relative sizes $(0.8,0.2)$, an algorithm that could find a set containing $2/3$ of the vertices from the large community and $1/3$ of the vertices from the small community would not satisfy the above above weak recovery criteria, while the algorithm produces nontrivial amounts of evidence on what communities the vertices are in. 
To be more specific, consider a two community SBM where each vertex is in community 1 with probability $0.99$, each pair of vertices in community 1 have an edge between them with probability $2/n$, while vertices in community 2 never have edges. Regardless of what edges a vertex has it is more likely to be in community 1 than community 2, so weak recovery according to the above definition is not impossible, but one can still divide the vertices into those with degree 0 and those with positive degree to obtain a non-trivial detection---see \cite{colin3cpam} for a formal counter-example. 

Using the normalized agreement fixes this issue. Weak recovery can then be defined as obtaining with high probability a weighted agreement of $$\tilde{A}(X,\hat{X}(G)) = 1/k + \Omega_n(1),$$ and this applies to the general SBM. Another definition of weak recovery that seems easier to manipulate and that implies the previous one is as follows; note that this definition requires a single partition even for the general SBM. 

\begin{definition}\label{detect-gen}
Weak recovery (or detection) is solvable in $\mathrm{SBM}(n,p,W)$ if for $(X,G)\sim \mathrm{SBM}(n,p,W)$, there exists $\e>0$, $i,j \in [k]$ and an algorithm that takes $G$ as an input and outputs a partition of $[n]$ into two sets $(S,S^c)$ such that 
$$\pp\{ |\Omega_i\cap S|/|\Omega_i|-|\Omega_j\cap S|/|\Omega_j| \geq \epsilon \} = 1-o(1),$$
where we recall that $\Omega_i=\{u \in [n] : X_u=i\}$.
\end{definition}

In other words, an algorithm solves weak recovery if it divides the graph's vertices into two sets such that vertices from two different communities have different probabilities of being assigned to one of the sets. With this definition, putting all vertices in one community does not detect, since $|\Omega_i\cap S|/|\Omega_i|=1$ for all $i \in [k]$. Further, in the symmetric SBM, this definition implies Definition \ref{detect-symm} provided that we fix the output:

\begin{lemma}
If an algorithm solves weak recovery in the sense of Definition \ref{detect-gen} for a symmetric SBM, then it solves max-detection (or detection according to Decelle et al.\ \cite{decelle}), provided that we consider it as returning $k-2$ empty sets in addition to its actual output.
\end{lemma}
See \cite{colin3nips} for the proof. 
%\begin{proof}
%Let $(X,G) \sim \sbm(n,p,W)$ and $\hat{X}$ return $S$ and $S'$. There exists $\epsilon>0$ such that with high probability there exist $i$ and $j$ such that $|\Omega_i\cap S|/|\Omega_i|-|\Omega_j\cap S|/|\Omega_j| \geq \epsilon$. So, if we map $S$ to community $i$ and $S'$ to community $j$, the algorithm classifies at least 
%\[|\Omega_i\cap S|/n+|\Omega_j\cap S'|/n=|\Omega_j|/n+|\Omega_i\cap S|/n-|\Omega_j\cap S|/n\ge 1/k+\epsilon/k-o(1)\]
% of the vertices correctly with high probability.
%\end{proof}
%\Enote{In case people get disturbed by the 'empty sets', maybe mention that one can also pick random subsets appropriately}
%\Cnote{Why is that better? The Decelle definition implicitly uses putting all vertices in one community as a basepoint, so it does not seem like there is a problem with empty sets.}
The above extends to other weakly symmetric SBMs, i.e., those that have constant expected degrees, but not all. 
%However, as mentioned in the example above ,there are cases of asymmetrical SBMs of which the first detection may not be possible while the above still is.  
%Since we will focus on detection for weakly symmetric SBMs in the rest of the is note, we can now talk of a single notion of detection for and use Definition \ref{detect-gen} for convenience. 

Finally, note that our notion of weak recovery requires us to separate at least two communities $i,j \in [k]$. One may ask for a definition where two specific communities need to be separated:
\begin{definition}\label{detect-gen2}
Separation of communities $i$ and $j$, with $i,j \in [k]$, is solvable in $\mathrm{SBM}(n,p,W)$ if for $(X,G)\sim \mathrm{SBM}(n,p,W)$, there exists $\e>0$ and an algorithm that takes $G$ as an input and outputs a partition of $[n]$ into two sets $(S,S^c)$ such that 
$$\pp\{ |\Omega_i\cap S|/|\Omega_i|-|\Omega_j\cap S|/|\Omega_j|\geq \epsilon \} = 1-o(1).$$
\end{definition}
%The following results from \cite{colin3} for the separation property:
%There exists an algorithm $A(G,(\lambda'_1,\lambda'_2,...,\lambda'_{h'}))$ that outputs a partition of $G$'s vertices into two sets and runs in quasilinear time such that the following holds. Let $p$ and $W$ be parameters for $SBM(n,p,W)$ such that $p$ is constant and $W$ is a function of $n$ times a constant. Also, let $1\le i,j\le k$ such that $M$ has an eigenvector $w$ with corresponding eigenvalue greater than $\sqrt{\lambda_i}$ such that $w_i\ne w_j$. There exists $\epsilon>0$ such that for any $(\lambda'_1,...,\lambda'_{h'})$ with $|\lambda_m-\lambda'_m|<\epsilon$ for all $m$ such that $|\lambda_m|<\sqrt{\lambda_1}$, when $A(G,(\lambda'_1,\lambda'_2,...,\lambda'_{h'}))$ is run on $G$ drawn from $SBM(n,p,W)$, the expected difference between the fraction of vertices from community $i$ that $A$ assigns to one of the sets it outputs and the fraction of vertices from community $j$ that $A$ assigns to that set is at least $\epsilon$.

There are at least two additional questions that are natural to ask about SBMs, both can be asked for efficient or information-theoretic algorithms:
\begin{itemize}
\item {\bf Distinguishability (or testing):} Consider an hypothesis test where a random graph $G$ is drawn with probability $1/2$ from an SBM model (with same expected degree in each community) and with probability $1/2$ from an Erd\H{o}s-R\'enyi model with matching expected degree. Is is possible to decide with asymptotic probability $1- o(1)$ from which ensemble the graph is drawn? In particular, this is not possible if the total variation distance between the two ensembles is vanishing. This problem is sometimes called `detection', which partly explains why we prefer to use `weak recovery' in lieu of `detection' for the reconstruction problem discussed earlier. Distinguishability is further discussed in Section \ref{info}. 
\item {\bf Model learnability or parameter estimation:} Assume that $G$ is drawn from an SBM ensemble, is it possible to obtain a consistent estimator for the parameters? E.g., can we estimate $k,p,Q$ from a graph drawn from $\sbm(n,p,Q/n)$? This is further discussed in Section \ref{learn}.
\end{itemize}
The obvious implications are: exact recovery $\Rightarrow$ almost exact recovery $\Rightarrow$ partial recovery $\Rightarrow$ weak recovery $\Rightarrow$ distinguishability, and almost exact recovery $\Rightarrow$ learnability. Moreover, for symmetric SBMs with two symmetric communities: learnability $\Leftrightarrow$ weak recovery $\Leftrightarrow$ distinguishability, but these are broken for general SBMs; see Section \ref{learn}.

\section{SBM regimes and topology}\label{topo}
Before discussing when the various recovery requirements can be solved or not in SBMs, it is important to recall a few topological properties of the SBM graph. 
%\section{Logarithmic degree regime}
%\section{Constant degree regime}

When all the entries of $W$ are the same and equal to $w$, the SBM collapses to the Erd\H{o}s-R\'enyi model $G(n,w)$ where each edge is drawn independently with probability $w$. Let us recall a few basic results for this model derived mainly from \cite{ER-seminal2}:
\begin{itemize}
\item 
 $G(n,c\log(n) /n)$ is connected with high probability if and only if $c>1$,
\item 
$G(n,c /n)$ has a giant component (i.e., a component of size linear in $n$) if and only if $c>1$.
\item For $\delta < 1/2$, the neighborhood at depth $r=\delta \log_c n$ of a vertex $v$ in $G(n,c /n)$, i.e., $B(v,r)=\{u \in [n] : d(u,v) \leq r\}$ where $d(u,v)$ is the length of the shortest path connecting $u$ and $v$, tends in total variation to a Galton-Watson branching process of offspring distribution $\mathrm{Poisson}(c)$.
%\item The number of $m$-cycles in $G(n,c /n)$ tends in distribution to $\mathrm{Poisson}\left(\frac{c^m}{2m}\right)$.
\end{itemize}

For SSBM$(n,k,\pin,\pou)$, these results hold by essentially replacing $c$ with the average degree.
\begin{itemize}
\item For $a,b>0$, SSBM$(n,k,a \log n/n,b \log n/n)$ is connected with high probability if and only if $\frac{a+(k-1)b}{k}>1$ (if $a$ or $b$ is equal to 0, the graph is of course not connected). 
\item SSBM$(n,k,a/n,b/n)$ has a giant component (i.e., a component of size linear in $n$) if and only if $d:=\frac{a+(k-1)b}{k}>1$,
\item For $\delta < 1/2$, the neighborhood at depth $r=\delta \log_d n$ of a vertex $v$ in SSBM$(n,k,a/n,b/n)$ tends in total variation to a Galton-Watson branching process of offspring distribution $\mathrm{Poisson}(d)$ where $d$ is as above.
%\item The number of $m$-cycles in SSBM$(n,k,a/n,b/n)$ tends in distribution to $\mathrm{Poisson}\left(\frac{d^m}{2m}\right)$.
\end{itemize}

Similar results hold for the general SBM, at least for the case of a constant excepted degree. For connectivity, one has that $\sbm(n,p,Q \log n /n)$ is connected with high probability if 
\begin{align}
\min_{i \in [k]} \| (\diag(p)Q)_i\|_1 >1
\end{align}
and is not connected with high probability if $\min_{i \in [k]} \| (\diag(p)Q)_i\|_1 <1$, where $(\diag(p)Q)_i$ is the $i$-th column of $\diag(p)Q$.

These results are important to us as they already point regimes where exact or weak recovery are not possible. Namely, if the SBM graph is not connected, exact recovery is not possible (since there is no hope to label disconnected components with higher chance than $1/2$), hence exact recovery can take place only if the SBM parameters are in the logarithmic degree regime. In other words, exact recovery in SSBM$(n,k,a \log n/n,b \log n/n)$ is not solvable if $\frac{a+(k-1)b}{k}<1$. This is however unlikely to provide a tight condition, i.e., exact recovery is not equivalent to connectivity, and the next section will precisely investigate how much more than $\frac{a+(k-1)b}{k}>1$ is needed to obtain exact recovery. Similarly, it is not hard to see that weak recovery is not solvable if the graph does not have a giant component, i.e., weak recovery is not solvable in SSBM$(n,k,a/n,b/n)$ if $\frac{a+(k-1)b}{k}<1$, and we will see in Section \ref{weak} how much more is needed to go from the giant component to weak recovery. 

\section{Learning the model}\label{learn}

In this section we overview the results on estimating the SBM parameters by observing a realization of the graph. The estimation problem tends to be `easier' than the recovery problems, although some problems remain for the general sparse SBM, see Section \ref{const-est}. We consider first the case where degrees are diverging, where estimation can be obtained as a side result of universal almost exact recovery, and the case of constant degrees, where estimation can be performed without being able to recover the clusters but only above the weak recovery threshold.   
\subsection{Diverging degree regime}\label{learn1}
For diverging degrees, one can estimate the parameters by first solving almost exact recovery without knowing the parameters, and proceeding then to estimate the parameters by simply computing the clusters' cuts and volumes. This requires solving a potentially harder problem, but it turns out to be achievable:  
\begin{theorem}\cite{colin2nips}\label{learn_dense}
Given $\delta>0$ and for any $k\in \mathbb{Z}$, $p\in (0,1)^k$ with $\sum p_i=1$ and $0<\delta\le \min p_i$, and any symmetric matrix $Q$ with no two rows equal such that every entry in $Q^k$ is strictly positive (in other words, $Q$ such that there is a nonzero probability of a path between vertices in any two communities in a graph drawn from SBM$(n,p,Q/n)$), there exist $\epsilon(c)=O(1/\log(c))$ such that for all sufficiently large $\alpha$, the Agnostic-sphere-comparison algorithm detects communities in graphs drawn from SBM$(n,p,\alpha Q/n)$ with accuracy at least $1-e^{-\Omega(\alpha)}$ in $O_n(n^{1+\epsilon(\alpha)})$ time.
\end{theorem}
The algorithm used in this theorem (agnostic-sphere-comparison) is discussed in section \ref{almost} in the context of two communities and is based on comparing neighborhoods of vertices.
Note that assumption that $\delta$ is given in this theorem can be removed when $\alpha=\omega(1)$, i.e., when the degrees diverge. We then obtain: 
\begin{corollary}\cite{colin2nips}
The number of communities $k$, the community prior $p$ and the connectivity matrix $Q$ can be consistently estimated in quasi-linear time in SBM$(n,p,\omega(1) Q/n)$.
\end{corollary}

Note that for symmetric SBMs, certain algorithms such as SDPs or spectral algorithms discussed in Section \ref{tackle} can also be used to recover the communities without knowledge of the parameters, and thus to learn the parameters in the symmetric case. 
A different line of work has also studied the problem of estimating `graphons' \cite{choi,sbm_graphon,wolfe2} via block models, assuming regularity conditions on the graphon, such as piecewise Lipschitz, to obtain estimation guarantees. 
In addition, \cite{borgs_nips} considers private graphon estimation in the logarithmic degree regime, and obtains a non-efficient procedure to estimate `graphons' in an appropriate version of the $L_2$ norm. More recently, \cite{borgs_new} extends the type of results from \cite{wolfe3} to a much more general family of `graphons' and to sparser regimes (though still with diverging degrees) with efficient methods (based on degrees) and non-efficient methods (based on least square and least cut norm).

\subsection{Constant degree regime}\label{const-est}
In the case of the constant degree regime, it is not possible to recover the clusters fully, and thus estimation has to be done differently. The first paper that shows how to estimate tightly the parameter in this regime is \cite{Mossel_SBM1}, which is based on approximating cycle counts by nonbacktracking walks. An alternative method based on expectation-maximization using the Bethe free energy is also proposed in \cite{decelle} (without a rigorous analysis). 
\begin{theorem}\label{mossel_estim}\cite{Mossel_SBM1}
Let $G \sim \ssbm(n,2,a/n,b/n)$ such that $(a-b)^2> 2(a+b)$, and let $C_m$ be the number of $m$-cycles in $G$, $\hat{d}_n=2 |E(G)|/n$ be the average degree in $G$ and $\hat{f}_n= (2 m_n C_{m_n}- \hat{d}_n^{m_n})^{1/m_n}$ where $m_n=\lfloor \log^{1/4}(n) \rfloor$. Then $\hat{d}_n+\hat{f}_n$ and $\hat{d}_n-\hat{f}_n$ are consistent estimators for $a$ and $b$ respectively. Further, there is a polynomial time estimator to calculate $\hat{d}_n$ and $\hat{f}_n$.
\end{theorem}
This theorem is extended in \cite{colin3} for the symmetric SBM with $k$ clusters, where $k$ is also estimated.
The first step needed is the following estimate.
\begin{lemma}
%Let $C_m$ be the number of $m$-cycles in $\ssbm(n,2,a/n,b/n)$. If $m=o(\log^{1/4}(n))$, then 
%\begin{align}
%\E C_m \sim \Var C_m \sim \frac{(a+b)^k + (a-b)^k}{k 2^{k+1}}
%\end{align}
Let $C_m$ be the number of $m$-cycles in $\sbm(n,p,Q/n)$. If $m=o(\log\log(n))$, then 
\begin{align}
\E C_m \sim \Var C_m \sim \frac{1}{2m} \tr (\diag(p)Q)^m.
\end{align}
\end{lemma}
To see this lemma, note that there is a cycle on a given selection of $m$ vertices with probability 
\begin{align}
&\sum_{x_1,\ldots,x_m \in [k]}\frac{Q_{x_1,x_2}}{n} \cdot \frac{Q_{x_2,x_3}}{n} \cdot \ldots \cdot \frac{Q_{x_{m},x_1}}{n} \cdot p_{x_1}\cdot \ldots \cdot p_{x_m}  \\
&=\tr (\diag(p)Q/n)^m .
\end{align}
Since there are $\sim n^m/2m$ such selections, the first moment follows. The second moment follows from the fact that overlapping cycles do not contribute to the second moment. See \cite{Mossel_SBM1} for proof details for the 2-SSBM and \cite{colin3} for the general SBM. 

Hence, one can estimate $\frac{1}{2m} \tr (\diag(p)Q)^m$ for slowly growing $m$. In the symmetric SBM, this gives enough degrees of freedom to estimate the three parameters $a,b,k$. Theorem \ref{mossel_estim} uses for example the average degree ($m=1$) and slowly growing cycles to obtain a system of equations that allows to solve for $a,b$. This extends easily to all symmetric SBMs, and the efficient part follows from the fact that for slowly growing $m$, the cycle counts coincide with the nonbacktracking walk counts with high probability \cite{Mossel_SBM1}. Note that Theorem \ref{mossel_estim} provides a tight condition for the estimation problem, i.e., \cite{Mossel_SBM1} also shows that when $(a-b)^2\le 2(a+b)$ (which we recall is equivalent to the requirement for impossibility of weak recovery) the SBM is contiguous to the Erd\H{o}s-R\'enyi model with edge probability $(a+b)/(2n)$.

However, for the general SBM, the problem is more delicate and one has to first stabilize the cycle count statistics to extract the eigenvalues of $PQ$, and use weak recovery methods to further peel down the parameters $p$ and $Q$. Deciding which parameters can or cannot be learned in the general SBM seems to be a non-trivial problem. This is also expected to come into play in the estimation of graphons \cite{choi,sbm_graphon,borgs_nips}.

\chapter{Tackling the stochastic block model}\label{tackle}
In this section, we discuss how to tackle the problem of community detection for the various recovery requirements of Section \ref{all-defs}. One feature of the SBM is that it can (and has) been viewed from many different angles. In particular, we will pursue here the algebraic and information-theoretic (or statistical) interpretations, viewing the SBM:
\begin{itemize}
\item As a low-rank perturbation model: the adjacency matrix of the SBM has low rank in expectation; thus one may hope to characterize its behavior, e.g., its eigenvectors, as perturbations of its expected counterparts. For example, the expected adjacency matrix of a two-community SBM$(n,p,Q/n)$ has the form:

{\small
\begin{align*}
\E  A =& \begin{pmatrix}
 Q_{11} \ldots Q_{11}  & Q_{12} \ldots  Q_{12}  \\
 \vdots \phantom{\ldots} \vdots & \vdots \phantom{\ldots} \vdots \\
 Q_{11} \ldots Q_{11}  &  Q_{12} \ldots  Q_{12}    \\
% &&\\
 Q_{12} \ldots Q_{12}  & Q_{22} \ldots  Q_{22}  \\
 \vdots \phantom{\ldots} \vdots & \vdots \phantom{\ldots} \vdots \\
% \vdots \phantom{\ldots} \vdots & \vdots \phantom{\ldots \,\,\, \ldots} \vdots \\
Q_{12} \ldots Q_{12}   &   Q_{22} \ldots Q_{22}   \\
% &&\\
%  Q_{13} \ldots Q_{13}  & Q_{23} \ldots Q_{23} & Q_{33} \ldots Q_{33} \\
% \vdots \phantom{\ldots} \vdots & \vdots \phantom{\ldots} \vdots & \vdots \phantom{\ldots} \vdots \\
% Q_{13} \ldots Q_{13}  &  Q_{23} \ldots Q_{23}  &  Q_{33} \ldots Q_{33}  \\
\end{pmatrix}\\ 
 &\quad \underbrace{\phantom{Q_{12} \ldots Q_{12}}}_{np_1} \,\,\,\,  \underbrace{\phantom{Q_{22} \ldots Q_{22}}}_{np_2} 
\end{align*}
}

\item As a noisy channel: the SBM graph can be viewed as the output on a channel that takes the community memberships as inputs. In particular, this corresponds to a memoryless channel encoded with a sparse graph code as in Figure \ref{meet2}.

\begin{figure}[h]
%\vspace{-0.3cm}
\begin{center}
\includegraphics[scale=.35]{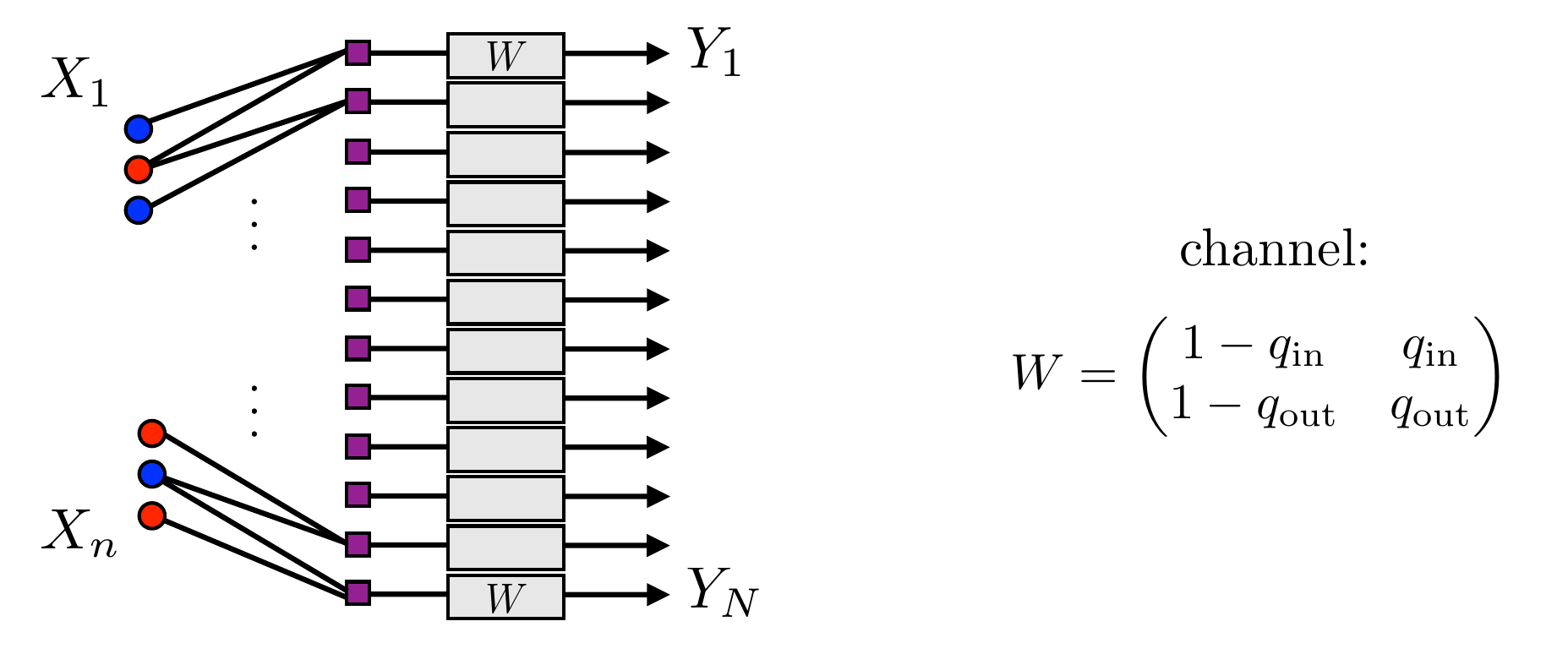}
%\vspace{-0.2cm}
\caption{The stochastic block model can be seen as an unorthodox code on a memoryless channel. In the above Figure, we illustrate the case of SSBM$(n,2,\pin,\pou)$. The information bits represent the binary community labels of the vertices. The encoder is a linear code that takes the XOR of any pair of bits, so that $N={n \choose 2}$ and the code is $R=n/N\sim 2/n$. The binary output for each channel represents the presence or absence of an edge between the corresponding vertices. The memoryless channel is determined by the matrix $W$.}
\label{meet2}
\end{center}
\vspace{-0.1cm}
\end{figure}

\end{itemize}

\section{The block MAP estimator}\label{tackle1}
A natural starting point (e.g., from viewpoint 2) is to resolve the estimation of $X$ from the noisy observation $G$ by taking the Maximum A Posteriori estimator. Upon observing $G$, if one estimates the community partition $\Omega=\Omega(X)$ with $\hat{\Omega}(G)$, the probability of error is given by
\begin{align}
P_e:= \pp\{\Omega \neq \hat{\Omega}(G)\} = \sum_{g} \pp\{ \hat{\Omega}(g) \neq \Omega | G=g\} \pp\{G=g\},
\end{align}
and an estimator $\hat{\Omega}_{\map}(\cdot)$ minimizing the above must minimize $\pp\{ \hat{\Omega}(g) \neq \Omega | G=g\}$ for every $g$.
To minimize $\pp\{ \hat{\Omega}(g) \neq \Omega  | G=g\}$, we must declare a reconstruction $\omega$ that maximizes the posterior distribution 
\begin{align}
\pp \{ \Omega= \omega | G=g \} \propto \pp \{ G=g |  \Omega= \omega\} \pp \{   \Omega= \omega\}.   \label{ap2}
\end{align}
%or equivalently 
%\begin{align}
%\sum_{x \in [k]^n: \Omega(x)=s} \pp \{ G= g | X=x \}  \prod_{i=1}^k p_{i}^{|\Omega_i(x)|}, \label{posterior}
%\end{align}
%and any such maximizer can be chosen arbitrarily. 

Consider now the strictly balanced SSBM, where $\pp \{   \Omega= \omega\}$ is the same for all equal size partitions. Then MAP is thus equivalent to the Maximum Likelihood estimator:
\begin{align}
\text{maximize } \pp \{ G=g |  \Omega= \omega\} \text{ over equal size partitions $\omega$.}  \label{newap}
\end{align}

In the two-community case, denoting by $N_{in}$ and $N_{out}$ the number of edges inside and across the clusters respectively, 
\begin{align}
\pp\{  G=g |\Omega=\omega\} &\propto  \left( \frac{\pou(1-\pin)}{\pin(1-\pou)} \right)^{N_{out}} .
\end{align}
Assuming $\pin\ge \pou$, we have $\frac{\pou(1-\pin)}{\pin(1-\pou)} \le 1$ and thus 
$$\text{MAP is equivalent to finding a min-bisection of $G$,}$$ 
i.e., a balanced partition with the least number of crossing edges. 

This brings us to a first question:
$$\text{Is it a good idea to use MAP, i.e, clusters from a min-bisection?}$$

Since MAP minimizes the probability of making an error for the reconstruction of the entire partition $\Omega$, it minimizes the error probability for exact recovery. Thus, if MAP fails in solving exact recovery, no other algorithm can succeed. In addition, MAP may not be optimal for weak recovery, since the most likely partition may not necessarily maximize the agreement. To see this, consider for example the uniform SSBM in a sparse regime with $a+b$ slightly above $2$. In this case, the graph has a giant component that contains less than half of the vertices, and there are various balanced partitions of the graph that have zero crossing edges (they separate the giant component from a collection of small disconnected components). Clearly, these are min-bisections, but they do not solve weak recovery. Nonetheless, weak recovery can still be solved in some of these cases:
\begin{lemma}
There exist $a,b$ such that weak recovery is solvable in $\ssbm(n,2,a/n,b/n)$ but block MAP fails to solve weak recovery. 
%the algorithm that outputs a random community assignment of maximum probability does not solve weak recovery.
\end{lemma}

\begin{proof}
Let $a=2.5$ and $b=0.1$. We have that $(a-b)^2>2(a+b)$, so there is an algorithm that solves weak recovery on $\ssbm(n,2,a/n,b/n)$. However, in a graph drawn from $\ssbm(n,2,a/n,b/n)$, with high probability, only about $42\%$ of the vertices are in the main component of the graph while the rest are in small components of size $O(\log n)$. So, one can partition the vertices of the graph into two equal sized sets with no edges between them by assigning every vertex in the main component and some suitable collection of small components to community $1$ and the rest to community $2$. However, the vertices in the main component are split approximately evenly between the two communities, and there is no way to tell which of the small components are disproportionately drawn from one community or the other, so for any $\epsilon>0$, each set returned by this algorithm will have less than $1/2+\epsilon$ of its vertices drawn from each community with probability $1-o(1)$.
\end{proof}
Note that such an argument is harder to establish if one restricts the min-bisection to the giant component (though we still conjecture that MAP can fail at weak recovery with this restriction). We summarize the two points obtained so far:
\begin{itemize}
\item Fact 1: If MAP does not solve exact recovery, then exact recovery is not solvable. 
\item Fact 2: Weak recovery may still be solvable when MAP does not solve weak recovery. 
\end{itemize}

\section{Computing block MAP: spectral and SDP relaxations}\label{achieve-exact}
Exactly resolving the maximization in \eqref{newap} requires comparing exponentially many terms a priori, so the MAP estimator may not always reveal the computational threshold for exact recovery. 
In fact, in the worst-case model, min-bisection is NP-hard, and approximations leave a polylogarithmic integrality gap \cite{feige_min}.
Various relaxations have been proposed for the MAP estimator. Here we review two of the main ideas. 

{\bf Spectral relaxations.}
Consider again the symmetric SBM with strictly balanced communities. 
Recall that MAP maximizes
\begin{align}
\max_{x \in \{+1,-1\}^n\atop{x^t 1^n = 0}} x^t A x ,
\end{align}
since this counts the number of edges inside the clusters minus the number of edges across the clusters, which is equivalent to the min-bisection problem (the total number of edges being fixed by the graph).
The general idea behind spectral methods is to relax the integral constraint to an Euclidean constraint on real valued vectors. 
This leads to looking for a maximizer of 
\begin{align}
\max_{x \in \mR^n : \| x\|_2^2=n  \atop{x^t 1^n = 0}} x^t A x . \label{spec2}
\end{align}
Without the constraint $x^t 1^n = 0$, the above maximization gives precisely the eigenvector corresponding to the largest eigenvalue of $A$. Note that $A 1^n$ is the vector containing the degrees of each node in $g$, and when $g$ is an instance of the symmetric SBM, this concentrates to the same value for each vertex, and $1^n$ is close to an eigenvector of $A$. Since $A$ is real symmetric, this suggests that the constraint $x^t 1^n=0$ leads the maximization \eqref{spec2} to focus on the eigenspace orthogonal to the first eigenvector, and thus to the eigenvector corresponding to the second largest eigenvalue. Thus it is reasonable to take the second largest eigenvector $\phi_2(A)$ of $A$ and round it to obtain an efficient relaxation of MAP:
\begin{align}
\hat{X}_{\mathrm{spec}} = 
\begin{cases}
1 &\text{ if } \phi_2(A) \ge 0, \\
2 &\text{ if } \phi_2(A) < 0 .
\end{cases}
\end{align}
We will discuss later on whether this is a good algorithm or not (in brief, it works well in the exact recovery regime but not in the weak recovery regime).
Equivalently, one can write the MAP estimator as a minimizer of  
\begin{align}
\min_{x \in \{+1,-1\}^n\atop{x^t 1^n = 0}} \sum_{1 \le i <j \le n} A_{ij} (x_i -x_j)^2 
\end{align}
since the above minimizes the size of the cut between two balanced clusters.
From simple algebraic manipulations, this is equivalent to looking for minimizers of  
\begin{align}
\min_{x \in \{+1,-1\}^n\atop{x^t 1^n = 0}} x^t Lx,
\end{align}
where $L$ is the Laplacian of the graph, i.e.,   
\begin{align}
L=D-A,
\end{align}
and $D$ is the degree matrix of the graph. With this approach $1^n$ is precisely an eigenvector of $L$ with eigenvalue 0, and the relaxation to a real valued vector leads directly to the second eigenvector of $L$, which can be rounded (positive or negative) to determine the communities. 
A third variant of such basic spectral approaches is to center $A$ and take the first eigenvector of $A- \frac{\pin+\pou}{2n} 1_n 1_n^t$ and round it. 

The challenge with such `basic' spectral methods is that, as the graph becomes sparser, the fluctuations in the node degrees  become more important, and this can disrupt the second largest eigenvector from concentrating on the communities (it may concentrate instead on large degree nodes). To analyze this, one may express the adjacency matrix as a perturbation of its expected value, i.e., 
\begin{align}
A&= \E A + (A-  \E A).
\end{align}
When indexing the first $n/2$ rows and columns to be in the same community, the expected adjacency matrix takes the following block structure
\begin{align}
\E  A = \begin{pmatrix}
\pin^{n/2 \times n/2} & \pou^{n/2 \times n/2} \\
\pou^{n/2 \times n/2} & \pin^{n/2 \times n/2}
\end{pmatrix},
\end{align}
where $\pin^{n/2 \times n/2}$ is the $n/2 \times n/2$ matrix with all entries equal to $\pin$. As expected, $\E A$ has two eigenvalues, the expected degree $(\pin+\pou)/2$ with the constant eigenvector, and $(\pin-\pou)/2$ with the eigenvector taking the same constant with opposite signs on each community. The spectral methods described above succeeds in recovering the true communities if the noise $Z=A-  \E A$ does not disrupt the first two eigenvectors of $A$ to be somewhat aligned with those of $\E A$. 
We will discuss when this takes place in Section \ref{exact}.

%See also Section \ref{challenge}. Theorems of random matrix theory allow to analyze this type of perturbations (see \cite{Vu-arxiv,bordenave}), most commonly when the noise is independent rather than for the specific noise occurring here, but a direct application does typically not suffice to achieve the exact recovery threshold.

{\bf SDP relaxations.} Another common relaxation can be obtained from semi-definite programming. 
We discuss again the case of two symmetric strictly balanced communities. 
%(as in \cite{abh}). 
%SDPs were also used in various works on the SBM such as \cite{new-xu,sbm-groth,levina,afonso_single,montanari_sen,adel_sbm}. 
The idea of SDPs is instead to lift the variables to change the quadratic optimization into a linear optimization (as for max-cut \cite{MXGoemans_DPWilliamson_1995}). Namely, since $\tr(AB)=\tr(BA)$ for any matrices of matching dimensions, we have  
\begin{align}
 x^t A x  = \tr (x^t A x) = \tr (A xx^t), \label{minbi}
\end{align}
hence defining $X:=xx^t$, we can write \eqref{minbi} as 
\begin{align}
\hat{X}_{\map}(g) = \argmax_{  X \succeq 0 \atop{ X_{ii} = 1, \forall i \in [n] \atop{ \rank X =1 \atop{  X 1^n = 0 }}}} \tr(AX). \label{minbi5}
\end{align}
Note that the first three constraints on $X$ force $X$ to take the form $xx^t$ for a vector $x \in \{+1,-1\}^n$, as desired, and the last constraint gives the balance requirement. The advantage of \eqref{minbi5} is that the objective function is now linear in the lifted variable $X$. The constraint $\rank X=1$ is now responsible for keeping the optimization hard. Hence, we simply remove that constraint to obtain an SDP relaxation:
\begin{align}
\hat{X}_{sdp}(g) = \argmax_{  X \succeq 0 \atop{ X_{ii} = 1, \forall i \in [n]  \atop{  X 1^n = 0 }}} \tr(AX). \label{minbib}
\end{align}
A possible approach to handle the constraint $X 1^n = 0 $ is to again use a centering of $A$. For example, one can replace the adjacency matrix $A$ by the matrix $B$ such that $B_{ij}=1$ if there is an edge between vertices $i$ and $j$, and $B_{ij}=-1$ otherwise. Using $-T$ for a large $T$ instead of $-1$ for non-edges would force the clusters to be balanced, and it turns out that $-1$ is already sufficient for our purpose. This gives another SDP:
\begin{align}
\hat{X}_{SDP}(g) = \argmax_{  X \succeq 0 \atop{ X_{ii} = 1, \forall i \in [n]  }} \tr(BX). \label{minbi2b}
\end{align}
We will further discuss the performance of SDPs in Section \ref{sdp2}.
In brief, they work well for exact recovery, and while they are suboptimal for weak recovery, they are not as senstive as vanilla spectral methods to degree variations. However, the complexity of  SDPs is significantly higher than that of spectral methods. We will also discuss how other spectral methods can afford optimality in both the weak and exact recovery regimes while preserving a quasi-linear time complexity. Notice however that we are putting the cart before the horse by talking about weak recovery now: we viewed spectral and SDP methods as relaxations of the MAP estimator, which is only an optimal estimator for exact recovery. These relaxations may still work for weak recovery, but the connection is less clear. Let us thus move to what would be the objective of merit for weak recovery.

\section{The bit MAP estimator}\label{tackle2}
If block MAP is not optimal for weak recovery, what is the right objective? To answer this more easily in the symmetric SBM, we have to in some way break the symmetry again, as done in the previous section using the partition function $\Omega(X)$. This is slightly more technical for weak recovery. We use a different trick to avoid uninteresting technicalities, and consider a weakly symmetric SBM. I.e., consider a two-community SBM with a Bernoulli prior given by $(p_1,p_2)$, $p_1 \ne p_2$, and connectivity $Q/n$ such that $PQ$ has two rows with the same sum. In other words, the expected degree of every vertex is the same (and weak recovery is non-trivial), but the model is slightly asymmetrical and one can determine the community labels from the partition. In this case, we can work with the agreement between the true labels and the algorithm reconstruction without use of the relabelling $\pi$, i.e., 
\begin{align}
A(X,\hat{X}(G))= \sum_{v=1}^n \1(X_v = \hat{X}_v(G)).
\end{align}
Consider now an algorithm that maximizes the expected agreement, i.e, 
\begin{align}
\E A(X,\hat{X}(G))= \sum_{v=1}^n \pp (X_v = \hat{X}_v(G)).
\end{align}
To solve weak recovery, one needs a non-trivial expected agreement, and to maximize the above, one has to maximize each term given by 
\begin{align}
\pp (X_i = \hat{X}_v(G))= \sum_{g}\pp (X_v = \hat{X}_v(g)|G=g)  \pp (G=g),
\end{align}
i.e., $\hat{X}_v(g)$ should take the maximal value of the function $ x_v \mapsto \pp (X_v = x_v |G=g)$. In other words, we need the marginal $\pp (X_v = \cdot |G=g)$. Note that in the symmetric SBM, this marginal is $1/2$, hence the need for the symmetry breaking.\footnote{There are different ways to break the symmetry in the symmetric SBM. One may reveal each vertex with some noisy genie; another option is to call community 1 the community that has the largest number of vertices among the $2\lfloor \sqrt{n} \rfloor+1$ largest degree vertices in the graph (we pick $2 \lfloor \sqrt{n} \rfloor +1$ to have an odd number and avoid ties).}

\section{Computing bit MAP: belief propagation}\label{bit-map}

How do we compute the marginal $\pp (X_v = x_v |G=g)$? By Bayes' rule, this requires the term $\pp (G=g|X_v = x_v)$, which can easily be obtained from the conditional independence of the edges given the vertex labels, and the marginal $\pp (G=g) = \sum_{x \in [2]^n} \pp (G=g|X = x)$, which is the non-trivial part.

Set $v_0$ to be a specific vertex in $G$. Let $v_1,... v_m$ be the vertices that are adjacent to $v_0$, and define the vectors $q_1,...,q_m$ such that for each community $i$ and vertex $v_j$, 
\[(q_j)_i=\pp(X_{v_j}=i|G\backslash \{v_0\}=g\backslash \{v_0\}).\]

Assume for a moment\footnote{This is where the symmetry breaking based on large degree vertices discussed in the previous footnote is convenient, as it allows to make the statement.} that, ignoring $v_0$, the probability distributions of $X_{v_1}, X_{v_2},..., X_{v_m}$ are asymptotically independent, i.e., for all $x_1,...,x_m$,
\begin{align}\pp \left(X_{v_1}=x_1,X_{v_2}=x_2,...,X_{v_m}
=x_m|G\backslash \{v_0\}=g\backslash \{v_0\}\right)\\
=(1+o(1))\prod_{i=1}^m \pp\left(X_{v_i}=x_i|G\backslash \{v_0\}=g\backslash \{v_0\}\right).
\end{align}
This is a reasonable assumption in the sparse SBM because the graph is locally tree-like, i.e., with probability $1-o(1)$, for every $i$ and $j$, every path between $v_i$ and $v_j$ in $G\backslash\{v_0\}$ has a length of $\Omega(\log n)$. So we would expect that knowing the community of $v_i$ would provide little evidence about the community of $v_j$.
Then, with high probability,  
\[\pp(X_{v_0}=i|G=g)=(1+o(1))\frac{p_i\prod_{j=1}^m (Q q_j)_i}{\sum_{i'=1}^k p_{i'}\prod_{j=1}^m (Q q_j)_{i'}}.\]

One can now iterate this reasoning. In order to estimate $P[X_v=i|G=g]$, one needs $P[X_{v_j}=i'|G\backslash \{v\}=g\backslash\{v\}]$ for all community labels $i'$ and all $v_j$ adjacent to $v$. In order to compute those, one needs $P(X_{v'}=i'|G\backslash\{v,v_j\}=g\backslash\{v,v_j\})$ for all $v_j$ adjacent to $v$, $v'$ adjacent to $v_j$, and every community label $i'$. To apply the formula recursively with $t$ layers of recursion, one needs an estimate of $P(X_{v_0}=i'|G\backslash \{v_1,...,v_t\}=g\backslash\{v_1,...,v_t\})$ for every path $v_0,...,v_t$ in $G$. The number of these paths is exponential in $t$, and so this approach would be inefficient. However, again due to the tree-like nature of the sparse SBM, it may be reasonable to assume that  
\begin{align}
&P(X_{v'}=i'|G\backslash\{v,v_j\}=g\backslash\{v,v_j\})\\
&= (1+o(1))P(X_{v'}=i'|G\backslash\{v_j\}=g\backslash\{v_j\}),
\end{align}
which should hold as long as there is no small cycle containing $v$, $v_j$, and $v'$. 

Therefore, using an initial estimate $(q_{v,v'})_i$ of $P(X_{v'}=i|G\backslash\{v\}=g\backslash\{v\})$ for each community $i$ and each adjacent $v$ and $v'$, we can iteratively refine our beliefs using the following algorithm which is essentially\footnote{One should normally also factor in the non-edges, but we ignore these for now as their effect is negligible in BP, although we will factor them back in when discussing linearized BP.} belief propagation (BP):

\vspace{1 cm}
\noindent
{\em Belief Propagation Algorithm (t,q, p,Q, G):}
\begin{enumerate}
\item Set $q^{(0)}=q$, where $q$ provides $(q_{v,v'})_i \in [0,1]$ for all $v,v' \in [n]$, $i \in [k]$; the initial belief that vertex $v'$ sends to vertex $v$ (one can set $q_{v',v}=q_{v'',v}$).

\item For each $0<t'< t$, each $v\in G$, and each community $i$, set
\[(q^{(t')}_{v,v'})_i=\frac{p_i\prod_{v'': (v',v'')\in E(G), v''\ne v} (Q q^{(t'-1)}_{v',v''})_i}{\sum_{i'=1}^k p_{i'}\prod_{v'': (v',v'')\in E(G), v''\ne v}(Q q^{(t'-1)}_{v',v''})_{i'}}.\]

\item For each $v\in G$ and each community $i$, set 
\[(q^{(t)}_{v})_i=\frac{p_i\prod_{v'': (v,v'')\in E(G)} (Q q^{(t-1)}_{v,v''})_i}{\sum_{i'=1}^k p_{i'}\prod_{v'': (v,v'')\in E(G)}(Q q^{(t-1)}_{v',v''})_{i'}}.\]

\item Return $q^{(t)}$.
\end{enumerate}
\vspace{1 cm}

This algorithm is efficient and terminates with a probability distribution for the community of each vertex given the graph, which seems to converge to the true distribution with enough iterations even with a random initialization. Showing this remains an open problem. Instead, we will discuss in Section \ref{linBP} how one can linearize BP in order to obtain a version of BP that can be analyzed more easily. This linearization of BP will further lead to a new spectral method on an operator called the nonbacktracking operator (see Section \ref{linBP}), which connects us back to spectral methods without the issues mentioned previously for the adjacency matrix in the weak recovery regime (i.e., the sensitivity to degree variations).

\chapter{Exact recovery for two communities}\label{exact-2}
Exact recovery for linear size communities has been one of the most studied problems for block models in its first decades. A partial list of papers is given by \cite{bui,dyer,boppana,snij,condon,mcsherry,bickel,choi,Vu-arxiv,chen-xu}. In this line of work, the approach is mainly driven by the choice of the algorithms, and in particular for the model with two symmetric communities. The results are as follows\footnote{Some of the conditions have been borrowed from attended talks and have not been double-checked.}:

{\footnotesize
\begin{center}
  \begin{tabular}{| l |  c |c   | }
  \hline
 Bui, Chaudhuri, & & $\pin = \Omega(1/n)$\\
Leighton, Sipser '84  & maxflow-mincut  & $ \pou=o(n^{-1-\frac{4}{(\pin+\pou)n}})$\\
\hline
Boppana '87  & spectral & $\frac{\pin-\pou}{\sqrt{\pin+\pou}} = \Omega(\sqrt{\log(n)/n})$\\
  \hline
 Dyer, Frieze '89  & degree min-cut & $\pin -\pou = \Omega(1)$\\
  \hline
Snijders, Nowicki '97  & EM & $\pin -\pou = \Omega(1)$\\
\hline
 Jerrum, Sorkin '98  &Metropolis & $\pin -\pou= \Omega(n^{-1/6+\epsilon})$\\
   \hline
 Condon, Karp '99 & augmentation algo. & $\pin -\pou= \Omega(n^{-1/2+\epsilon})$\\
  \hline
 Carson, Impagliazzo '01  & hill-climbing & $\pin-\pou= \Omega(\log^4(n)/\sqrt{n})$\\
 \hline
 McSherry '01  & spectral & $\frac{\pin-\pou}{\sqrt{\pin}} = \Omega(\sqrt{\log(n)/n})$\\
 \hline
 Bickel, Chen '09  & N-G modularity & $\frac{\pin-\pou}{\sqrt{\pin+\pou}} = \Omega(\log(n)/\sqrt{n})$\\
 \hline
 Rohe, Chatterjee, Yu '11  & spectral & $\pin -\pou= \Omega(1)$\\
  \hline
% Vu '14  & spectral (SVD) & $(p - q)/\sqrt{p} \geq \Omega (\sqrt{\log(n)/n })$\\
%  \hline
%Choi, Wolfe, Airoldi & ML & $p+q=\Omega(\log^3(n)/n)$ $(p -q)/\sqrt{p+q} = \Omega(1/\sqrt{n})$\\
% \hline
  \end{tabular}
\end{center}
}
More recently, \cite{Vu-arxiv} obtained a result for a spectral algorithm that works in the regime where the expected degrees are logarithmic, rather than poly-logarithmic as in \cite{mcsherry,choi}, extending results also obtained in \cite{jiaming}. 
Note that exact recovery requires the node degrees to be at least logarithmic, as discussed in Section \ref{topo}.
Thus the results of \cite{Vu-arxiv} are tight in the scaling, and the first to apply in such generality, but as for the other results in Table 1, these results do not reveal the phase transition.
%The result of Vu does however not reveal the fundamental limit for exact recovery. 
%The next result shows  exact recovery has also a phase transition that extends the connectivity one.
% and that is governed by an $f$-divergence reminiscent of Shannon's coding theorem:
The fundamental limit for exact recovery was derived first for the case of symmetric SBMs with two communities:
\begin{theorem}\label{exact1} \cite{abh_arxiv,mossel-consist}
Exact recovery in $\ssbm \left(n,2,a \frac{\log(n)}{n},b \frac{\log(n)}{n} \right)$ is solvable and efficiently so if $|\sqrt{a} - \sqrt{b}| > \sqrt{2}$ and unsolvable if $|\sqrt{a} - \sqrt{b}| < \sqrt{2}$.  
\end{theorem}
\noindent
A few remarks regarding this result:
\begin{itemize}
\item At the threshold, one has to distinguish two cases: if $a,b>0$, then exact recovery is solvable (and efficiently so) if $|\sqrt{a} - \sqrt{b}| = \sqrt{2}$, as first shown in \cite{mossel-consist}. If $a$ or $b$ are equal to 0, exact recovery is solvable (and efficiently so) if $\sqrt{b}>\sqrt{2}$ or $\sqrt{a}>\sqrt{2}$ respectively, and this corresponds to connectivity.
\item Theorem \ref{exact1} provides a necessary and sufficient condition for exact recovery, and covers all cases in $\ssbm(n,2,\pin,\pou)$ were $\pin$ and $\pou$ may depend on $n$ but not be asymptotically equivalent (i.e., $\pin/\pou \nrightarrow 1$). For example, if $\pin=1/\sqrt{n}$ and $\pou=\log^3(n)/n$, which can be written as $\pin= \frac{\sqrt{n}}{\log n} \frac{\log n}{n}$ and $\pou=\log^2 n \frac{\log n}{n}$, then exact recovery is trivially solvable as $|\sqrt{a} - \sqrt{b}|$ goes to infinity. If instead $\pin/\pou \to 1$, then one needs to look at the second order terms. This is covered by \cite{mossel-consist} for the 2 symmetric community case, which shows that for $a_n,b_n=\Theta(1)$, exact recovery is solvable if and only if $((\sqrt{a_n}-\sqrt{b_n})^2-1) \log n  + \log \log n /2 = \omega(1)$.
\item Note that $|\sqrt{a} - \sqrt{b}| > \sqrt{2}$ can be rewritten as $\frac{a+b}{2} > 1 + \sqrt{ab}$ and recall that $\frac{a+b}{2} > 2$ is the connectivity requirement in SSBM. As expected, exact recovery requires connectivity, but connectivity is not sufficient. The extra term $\sqrt{ab}$ is the `over-sampling' factor needed to go from connectivity to exact recovery, and the connectivity threshold can be recovered by considering the case where $b=0$. An information-theoretic interpretation of Theorem \ref{exact1} is also discussed in Section \ref{exact}.  
\end{itemize}

\section{Warm up: genie-aided hypothesis test}\label{genie1}
Before discussing exact recovery in the SBM, we discuss a simpler problem which will turn out to be crucial to understanding exact recovery. Namely, exact recovery with a genie that reveals all the vertices labels except for a few. 
By `a few' we really mean here one or two. 
If one works with the strictly balanced model for the community prior, then it is not interesting to reveal all vertices but a single one, as this one is forced to take the value of the community that does not have exactly $n/2$ vertices. In this case one should isolate two vertices. 
If one works with the Bernoulli model for the community prior, then one can isolate a single vertex and it is already non-trivial to decide for the labelling of that vertex given the others.

To further clean up the problem, assume for now that we have a single vertex (say vertex 0) that needs to be labelled, with $n/2$ vertices revealed in each community, i.e., assume that we have a model with two communities of size exactly $n/2$ and an extra vertex that can be in each community with probability $1/2$. 

To minimize the probability of error for this vertex we need to use the MAP estimator that picks $u$ maximizing  
\begin{align}
\pp \{X_0 = u | G=g, X_{\sim 0} = x_{\sim 0} \}. \label{targ}
\end{align}
Note that the above probability depends only on the number of edges that vertex $0$ has with each of the two communities; denoting by $N_1$ and $N_2$ these edge counts, we have   
\begin{align}
&\pp \{X_0 = u | G=g, X_{\sim 1}= x_{\sim 1} \} \\
&= \pp \{X_0 = u | N_1(G,X_{\sim 0})=N_1(g,x_{\sim 0}) \} \\
 &\propto \pp \{N_1(G,X_{\sim 0})=N_1(g,x_{\sim 0})| X_0=u \}
\end{align}
This gives an hypothesis test with two hypotheses corresponding to the two values that vertex $0$ can take, with equiprobable prior and distributions for the observable $N_1=N_1(G,X_{\sim 0})$ given by 
\begin{align}
\text{Genie-aided hypothesis test: }
\begin{cases}
X_0=1:  N_1 \sim \bin(n/2,\pin) \\
X_0=2: N_1 \sim \bin(n/2,\pou)
\end{cases}
\end{align}
The probability of error of the MAP test is then given by\footnote{Ties can be broken arbitrarily; assume that an error is declared in case of ties to simplify the expressions.}
\begin{align}
P_e := \pp \{ \mathrm{Bin}(n/2,\pin) \le \mathrm{Bin}(n/2,\pou) \}. \label{pe}
\end{align}
This is the probability that a vertex has more ``friends'' in the other community than its own.
We have the following key estimate. 
\begin{lemma}\label{ch1}
Let $\pin= a \log(n)/n$, $\pou= b \log(n)/n$. The probability of error of the genie-aided hypothesis test is given by 
\begin{align}
\pp \{ \mathrm{Bin}(n/2,\pin) \le \mathrm{Bin}(n/2,\pou) \}= n^{-\left(\frac{\sqrt{a}-\sqrt{b}}{\sqrt{2}} \right)^2 + o(1)}.
\end{align}
\end{lemma}
\begin{remark}\label{same}
The same equality holds for $\pp \{  \mathrm{Bin}(n/2,\pin) + O(1) \le \mathrm{Bin}(n/2,\pou) \}$; these are all special cases of Lemma 1 in \cite{kaiz}, initially proved in \cite{abh}. 
\end{remark}
The next result, which we will prove in the next section, reveals why this hypothesis test is crucial.
\begin{theorem}
Exact recovery in $\ssbm(n,a \log(n)/n,b \log(n)/n)$ is 
\begin{align}
\begin{cases}
\text{solvable if  $nP_e \to 0$} \\
\text{unsolvable if $nP_e \to \infty$}. 
\end{cases}
\end{align}
\end{theorem}
In other words, when the probability of error of classifying a single vertex when all others are revealed scales sublinearly, one can classify all vertices correctly whp, and when it scales supperlinearly, one cannot classify all vertices correctly whp.

%Let us consider the balanced model. We $n/2-1$ vertices in each community that are revealed and two vertices picked arbitrarily, say vertices 1 and 2, that we want to label. To preserve the balanceness, there are only two options to consider; the vertices must take different labels. Exploiting the symmetry of the model, we can further simplify this as follows.
%
%To minimize the probability of making an error at exact recovery for these two vertices, i.e., that we do not label both vertices correctly, we again need to use the MAP estimator that declares the most likely labels based on the posterior distribution:
%\begin{align}
%\pp \{X_1 = x_1, X_2 = x_2 | G=g, X_{\sim (1,2)} = x_{\sim (1,2)} \}, \eqref{targ}
%\end{align}
%and we make an error if it appears to be more likely to swap these two labels.
%Define the bad event $Bad$ such that 
%\begin{align}
%G, X_{\sim (1,2)} \in Bad
%\end{align}
%refers to the event that $G, X_{\sim (1,2)}$ take values that lead to \ref{targ} being maximized by $X_2=x_1$ and $X_1=x_2$ (without a tie). Note that  

\section{The information-theoretic threshold}
In this section, we establish the information-theoretic threshold for exact recovery in the two-community symmetric SBM with the uniform prior. That is, we assume that the two communities have size exactly $n/2$, where $n$ is even, and are uniformly drawn with that constraint. 

Recall that the MAP estimator for this model picks a min-bisection (see Section \ref{tackle1}), i.e., a partition of the vertices in two balanced groups with the least number of crossing edges (breaking ties arbitrarily). 
We will thus investigate when this estimator succeeds/fails in recovering the planted partition. 
Recall also that we work in the regime where 
\begin{align}
\pin = a \, \frac{\log n}{n}, \quad \pou = b \, \frac{\log n}{n} 
\end{align}
where the logarithm is in natural base, $a$, $b$ are two positive constants. 
%Studying this regime will allow to cover most cases, expect when $A/B \to 1$. The main result is as follows.
%\begin{theorem}
%For any $a,b \neq 0$, exact recovery is solvable (and efficiently so) if and only if $$|\sqrt{a}-\sqrt{b}| \ge \sqrt{2}.$$ 
%\end{theorem}
%We next prove this theorem, leaving aside the case where $|\sqrt{a}-\sqrt{b}| = \sqrt{2}$. 
%Note that when $b=0$, exact recovery is solvable if and only if $b>2$, which is the requirement for having connected clusters, and when $a=0$, exact recovery is solvable if and only if $a>2$, which is the requirement for having a connected bipartite graph. 

\subsection{Converse}
First recall that the term `converse' is commonly used in information theory to refer to the impossibility part of a result, i.e., when exact recovery cannot be solved in this case. Note next that exact recovery cannot be solved in a regime where the graph is disconnected with high probability, because two disconnected components cannot be correctly labelled with a probability tending to one. So this gives us already a simple condition: 

\begin{lemma}[Disconnectedness]
Exact recovery is not solvable if 
\begin{align}
\frac{a+b}{2}<1.
\end{align}
\end{lemma}
\begin{proof}
Under this condition, the graph is disconnected with high probability \cite{ER-seminal2}. 
\end{proof}
As we will see, this condition is not tight and exact recovery requires more than connectivity:
\begin{theorem}
Exact recovery is not solvable if 
\begin{align}
\frac{a+b}{2}<1 +\sqrt{ab} \quad \Longleftrightarrow \quad |\sqrt{a}-\sqrt{b}| < \sqrt{2}.
\end{align}
\end{theorem}

We will now describe the main obstruction for exact recovery.
First assume without loss of generality that the planted community is given by 
\begin{align}
x_0=(1,\dots,1,2,\dots,2),
\end{align}
with resulting communities
\begin{align}
C_1=[1:n/2], \quad C_2=[n/2:n],
\end{align} 
and let 
\begin{align}
G \sim P_{G|X}(\cdot | x_0)
\end{align}
be the SBM graph generated from this planted community assignment. 

\begin{definition}
We define the set of bad pairs of vertices by
\begin{align}
\mathcal{B}(G):=\{ (u,v) : u \in C_1, v \in  C_2, P_{G|X}(G|x_0) \le P_{G|X}(G|x_0[u \leftrightarrow v]) \},
\end{align}
where $x_0[u \leftrightarrow v]$ denotes the vector obtained by swapping the values of coordinates $u$ and $v$ in $x_0$.
\end{definition}
\begin{lemma}
Exact recovery is not solvable if $\mathcal{B}(G)$ is non-empty with non-vanishing probability.  
\end{lemma}
\begin{proof}
If there exists $(u,v)$ in $\mathcal{B}(G)$, we can swap the coordinates $u$ and $v$ in $x_0$ and increase the likelihood of the partition, thus obtaining a different partition than the planted one that is as likely as the planted one, and thus a probability of error of at least $1/2$.  
\end{proof}
We now examine the condition $P_{G|X}(G|x_0) \le P_{G|X}(G|x_0[u \leftrightarrow v])$.
This is a condition on the number of edges that vertex $u$ and $v$ have in each of the two communities. 
First note that an edge between vertex $u$ and $v$ stays in the cut if the two vertices are swapped. 
So the likelihood can only vary based on the number of edges that $u$ has in its community and in the other community ignoring $v$, and similarly for $v$. 
\begin{definition}
For a vertex $u$, define $d_+(u)$ and $d_-(u)$ as the number of edges that $u$ has in its own and other community respectively.  
For vertices $u$ and $v$ in different communities,
define $d_-(u \setminus v)$ as the number of edges that a vertex $u$ has in the other community ignoring vertex $v$.
\end{definition}
We then have
\begin{align}
&P_{G|X}(G|x_0) \le P_{G|X}(G|x_0[u \leftrightarrow v])  \\
\Longleftrightarrow \quad &d_+(u)+d_+(v) \le d_-(u \setminus v)+d_-(v \setminus u).
\end{align}
We can now define the set of bad vertices (rather than bad pairs) in each community. 
\begin{definition}
\begin{align}
\mathcal{B}_i(G):=\{ u : u \in C_i, d_+(u) \le d_-(u)-1  \}, \quad i=1,2.
\end{align}
\end{definition}
\begin{lemma}
If $\mathcal{B}_1(G)$ is non-empty with probability $1/2 + \Omega(1)$, then $\mathcal{B}(G)$ is non-empty with non-vanishing probability. 
\end{lemma}
\begin{proof}
If $u \in C_1$ and $v \in C_2$ are such that $d_+(u) \le d_-(u)-1$ and $d_+(v) \le d_-(v)-1$, then 
$d_+(u)+d_+(v) \le d_-(u)+d_-(v)-2$, and since $d_-(u) \le d_-(u \setminus v) + 1$, this implies 
$d_+(u)+d_+(v) \le d_-(u \setminus v)+d_-(v \setminus u)$.
Therefore 
\begin{align}
\pp\{ \exists (u,v) \in \mathcal{B}(G) \} &\ge \pp\{ \exists u \in \mathcal{B}_1(G),  \exists v \in \mathcal{B}_2(G) \}.
\end{align}
Using the union bound and the symmetry in the model, we have 
\begin{align}
 \pp\{ \exists u \in \mathcal{B}_1(G),  \exists v \in \mathcal{B}_2(G) \}& \ge 2 \pp\{ \exists u \in \mathcal{B}_1(G) \} - 1.
\end{align}
\end{proof}
We can now see the theorem's bound appearing: the probability that a given vertex is bad is essentially the genie-aided hypothesis test of previous section, which has a probability of $n^{-\left(\frac{\sqrt{a}-\sqrt{b}}{\sqrt{2}} \right)^2 + o(1)}$, and there are order $n$ vertices in each community, so under ``approximate independence,'' there should exists a bad vertex with probability $n^{1-\left(\frac{\sqrt{a}-\sqrt{b}}{\sqrt{2}} \right)^2 + o(1)}$ which gives the theorem's bound. We now handle the ``approximate independence'' part. Recall that if $Z$ is a positive random variable with finite variance, $\pp\{Z>0\}\ge (\E Z)^2/\E Z^2$ since by Cauchy-Schwarz $(\E Z)^2 = (\E Z \1(Z>0) )^2 \le \E Z^2 \pp\{Z>0\}$. This also implies $\pp\{Z=0\} \le \mathrm{Var}(Z)/\E Z^2$, or the weaker form $\pp\{Z=0\} \le \mathrm{Var}(Z)/(\E Z)^2$ (which can also be proved by Markov's inequality since $\pp\{Z=0\} \le \pp\{(Z - \E Z)^2 \ge (\E Z)^2\}$). 
\begin{lemma}
If $\sqrt{a}-\sqrt{b} <\sqrt{2}$, then 
\begin{align}
\pp\{ \exists u \in \mathcal{B}_1(G) \} = 1-o(1).
\end{align}
\end{lemma}
\begin{proof}
We have 
\begin{align}
\pp\{ \exists u \in \mathcal{B}_1(G) \} &= 1- \pp\{ \forall u \in C_1, u \notin \mathcal{B}_1(G) \}
\end{align}
If the events $\{u \notin \mathcal{B}_1(G)\}_{u \in C_1}$ were pairwise independent, then we would be done.
The technical issue is that for two vertices $u$ and $v$, the events are not exactly independent since the vertices can share an edge. This does not however create significant dependencies. 
Let 
\begin{align}
B_u:= \1 \{ d_+(u) \le d_-(u)-1  \}.
\end{align}
By the second moment bound, 
%\begin{align}
%\pp\{ \forall u \in C_1, u \notin \mathcal{B}_1(G) \} &= \pp\{ \sum_{u=1}^{n/2} B_u = 0 \}\\
% &\le \frac{\Var \sum_{u=1}^{n/2}  B_u}{(\E \sum_{u=1}^{n/2}  B_u)^2}
%\end{align}
%thus
\begin{align}
&\pp\{ \exists u \in \mathcal{B}_1(G) \} \ge \frac{(\E \sum_{u=1}^{n/2}  B_u)^2}{\E (\sum_{u=1}^{n/2}  B_u)^2} =\\
& \frac{((n/2)  \pp\{B_1=1\})^2 }{ (n/2)  \pp\{B_1=1\} + (n/2) (n/2 -1) \pp\{B_1=1,B_2=1\}}
\end{align}
and the last bound tends to 1 if the following three conditions hold
\begin{align}
n \pp\{B_1=1\} &= \omega(1),\\
\frac{\pp\{B_1=1| B_2=1\}}{\pp\{B_1=1\} } &=1+o(1).%,\\
%\frac{\pp\{B_1=1| B_{n/2+1}=1\}}{\pp\{B_1=1\} } &=1+o(1).
\end{align}
The first condition follows from the genie-aided hypothesis test\footnote{In this case, one of the two Binomial random variables has $n-1$ trails rather than $n$, with a trial replace by 1, which makes no difference in the result as mentioned in Remark \ref{same}.} and reveals the bound in the theorem; the second condition amounts to say that $B_1$ and $B_2$ are asymptotically independent. 

We now show the second condition. We have
\begin{align}
&\pp\{B_1=1| B_2=1\} = \pp\{ d_+(1) \le d_-(1)-1 | d_+(2) \le d_-(2)-1 \} \\
&= \pp\{ B(n/2-2,\pin) + B_{1,2} \le B(n/2,\pou)-1 \\ & \qquad \,\, | B'(n/2-2,\pin) + B_{1,2} \le B'(n/2,\pou)-1 \} \label{12}
\end{align}
where $B(m,C)$ or $B'(m,C)$ denotes a Binomial random variable with $m$ trials and success probability $C$, $B_{1,2}$ is Bernoulli$(\pin)$, and the five different random variables appearing in \eqref{12} are mutually independent. Thus, 
\begin{align}
&\pp\{B_1=1| B_2=1\}   \\
&= \sum_{b=0,1} \pp\{ B(n/2-2,\pin) + b \le B(n/2,\pou)-1  \\ & \qquad \qquad | B'(n/2-2,\pin) + b \le B'(n/2,\pou)-1 , B_{1,2}=b\} \\
&\cdot \pp\{ B_{1,2}=b | B'(n/2-2,\pin) + B_{1,2} \le B'(n/2,\pou) -1\}\\
&= \sum_{b=0,1} \pp\{ B(n/2-2,\pin) + b \le B(n/2,\pou)-1 \}  \\ & \qquad  \cdot \pp\{ B_{1,2}=b | B'(n/2-2,\pin) + B_{1,2} \le B'(n/2,\pou)-1 \}. \label{b12}
\end{align}
Let 
\begin{align}
\Delta := B(n/2-2,\pou) - B(n/2-2,\pin)-1.
\end{align}
We have 
\begin{align}
\pp\{ B_{1,2}=b |  B_{1,2} \le \Delta \} = \frac{\pp\{ B_{1,2}=b \} \pp\{ b \le \Delta \}}{\sum_{b=0,1} \pp\{ B_{1,2}=b \} \pp\{ b \le \Delta \} }
\end{align}
and since $\pp\{ B_{1,2}=1 \} \le \pp\{ B_{1,2}=0 \}= 1-o(1)$ and $\pp\{ 1 \le \Delta \} \le \pp\{ 0 \le \Delta \}$, we have 
\begin{align}
&\sum_{b=0,1} \pp\{ B_{1,2}=b \} \pp\{ b \le \Delta \} = \pp\{ 0 \le \Delta \} (1+o(1)) ,\\
&\pp\{ B_{1,2}=0 |  B_{1,2} \le \Delta \} = \pp\{ B_{1,2}=0 \} (1+o(1))=1+o(1) ,\\
&\pp\{ B_{1,2}=1 |  B_{1,2} \le \Delta \}  
%\frac{ \pp\{ B_{1,2}=1 \} \pp\{ 1 \le \Delta \}}{\pp\{ 0 \le \Delta \}}(1+o(1))
=o(1)
\end{align}
Thus 
\begin{align}
 \eqref{b12} =  \pp\{ 0 \le \Delta \} (1+o(1)) + \pp\{ 1 \le \Delta \} o(1) =  \pp\{ 0 \le \Delta \} (1+o(1)).
\end{align}
On the other hand, for a random variable $B'$ that is Bernoulli$(\pin)$ and independent of $\Delta$, we have 
\begin{align}
 &\pp\{B_1=1\} = \pp\{ B(n/2-2,\pin) + B' \le B(n/2,\pou) -1 \} \\
 &=  \pp\{B'  \le \Delta \} = \pp\{0  \le \Delta \} (1+o(1)) 
\end{align}
Thus, 
\begin{align}
\frac{\pp\{B_1=1| B_2=1\}}{\pp\{B_1=1\} } &=1+o(1),
\end{align}
which concludes the proof. 
\end{proof}

\subsection{Achievability}\label{ach1}
The next result shows that the previous bound is tight. 
\begin{theorem}
Exact recovery is  solvable if 
\begin{align}
|\sqrt{a}-\sqrt{b}| > \sqrt{2}.
\end{align}
\end{theorem}
We will discuss the boundary case $|\sqrt{a}-\sqrt{b}| > \sqrt{2}$ later on (it is still possible to solve exact recovery in this case as long as both $a$ and $b$ are non-zero).
To prove this theorem, one can proceed with different approaches:

\begin{enumerate}
\item {\it Showing k-swaps are dominated by 1-swaps. } The converse shows that below the threshold, there exists a bad pair of vertices in each community that can be swapped (and thus placed in the wrong community) while increasing the likelihood (i.e., reducing the cut). To show that the min-bisection gives the planted bisection, we need to show instead that there is no possibility to swap $k$ vertices from each community and reduce the cut {\it for any} $k \in \{1,\dots,n/4\}$ (we can use $n/4$ because the communities have size $n/2$). That is, above the threshold, 
\begin{align}
&\pp\{ \exists S_1 \subseteq C_1, S_2 \subseteq C_2 : \\& \quad |S_1|=|S_2|, d_+(S_1) + d_+(S_2) \le d_-(S_1\setminus S_2) + d_-(S_2 \setminus S_1)  \} \\&= o(1).
\end{align}
For $T_1 \subseteq C_1, T_2 \subseteq C_2$ such that $|T_1|=|T_2| = k$, define 
\begin{align}
P_e(k):=\pp\{ |d_+(T_1) + d_+(T_2) \le d_-(T_1 \setminus T_2) + d_-(T_2 \setminus T_2)  \} .
\end{align}
Then, by a union bound, 
\begin{align}
&\pp\{ \exists S_1 \subseteq C_1, S_2 \subseteq C_2 : \\& \quad |S_1|=|S_2|, d_+(S_1) + d_+(S_2) \le d_-(S_1 \setminus S_2) + d_-(S_2 \setminus S_1)  \} \\
&= \pp\{ \exists k \in [n/4], S_1 \subseteq C_1, S_2 \subseteq C_2 : |S_1|=|S_2|=k, \\
&\phantom{\pp\{ \exists k \in [n/4], } d_+(S_1) + d_+(S_2) \le d_-(S_1) + d_-(S_2)  \} \\
&\le \sum_{k=1}^{n/4} {n/4 \choose k} {n/4 \choose k} P_e(k)\\
&= (n/4)^2 P_e(1) + R 
\end{align}
where $R :=\sum_{k=2}^{n/4} {n/4 \choose k} {n/4 \choose k} P_e(k)$.
Note that $P_e(1)$ behaves like the error probability of the genie-aided test squared (we look at two vertices instead of one), and one can show that the product $n^2 P_e(1)$ is vanishing above the threshold. So it remains to show that the reminder $R$ is also vanishing, and in fact, one can show that $R = O(n^2 P_e(1))$, i.e., the first term (1-swaps) dominates the other terms (k-swaps). This approached is used in \cite{abh}.
\item {\it Using graph-splitting.} The technique of graph-splitting is developed in \cite{abh,colin1focs} to allow for multi-round methods, where solutions are successively refined. The idea is to split the graph $G$ into two new graphs $G_1$ and $G_2$ on the same vertex set, by throwing each edge independently from $G$ to $G_1$ with probability $\gamma$ and keeping the other edges in $G_2$. In a sparse enough regime, such as logarithmic degrees, one can further treat the two graphs as essentially independent SBMs on the same planted community. Taking $\gamma=\log\log(n)/\log(n)$, one obtains for $G_1$ an SBM with degrees that scale with $\log\log(n)$, and it is not too hard to show that almost exact recovery can be solved in such a diverging-degree regime. One can then use the almost exact clustering obtained on $G_1$ and refine it using the edges of $G_2$, using a genie-aided test for each vertex, where the genie is not an exact genie as in Section \ref{genie1}, but an almost-exact genie obtained from $G_1$. One then shows that the almost-exact nature of the genie does not change the outcome, and the same threshold emerges. This approach is discussed in more details in Section \ref{exact} when we consider exact recovery in the general SBM. 
\item {\it Using the spectral algorithm.} While it is not necessary to use an efficient algorithm to establish the information-theoretic threshold, the spectral algorithm offers a nice algebraic intuition to the problem. This approach is discussed in detail in next section.  
\end{enumerate}

\section{Achieving the threshold}
\subsection{Spectral algorithm}
In this section, we show that the vanilla spectral algorithm discussed in Section \ref{achieve-exact} achieves the exact recovery threshold. The proof is based on \cite{kaiz}. Recall that the algorithm is a relaxation of the min-bisection, changing the integral constraint on the community labels to an Euclidean-norm constraint.
%\begin{align}
%\max_{x \in \mR^n : \| x\|_2^2=n  \atop{x^t 1^n = 0}} x^t A x .
%\end{align}
This suggest that taking the second largest eigenvector of $A$, i.e., the eigenvector corresponding to the second largest eigenvalue of $A$, and rounding it, gives a plausible reconstruction. 
Techniques from random matrix theory are naturally relevant here, as used in various works such as \cite{vu-rmt,rao-newman,vu3,Vu-arxiv,prout_colt,wein1} and references therein.

%\begin{align}
%X_{\mathrm{spec}}(G) = 
%\begin{cases}
%1 &\text{ if } \phi_2(A) \ge 0 \\
%2 &\text{ if } \phi_2(A) < 0 
%\end{cases}
%\end{align}
%gives a plausible reconstruction. 

For the rest of the section, we write $p := \pin$ and $q := \pou$ for simplicity. Denote by $A'$ the adjacency matrix of the graph with self-loops added with probability $p$ for each vertex. 
Therefore, 
\begin{align}
\E A' = n\frac{p+q}{2}  \bar{\phi}_1 \bar{\phi}_1^t + n\frac{p-q}{2} \bar{\phi}_2 \bar{\phi}_2^t
\end{align}
where 
\begin{align}
\bar{\phi}_1&= 1^n/\sqrt{n},\\
\bar{\phi}_2&= [(-1)^{1\{X_1=1\}}, \ldots, (-1)^{1\{X_n=1\}}] / \sqrt{n}.
\end{align}
In words, $\bar{\phi}_2$ is a vector whose signs indicate the assignment of $X$. We will work with the ``centered'' adjacency matrix\footnote{Strictly speaking, obtaining the centered adjacency matrix requires that $p+q$ is known. A better approach is to replace $p+q$ by an estimate, or use the second eigenvector of $A'$ directly as analyzed in \cite{kaiz}.} which we denote by $A$ in this section (with an abuse of notation compared to previous sections) where we subtract the top expected eigenvector:
\begin{align}
A:= A' -  n\frac{p+q}{2}  \bar{\phi}_1 \bar{\phi}_1^t.
\end{align}
%{\color{red} KW: using $A$ and $B$ may cause confusions, as the centered adjacency matrix is also denoted by $A$. How about using $p$ and $q$ instead?}
The slight advantage is that $A$ is now rank 1 in expectation:
\begin{align}
\bar{A}:=\E A&= n\frac{p-q}{2} \bar{\phi}_2 \bar{\phi}_2^t \\
&= \bar{\lambda} \bar{\phi} \bar{\phi}^t 
\end{align}
where we renamed 
\begin{align}
\bar{\lambda} &:=\frac{(a-b)\log (n)}{2}\\
\bar{\phi} &:= \bar{\phi}_2
\end{align}

We now want to show that the top eigenvector $\phi$ of $A$ has all its components aligned with $\bar{\phi}$ in terms of signs (up to a global flip). Define, for $i \in [n]$,
\begin{align}
X_{\mathrm{spec}}(i) = 
\begin{cases}
1 &\text{ if } \phi(i) \ge 0, \\
2 &\text{ if } \phi(i) < 0 .
\end{cases}
\end{align}

\begin{theorem}\label{main_spec}
The spectral algorithm that outputs $X_\mathrm{spec}$ solves exact recovery when  
\begin{align}
|\sqrt{a}-\sqrt{b}| > \sqrt{2}.
\end{align}
\end{theorem}
Note that the algorithm runs in polynomial time in $n$, at most $n^3$ counting loosely and less using the sparsity of $A$. 
Also note that we do not have to worry about the case where the resulting community is not balanced, as this enters the vanishing error probability. 
Another way to write the theorem is as follows:
\begin{theorem}\label{spec_thm}
In the regime $p = a \, \frac{\log n}{n}$, $q = b \, \frac{\log n}{n}$, $a \neq b$, $a,b > 0$, $|\sqrt{a}-\sqrt{b}| \neq  \sqrt{2}$, 
\begin{align}
\pp \{ X_{\mathrm{spec}} \equiv X_{\mathrm{MAP}} \} = 1-o(1) \text{ whenever } \pp \{ X \equiv X_{\mathrm{MAP}} \} = 1-o(1),
\end{align}
i.e., the spectral and MAP estimators are equivalent whenever MAP succeeds at recovering the true $X$ (we use $x\equiv y$ in the above for $x \in \{y,y^c\}$ where $y^c$ flips the components of $y$, due to the usual symmetry).
\end{theorem}
This an interesting phenomenon that seems to take place for more than one problem in combinatorial statistics, see for example discussion in \cite{kaiz}.

We now proceed to prove Theorem \ref{spec_thm}.
We will break the proof in several parts. A first important result due to \cite{ofek} gives a bound on the norm of $A - \bar{A}$ (see \cite{kaiz} for a proof):
\begin{lemma}\label{feige2}\cite{ofek}
For any $a,b > 0$, there exist $c_1,c_2 >0$ such that 
\begin{align}
\pp\{ \| A - \bar{A} \|_2 \ge c_1\sqrt{\log(n)} \} \le c_2 n^{-3}.
\end{align}
\end{lemma}
This implies a first reassuring fact, i.e., the first eigenvalue $\lambda$ of $A$ is in fact asymptotic to $\bar{\lambda}$ in our regime. This follows from the Courant-Fisher Theorem:
\begin{lemma}[Courant-Fisher or Weyl's Theorem]
The following holds surely,
\begin{align}
| \lambda - \bar{\lambda} | \le \| A - \bar{A} \|_2.
\end{align}
\end{lemma}
This implies that $\lambda \sim \bar{\lambda}$ with high probability, since $\bar{\lambda} \asymp \log (n)$. However, this does not imply anything for the eigenvectors yet. A classical result to convert bounds on the norm $A - \bar{A}$ to eigenvector alignments is the Davis-Kahan Theorem. Below we present a user-friendly version based on the Lemma 3 in \cite{kaiz}.
\begin{theorem}[Davis-Kahan Theorem]\label{thm-DK}
Suppose that $\bar{H}=\sum_{j=1}^n \bar{\mu}_j \bar{u}_j \bar{u}_j^t$ and $H=\bar{H}+E$, where $\bar{\mu}_1 \geq \cdots\geq \bar{\mu}_n$, $\| \bar{u}_j \|_2=1$ and $E$ is symmetric. Let $u_j$ be an eigenvector of $H$ corresponding to its $j$-th largest eigenvalue, and $\Delta = \min\{ \bar{\mu}_{j-1} - \bar{\mu}_j , \bar{\mu}_j - \bar{\mu}_{j+1} \}$, where we define $\bar{\mu}_0 = +\infty$ and $\bar{\mu}_{n+1} = -\infty$. We have
\begin{align}
\min_{s=\pm 1} \| s u_j - \bar{u}_j \|_2 \lesssim \frac{ \| E\|_2 }{\Delta}.
\end{align}
In addition, if $\|E \|_2\leq \Delta$, then
\begin{align}
\min_{s=\pm 1} \| s u_j - \bar{u}_j \|_2 \lesssim \frac{ \| E \bar{u}_j \|_2 }{ \Delta },
\end{align}
where both $\lesssim$ only hide absolute constants.
\end{theorem}
\begin{corollary}
For any $a,b > 0$, with high probability, 
\begin{align}
|\langle \phi,\bar{\phi} \rangle | = 1-o(1).
\end{align}
\end{corollary}
While this gives a strong alignment, it does not give any result for exact recovery. One can use a graph-splitting step to leverage this result into an exact recovery result by using a cleaning phase on the eigenvector (see item 2 in Section \ref{ach1}). Interestingly, one can also proceed directly using a sharper spectral analysis, and show that the sign of the eigenvector $\phi$ directly recovers the communities. This was done in \cite{kaiz} and is covered below.

A first attempt would be to show that $\phi$ and $\bar{\phi}$ are close enough in each component, i.e., that with high probability, 
\begin{align}
&|\phi_i - \bar{\phi}_i | \stackrel{?}{<} |\bar{\phi}_i |, \quad \forall i \in [n] \label{hope} \\
\Longleftrightarrow \quad  &\|\phi  - \bar{\phi} \|_{\infty} \stackrel{?}{<} 1/\sqrt{n}
\end{align}
or $\|\phi - (-\bar{\phi}) \|_{\infty} \stackrel{?}{<} 1/\sqrt{n}$ since we must allow for a global flip due to symmetry.
This would imply that rounding the components of $\phi$ to their signs would produce with high probability the same signs as $\bar{\phi}$ (or $-\bar{\phi}$), which solves exact recovery. 

Unfortunately the above inequality does not hold all the way down to the exact recovery threshold, which makes the problem more difficult. However, note that it is not necessary to have \eqref{hope} in order to obtain the correct sign by rounding $\phi$: one can have a large gap for $|\phi_i - \bar{\phi}_i |$ which still produces the good sign as long as this gap in ``on the right side," i.e., $\phi_i$ can be much larger than $\bar{\phi}_i$ if $\bar{\phi}_i$ is positive and $\phi_i$ can be much smaller than $\bar{\phi}_i$ if $\bar{\phi}_i$ is negative (or the reverse statement for $-\bar{\phi}$). This is illustrated in Figure \ref{fig:intro} and is shown below. 

\begin{figure}[h!]
	\centering
	\includegraphics[scale=0.35]{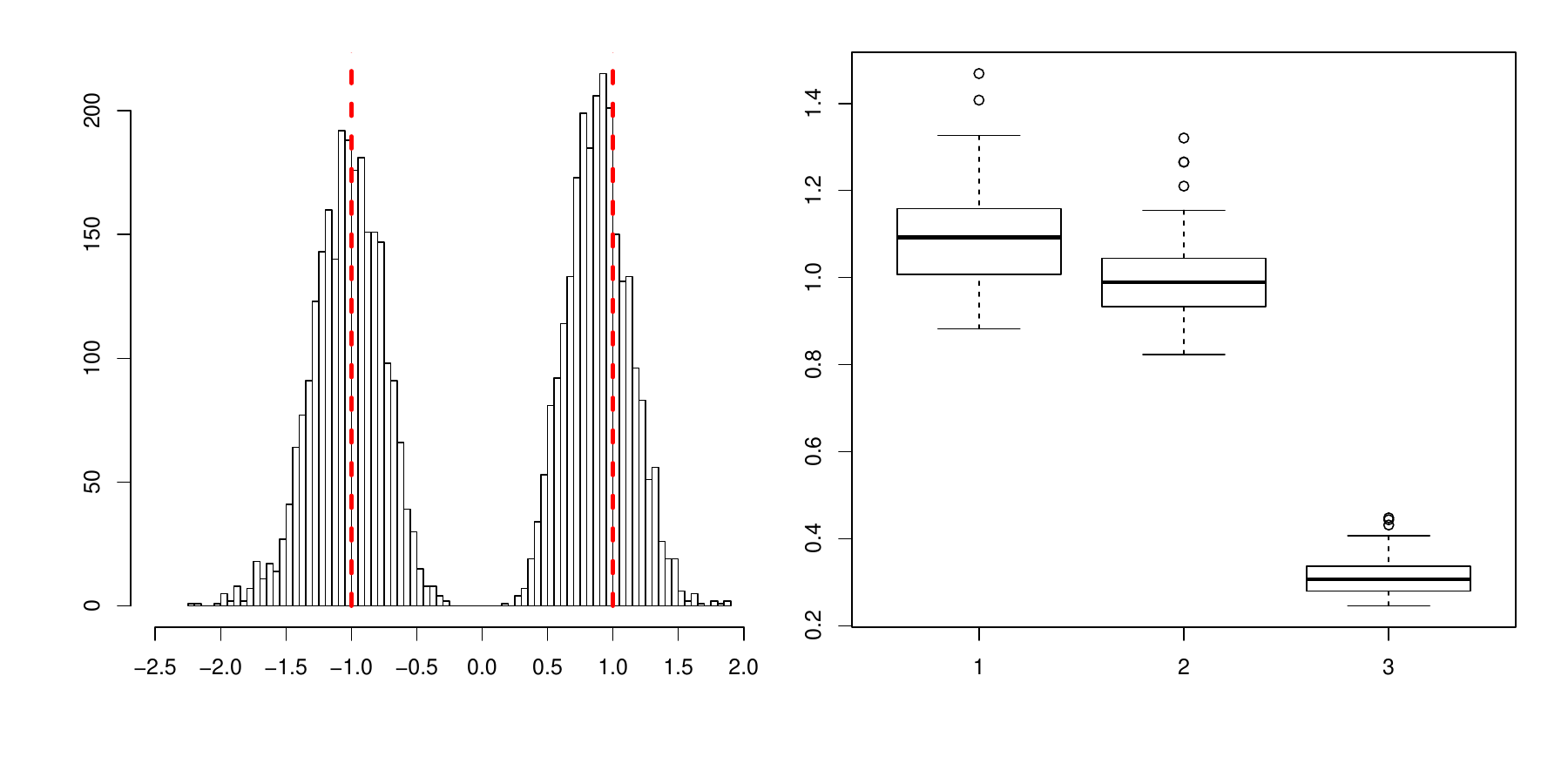}
	\caption{For the original uncentered adjacency matrix $A'$, the plot shows the second eigenvector $u_2$ and its first-order approximation in SBM. For the centered adjacency matrix as analyzed here, the top eigenvector shows the same phenomenon. The left plot: The histogram of coordinates of $\sqrt{n} \,u_2$ computed from a single realization of adjacency matrix $A$, where $n$ is $5000$, $a$ is $4.5$ and $b$ is $0.25$. Under this regime, exact recovery is expected, and indeed coordinates form two well-separated clusters. The right plot: boxplots showing three different distance/errors (up to sign) over $100$ realizations. (1) $ \sqrt{n}\,  \| u_2 - u_2^* \|_{\infty}$, (2) $\sqrt{n}\,  \| Au_2^* / \lambda_2^* - u_2^*  \|_{\infty}$, (3) $\sqrt{n}\,  \|   u_2 - Au_2^* / \lambda_2^*  \|_{\infty}$. These boxplots show that $Au_2^* / \lambda_2^*$ is a good approximation of $u_2$ under $\ell_{\infty}$ norm even though $\| u_2 - u_2^* \|_{\infty}$ may be large.  }\label{fig:intro}
\end{figure}

The main idea is to show that the components of $\phi$ are well-approximated by the components of $A \bar{\phi}/\bar{\lambda}$ (rather than $\bar{\phi}$), i.e., 
\begin{align}
\phi = A {\phi}/{\lambda} \approx  A \bar{\phi}/\bar{\lambda} = \bar{\phi} + (A-\E A)\bar{\phi}/\bar{\lambda}. 
\end{align}
Formally:
\begin{theorem}\label{thm-approx}
Let $a > b>0$. There exist constants $C$ and $c$ such that for sufficiently large $n$,
\begin{align}
\pp \left( \min_{s = \pm 1} \| s \phi - A\bar{\phi}/\bar{\lambda} \|_{\infty} \leq \frac{c}{ \sqrt{n} \log \log n} \right) \geq 1- Cn^{-2}.
\end{align}
\end{theorem}
Note that $A\bar{\phi}$ is familiar to us: it contains exactly the random variable entering the error event for the genie-aided hypothesis test in \eqref{pe}, i.e., 
\begin{align}
\sqrt{n}\, (A\bar{\phi})_i \sgn(\bar{\phi}_i) \stackrel{(d)}{=} \mathrm{Bin}(n/2,p) - \mathrm{Bin}(n/2,q).
\end{align}
This is because 
\begin{align}
\sqrt{n}\, (A\bar{\phi})_i \sgn(\bar{\phi}_i) = \sqrt{n}\, (A'\bar{\phi})_i \sgn(\bar{\phi}_i) = \sum_{j=1}^n A'_{ij} \sgn( \bar{\phi}_i \bar{\phi}_j),
\end{align}
and each $A'_{ij}$ is an independent Bernoulli variable whose success probability depends on whether $X_i = X_j$. We know from Lemma \ref{ch1} has probability $n^{-1 -\Omega(1)}$ to be negative (i.e., to move to ``the other side'') above the exact recovery threshold. Since $A\bar{\phi}$ is normalized by $\bar{\lambda}$, we will use the stronger version of  Lemma \ref{ch1} mentioned in Remark \ref{same}. We now give the proof for Theorem \ref{main_spec}, and then proceed to proving Theorem \ref{thm-approx}.

\begin{proof}[Proof of Theorem \ref{main_spec}.]
Define the following events: 
%{\color{red} KW: The sign flip applies to $\phi$, not $A \bar{\phi} / \bar{\lambda}$. Hence in the event $\mathcal{E}_1$ there is no need to consider $\sgn( A\bar{\phi}/\bar{\lambda} ) = - \sgn( \bar{\phi} )$ as with high probability it will not happen.  Hence $\mathcal{E}_1$ could be conveniently rewritten as
%$\left\{ \min_{i\in[n]} (A \bar{\phi}/\bar{\lambda})_i \bar{\phi}_i \geq \frac{2\varepsilon}{(a-b)\sqrt{n}}\right\}$, which facilitates the application of binomial tail bounds.
%}
\begin{align}
\mathcal{E}_1&:=\left\{ \min_{i\in[n]} (A \bar{\phi}/\bar{\lambda})_i \sgn(\bar{\phi}_i) \geq \frac{2\varepsilon}{(a-b)\sqrt{n}}\right\} \\
\mathcal{E}_2&:=\left\{ \min_{s = \pm 1} \| s \phi - A\bar{\phi}/\bar{\lambda} \|_{\infty} \leq \frac{c}{ \sqrt{n} \log \log n}\right\}.
\end{align}
Note that Theorem \ref{thm-approx} says that $\mathcal{E}_2$ takes place with high probability.
Therefore if $\mathcal{E}_1$ takes place with high probability as well, the event
\begin{align}
\sgn(\phi)= \pm \sgn(\bar{\phi})
\end{align}
must take place with high probability, because entrywise (up to a global flip) $\phi$ is at distance $O(\frac{1}{ \sqrt{n} \log \log n})$  from $A\bar{\phi}/\bar{\lambda}$, and $A\bar{\phi}/\bar{\lambda}$ is at distance $\Omega(\frac{1}{ \sqrt{n}})$ from the origin, so the noise that takes $\bar{\phi}$ to $\phi$ cannot make $\phi$ across the origin and change sign since $|\bar{\phi_i}|=1/\sqrt{n}$ (though it could take $\phi$ far on the other side). 

We now show that $\mathcal{E}_1$ takes place with high probability.
We have 
\begin{align}
&\pp \left( \mathcal{E}_1 \right)  \\
&\ge  1-\sum_{i=1}^n\, \pp\left( (A \bar{\phi}/\bar{\lambda})_i \sgn(\bar{\phi}_i) < \frac{2\varepsilon}{(a-b)\sqrt{n}} \right) \\
&=  1- n\, \pp\left( (A \bar{\phi}/\bar{\lambda})_1 \sgn(\bar{\phi}_1) < \frac{2\varepsilon}{(a-b)\sqrt{n}} \right).
%&= 1-o(1).
\end{align}
%{\color{red} KW: with the new $\mathcal{E}_1$ we can use union bounds to write
%\begin{align}
%&\pp \left( \mathcal{E}_1 \right)  = \pp \left( \min_{i\in[n]} (A \bar{\phi}/\bar{\lambda})_i \bar{\phi}_i \geq \frac{2\varepsilon}{(a-b)\sqrt{n}} \right)
%\geq 1-n \pp \left( (A \bar{\phi}/\bar{\lambda})_1 < \frac{2\varepsilon}{(a-b)\sqrt{n}} \right)
%\end{align}
%and then apply the tail bound.}

When $\sqrt{a} - \sqrt{b} > \sqrt{2}$, we can choose some $\e = \e (a,b)>0$ such that $ (\sqrt{a} - \sqrt{b})^2/2 - \varepsilon \log(a/b)/2 > 1$. Thus from (a strengthening of) the genie-aided error bound, 
\begin{align}
\pp \left(  |(A \bar{\phi}/\bar{\lambda})_1| \le \frac{2\varepsilon}{(a-b)\sqrt{n}} \right)  \le n^{- (\sqrt{a} - \sqrt{b})^2/2 + \varepsilon \log(a/b)/2 } = n^{-1-\Omega(1)}
\end{align}
and 
%a fortiori, 
%\begin{align}
%\pp \left(   \sgn( A\bar{\phi}/\bar{\lambda} ) = \pm \sgn( \bar{\phi} ) \right) =1 -o(1),
%\end{align}
%since each component of $A\bar{\phi}/\bar{\lambda}$ is at least strictly positive or strictly negative (depending on the global flip on $\bar{\phi}$) with probability $n^{-1-\Omega(1)}$. 
therefore, 
\begin{align}
&\pp \left( \mathcal{E}_1 \right)  =1-o(1).
\end{align}
\end{proof}
We now proceed with the proof of Theorem \ref{thm-approx}.
\begin{proof}[Proof of Theorem \ref{thm-approx}.]
To simplify the notation, assume that $\phi$ is chosen such that $\phi^t\bar{\phi} \ge 0$, so that we can remove the sign variable $s$. We want to obtain a bound on 
%\begin{align}
%\| \phi - A\bar{\phi}/\bar{\lambda} \|_{\infty} &=  \| A\phi/\lambda - A\bar{\phi}/\bar{\lambda} \|_{\infty}   \\
%&= \| A\phi/\lambda - A\bar{\phi}/\lambda +A\bar{\phi}/\lambda - A\bar{\phi}/\bar{\lambda} \|_{\infty} \\
%&\le \| A\phi/\lambda - A\bar{\phi}/\lambda \|_{\infty}  + \| A\bar{\phi}/\lambda - A\bar{\phi}/\bar{\lambda} \|_{\infty} \\
%&= \frac{1}{\lambda} \| A(\phi - \bar{\phi}) \|_{\infty}  +  \frac{|\bar{\lambda} - \lambda|}{\lambda \bar{\lambda}} \| A\bar{\phi} \|_{\infty}.
%\end{align}
\begin{align}
\| \phi - A\bar{\phi}/\bar{\lambda} \|_{\infty} &= \| \phi -  A \phi/
\bar{\lambda} +A \phi/\bar{\lambda} - A\bar{\phi}/\bar{\lambda} \|_{\infty} \\
&\le \| \phi  -  A \phi/ \bar{\lambda}  \|_{\infty} + \| A \phi/\bar{\lambda} - A\bar{\phi}/\bar{\lambda} \|_{\infty} \\
&= \frac{|\bar{\lambda} - {\lambda}|}{\bar{\lambda}} \| \phi \|_{\infty} +   \frac{1}{ \bar{\lambda}} \| A(\phi - \bar{\phi}) \|_{\infty} . \label{broken}
\end{align}
Let us define the event 
\begin{align}
\mathcal{E}:= \{\| A - \bar{A} \|_2 \le c_1\sqrt{\log(n)} \}. \label{event}
\end{align}
Recall that Weyl's theorem gives $| \lambda - \bar{\lambda} | \le \| A - \bar{A} \|_2$, and thus under the event $\mathcal{E}$, which takes place with high probability by Lemma \ref{feige2}, we must have $| \lambda - \bar{\lambda} |  = O(\sqrt{\log(n)})$. Given that $\bar{\lambda} = \Theta(\log(n))$, we must have under $\mathcal{E}$  
\begin{align}
\bar{\lambda}/2 \le \lambda \le 2\bar{\lambda} \label{order1}
\end{align}
and a fortiori $\frac{|\bar{\lambda} - {\lambda}|}{\bar{\lambda}}=O(1/\sqrt{\log(n)})$. Therefore, under $\mathcal{E}$, we have that the first term in \eqref{broken} is bounded as 
\begin{align}
\frac{|\bar{\lambda} - {\lambda}|}{\bar{\lambda}} \| \phi \|_{\infty} \le O(\| \phi \|_{\infty}/\sqrt{\log(n)}). \label{plug3}
\end{align}
Before worrying about estimating $\| \phi \|_{\infty}$, we move to the second term in \eqref{broken}.
One difficulty in estimating $\| A(\phi - \bar{\phi}) \|_{\infty}$ is that $A$ and $(\phi - \bar{\phi})$ are dependent since $\phi$ is an eigenvector of $A$. Thus, to bound the $m$-th component of $A(\phi - \bar{\phi})$, namely $A_m(\phi - \bar{\phi})$, where $A_m$ is the $m$-row of $A$, we cannot directly use a concentration result that applies to expressions of the kind $A_m w$ where $w$ is an independent test vector. To decouple the dependencies, we use a leave-one-out technique, as used for example in \cite{bean2013optimal}, \cite{JMo15}, and \cite{ZhoBou17}. 

Define $n$ auxiliary
matrices $\{A^{(m)}\}_{m=1}^{n}\subseteq\mathbb{R}^{n\times n}$ as
follows: for any $m\in[n]$, let 
\[
(A^{(m)})_{ij}=A_{ij} \delta_{\{i\neq m,j\neq m\}},\qquad \forall i,j\in[n]
\]
where $\delta_A$ is the indicator function on the event $A$. Therefore, $A^{(m)}$ is obtained from $A$ by zeroing out the $m$-th row and column. Let $\phi^{(m)}$ be the leading eigenvector of $A^{(m)}$. Again, $\phi^{(m)}$ is chosen such that $(\bar{\phi})^{t} \phi^{(m)}\geq 0$. Denoting the $m$th row vector as $A_m$, we can write 
\begin{align}
(A (\phi - \bar{\phi}))_m = A_m (\phi - \bar{\phi}) = A_m (\phi - \phi^{(m)}) + A_m (\phi^{(m)} -\bar{\phi})
\end{align}
and thus 
\begin{align}
| (A (\phi - \bar{\phi}))_m | &\le | A_m (\phi - \phi^{(m)}) | + | A_m (\phi^{(m)} -\bar{\phi}) | \\
& \le  \| A_m \|_2 \| \phi - \phi^{(m)} \|_2 + | A_m (\phi^{(m)} -\bar{\phi}) | \\
& \le  \underbrace{\| A \|_{2 \to \infty} \| \phi - \phi^{(m)} \|_2}_{T_1} + \underbrace{| A_m (\phi^{(m)} -\bar{\phi}) |}_{T_2} \label{dec}
\end{align}
where we used the Cauchy-Schwarz Inequality in the first inequality, and $ \| A \|_{2 \to \infty} : = \max_{m \in [n]} \|A_m \|_2$ in the second inequality. The point of introducing $A^{(m)}$ is that for the second term in  \eqref{dec}, we have that $A_m$ and $(\phi^{(m)} -\bar{\phi})$ are now independent. Thus this term can be tackled with concentration results. We now handle both terms in \eqref{dec}.
Recall the definition of the event $\mathcal{E}$ in \eqref{event}.

{\it Claim 1: Under $\mathcal{E}$, $T_1=\| A \|_{2 \to \infty} \| \phi - \phi^{(m)} \|_2 = O(\sqrt{\log(n)} \| \phi \|_{\infty})$.} 

To prove this claim, assume that $\mathcal{E}$ takes place. To bound $\| u - \phi^{(m)}  \|_2$, we will view $A^{(m)}$ as a perturbation of $A$ and apply the Davis-Kahan Theorem (Theorem \ref{thm-DK}).

We first show that
\begin{align}
\mathcal{E} \subseteq \left\{ \| \phi^{(m)} - \phi \|_2 = \min_{s=\pm 1} \| s \phi^{(m)} - \phi \|_2 \right\}.
\label{ineq-loop-event-1}
\end{align}
Note that
\begin{align}
\| A^{(m)} - A \|_2 \leq \| A^{(m)} - A \|_F \leq \sqrt{2}\| A \|_{2\to\infty},
\label{ineq-loop-adj}
\end{align}
and
\begin{align}
\| A \|_{2\to\infty} \leq \| A  - \bar{A}\|_{2} + \| \bar{A}\|_{2\to\infty} \lesssim \sqrt{\log n} + \frac{ \log n}{\sqrt{n}} \lesssim \sqrt{\log n}.\label{ineq-row}
\end{align}
Therefore, (\ref{ineq-loop-adj}) and (\ref{ineq-row}) imply that
\begin{align}
\| A^{(m)} - \bar{A} \|_2 \leq \| A^{(m)} - A \|_2 + \| A - \bar{A} \|_2  \lesssim \sqrt{\log n}.
\label{ineq-loop-matrix}
\end{align}
By definition we have $(\bar{\phi})^T \phi \geq 0$ and $(\bar{\phi})^T \phi^{(m)} \geq 0$. Using Theorem \ref{thm-DK}, we have 
\begin{align}
&\| \phi^{(m)} - \bar{\phi} \|_2=
\min_{s=\pm 1}\| s \phi^{(m)} - \bar{\phi} \|_2 \lesssim \| A^{(m)} - \bar{A} \|_2 / \bar{\lambda} \lesssim 1/\sqrt{\log n},
\label{ineq-loop-L2}
\\
&\| \phi- \bar{\phi} \|_2=
\min_{s=\pm 1}\| s \phi - \bar{\phi} \|_2 \lesssim \| A - \bar{A} \|_2 / \bar{\lambda} \lesssim 1/\sqrt{\log n}.
\end{align}
and thus 
\begin{align}
\| \phi^{(m)} - \phi \|_2 \lesssim 1/\sqrt{\log n}.
\label{ineq-loop}
\end{align}
When $n$ is large enough, we have that $\| \phi^{(m)} - \phi \|_2 \leq 1$, which implies (\ref{ineq-loop-event-1}) since $\{\pm 1\} \ni s \mapsto \| s \phi^{(m)} - \phi \|_2 = 2- s\langle \phi^{(m)} , \phi \rangle$ has its minima  below 1 for $s=1$.  

By Weyl's inequality, $\max_i | \lambda_i(A) - \lambda_i(\bar{A}) | \le \| A - \bar{A} \|_2$. Recall that we are under $\mathcal{E}$, thus 
\begin{align}
\lambda_1(A) - \lambda_2(A) \geq \bar{\lambda} - 2 \| A-\bar{A} \|_2 \gtrsim \bar{\lambda} \gtrsim \log n.
\label{ineq-loop-gap}
\end{align}
Moreover from (\ref{ineq-loop-matrix}) we have that for $n$ large enough,
\begin{align}
\| A^{(m)} - A \|_2 < \lambda_1(A) - \lambda_2(A). \label{cond2}
\end{align}

Therefore (\ref{cond2}) satisfies the conditions for Theorem \ref{thm-DK}. Combined with (\ref{ineq-loop-event-1}), this yields
\begin{align}
\| \phi^{(m)} - \phi \|_2 &=
\min_{s=\pm 1} \| s \phi^{(m)} - \phi \|_2\lesssim \frac{ \|( A^{(m)} - A )\phi \|_2 }{\lambda_1(A) - \lambda_2(A)} \\
& \lesssim \frac{ \|( A - A^{(m)} )\phi \|_2 }{ \bar{\lambda} }.
\label{ineq-loop-DK}
\end{align}
Note that $(( A - A^{(m)} )\phi )_m =A_{m} \phi = \lambda \phi_m$ and
$(( A - A^{(m)} )\phi )_i = A_{im} \phi_m$ for $i\neq m$. 
By (\ref{order1}) and (\ref{ineq-row}),
\begin{align}
\| ( A - A^{(m)} )\phi\|_2 &= \left( \lambda^2 |\phi_m|^2 + \sum_{i\neq m} A_{im}^2 \phi_m^2 \right)^{1/2} \\
& \leq |\phi_m| \sqrt{\lambda^2 + \| A \|_{2\to\infty}^2 } \, \lesssim |\phi_m| \bar{\lambda} .
\end{align}
Using this with (\ref{ineq-loop-DK}), we have that there exists $C_1>0$ such that  
\begin{align}
 &\| \phi^{(m)} - \phi \|_2 \le C_1  |\phi_m| \le C_1  \| \phi \|_{\infty} , \quad \forall m \in [n].
\label{ineq-loop-event}
\end{align}
%\begin{align}
%\mathcal{E} \subseteq \bigcap_{m=1}^n \left \{ \| \phi^{(m)} - \phi \|_2 \leq C_1 |\phi_m| \right\}
%\subseteq \bigcap_{m=1}^n \left \{ \| \phi^{(m)} - \phi \|_2 \leq C_1 \| \phi \|_{\infty} \right\}.
%\label{ineq-loop-event}
%\end{align}
Finally, from (\ref{ineq-row}) and (\ref{ineq-loop-event}), there exists $C_2>0$ such that  
\begin{align}
T_1=\| A \|_{2 \to \infty} \| \phi - \phi^{(m)} \|_2 \le C_2 \sqrt{\log(n)} \| \phi \|_{\infty} ,
\label{ineq-loop-event-0}
\end{align}
which proves Claim 1.

{\it Claim 2: Under $\mathcal{E}$, $T_2 =| A_m (\phi^{(m)} -\bar{\phi}) |= O(\log(n)\| \phi \|_{\infty} / \log\log(n))$.} 

To prove this claim, we work again under $\mathcal{E}$ and use a concentration bound. This is where we exploit the independence between $\phi^{(m)} - \bar{\phi}$ and $\{ A_{mi} \}_{i=1}^n$ to control $|A_{m}(\phi^{(m)}-\bar{\phi})|$. From concentration bounds, namely taking $w = \phi^{(m)} - \bar{\phi}$, $\{ X_i \}_{i=1}^n=\{ {A'}_{mi} \}_{i=1}^n$ (note that $A_{mi}- \bar{A}_{mi} = {A'}_{mi}- \E {A'}_{mi}$), $p=(a\vee b) \frac{\log n}{n}$ and $\alpha=3/(a\vee b)$ in Lemma 10 from \cite{kaiz}, we get
\begin{align*}
	&\pp\left( \left| \sum_{i=1}^{n}( A_{mi}- \bar{A}_{mi} )( \phi^{(m)} - \bar{\phi} )_i \right| <
	\frac{[2(a\vee b) + 3]\log n\cdot \|\phi^{(m)} - \bar{\phi}\|_{\infty}}{ 1 \vee \log \left( \frac{\sqrt{n} \|\phi^{(m)} - \bar{\phi}\|_{\infty} }{\| \phi^{(m)} - \bar{\phi} \|_2}  \right) } 
	\right)  \\
	&> 1- 2n^{-3}.
\end{align*}
Had we applied Bernstein's inequality directly, we would get a loose upper bound that contains the term $\log n\cdot \|\phi^{(m)} - \bar{\phi}\|_{\infty}$ without the denominator, which turns out to be critical for the analysis. In fact, as we show below, the denominator will give us an additional $\log \log n$ factor.

For a scalar $C_3$, define 
\begin{align}
&\mathcal{E}^{(m)}_0:= \notag \\
& \{ \left| A_{m}( \phi^{(m)} - \bar{\phi} ) \right| < \frac{C_3 \log n}{\sqrt{n}}
(
\| \phi^{(m)} - \bar{\phi} \|_2 +
\frac{ \sqrt{n} \|\phi^{(m)} - \bar{\phi}\|_{\infty}}{ 1 \vee \log \left( \frac{\sqrt{n} \|\phi^{(m)} - \bar{\phi}\|_{\infty} }{\| \phi^{(m)} - \bar{\phi} \|_2}  \right) } \label{longbnd}
) \}, \\
&\mathcal{E}_0:= \bigcap_{m=1}^n \mathcal{E}^{(m)}_0.
\end{align} 
Since
\begin{align*}
\left| \sum_{i=1}^{n} \bar{A}_{mi}( \phi^{(m)} - \bar{\phi} )_i \right| \leq \| \bar{A} \|_{2\to\infty} \| \phi^{(m)} - \bar{\phi} \|_2
\lesssim  \frac{\log n}{\sqrt{n}} \| \phi^{(m)} - \bar{\phi} \|_2,
\end{align*}
there exists a choice of $C_3 >0$ in the definition of $\mathcal{E}^{(m)}_0$ such that
\begin{align}
\pp (\mathcal{E}^{(m)}_0) &> 1- 2n^{-3},
\end{align}
and from the union bound, 
\begin{align}
\pp  ( \mathcal{E}_0^c ) &\leq 2n^{-2}.
\end{align}
From \eqref{ineq-loop}, we have
\begin{align}
\mathcal{E} \subseteq \bigcap_{m=1}^n \left\{ 
\| \phi^{(m)} - \bar{\phi} \|_2 \leq C_4 /\sqrt{\log n}
\right\}
\label{ineq-concentration-event}
\end{align}
for some constant $C_4 \geq 1$. Define $\gamma = C_4/\sqrt{\log n}$, which is smaller than 1 for large enough $n$, and define function $h(z) = [1 \vee \log(1/z)]^{-1}$ where $z \in (0,1]$. It is easily checked that $h(z)$ is non-decreasing and $h(z)/z$ is non-increasing, so $h(z) \le h(\gamma) \vee h(\gamma)z/\gamma$. Setting $z = \frac{\| \phi^{(m)} - \bar{\phi} \|_2}{\sqrt{n} \|\phi^{(m)} - \bar{\phi}\|_{\infty} }$, we can use this inequality to simplify the bound in \eqref{longbnd} and obtain
%todo (see Lemma \ref{lem-elementary} in Appendix \ref{}) 
\begin{align*}
\mathcal{E} \cap \mathcal{E}^{(m)}_0
\subseteq \left\{ 
\left| A_{m}( \phi^{(m)} - \bar{\phi} ) \right| <2 C_3\log n \frac{  \|\phi^{(m)} - \bar{\phi}\|_{\infty} \vee \frac{1}{\sqrt{n}} }{\log (1/\gamma)} \right\}.
\end{align*}
Note that under $\mathcal{E}$, (\ref{ineq-loop-event}) leads to
\begin{align*}
\| \phi^{(m)} - \bar{\phi} \|_{\infty} & \leq \| \phi^{(m)} - \phi \|_{\infty} + \| \phi - \bar{\phi} \|_{\infty} \\ &\leq \| \phi^{(m)} - \phi \|_{2} + \| \phi \|_{\infty} + \| \bar{\phi}\|_{\infty} \notag\\
&\leq (C_1+1) \|\phi \|_{\infty} + \frac{1}{\sqrt{n}}
\leq (C_1+2) \|\phi \|_{\infty},
\end{align*}
where we use $\| \phi \|_{\infty} \geq 1/\sqrt{n}$. Hence, recalling that $\gamma = C_4/\sqrt{\log n}$, we have that under $\mathcal{E} \cap \mathcal{E}_0$ there exists $C_5>0$ such that
\begin{align}
  \left|  A_{m}( \phi^{(m)} - \bar{\phi} ) \right| < \frac{ C_5 \log(n) \| \phi \|_{\infty} }{\log \log n }, \quad \forall m \in [n],
\label{ineq-concentration-event-2}
\end{align}
which implies Claim 2.

Let us use now plug the bounds from Claim 1 and 2 in \eqref{dec}, which gives under $\mathcal{E} \cap \mathcal{E}_0$, 
\begin{align}
\frac{1}{ \bar{\lambda}} \| A(\phi - \bar{\phi}) \|_{\infty}  \lesssim \frac{\|\phi \|_{\infty}}{ \log\log(n)}.
\end{align}
Putting this back in \eqref{broken} together with \eqref{plug3}, we obtain under $\mathcal{E} \cap \mathcal{E}_1$,
\begin{align}
\| \phi - A\bar{\phi}/\bar{\lambda} \|_{\infty} \lesssim  \frac{\|\phi \|_{\infty}}{ \log\log(n)}.
\label{ineq-penultimate}
\end{align}
We are now left to show that $\| \phi \|_{\infty} = O ( 1/\sqrt{n})$, which implies the theorem's statement since $\mathcal{E} \cap \mathcal{E}_0$ takes place with probability $1- \Omega(n^{-2})$. 

{\it Claim 3: $\| \phi \|_{\infty} = O ( 1/\sqrt{n})$.} 

To prove this claim, note that from (\ref{ineq-penultimate}) and $\|\phi \|_{\infty} \leq \| \phi - A\bar{\phi}/\bar{\lambda} \|_{\infty} + \| A\bar{\phi}/\bar{\lambda}\|_{\infty}$, we have $\| \phi \|_{\infty} \lesssim \| A\bar{\phi}\|_{\infty}/\bar{\lambda}$. Hence it remains to show that $\| A\bar{\phi}\|_{\infty} \lesssim \frac{\log n}{\sqrt{n}}$. Observe that $|A_{m} \bar{\phi}|
\leq \| \bar{\phi} \|_{\infty} \sum_{i=1}^{n} |A_{mi}| = \sum_{i=1}^{n} |A_{mi}| /\sqrt{n}$ and
\begin{align*}
&|A_{mi}|=\left |{A'}_{mi} - \frac{a-b}{2n}\log n \right |
=\begin{cases}
& | 1- \frac{a-b}{2n}\log n|,\text{ if }{A'}_{mi}=1\\
& | \frac{a-b}{2n}\log n |,\text{ if }{A'}_{mi}=0
\end{cases} \\
& \leq 2A_{mi}^2 + \frac{|a-b|}{2n}\log n,
\end{align*}
where the last inequality holds for $n$ large enough. Moreover from (\ref{ineq-row}), 
\[
\sum_{i=1}^{n} |A_{mi}| \leq \sum_{i=1}^{n} 2A_{mi}^2 + \frac{|a-b|}{2}\log n
\leq 2\|A\|_{2\to\infty}^2 + \frac{|a-b|}{2}\log n \lesssim \log n.
\]
Therefore, $\left\| A \bar{\phi} \right\|_{\infty} \lesssim \frac{ \log n  }{\sqrt{n}}$, which proves Claim 3 and the theorem. \end{proof}

\subsection{SDP algorithm}\label{sdp2}
Recall that an SDP was derived in Section \ref{achieve-exact}, lifting the variables to change the quadratic optimization into a linear optimization, and removing the rank-one constraint to obtain 
%(as for max-cut \cite{MXGoemans_DPWilliamson_1995}). Since $\tr(AB)=\tr(BA)$ for any matrices of matching dimensions, we have  
%\begin{align}
% x^t A x  = \tr (x^t A x) = \tr (A xx^t), \label{minbi2}
%\end{align}
%hence defining $X:=xx^t$, we can write \eqref{minbi2} as 
%\begin{align}
%\hat{X}_{\map}(g) = \argmax_{  X \succeq 0 \atop{ X_{ii} = 1, \forall i \in [n] \atop{ \rank X =1 \atop{  X 1^n = 0 }}}} \tr(AX). \label{minbi6}
%\end{align}
%The first three constraints on $X$ force $X$ to take the form $xx^t$ for a vector $x \in \{+1,-1\}^n$, as desired, and the last constraint gives the balance requirement. The advantage of \eqref{minbi6} is that the objective function is linear in the lifted variable $X$. The constraint $\rank X=1$ is responsible for hardness of the optimization. We hence simply remove that constraint to obtain an SDP relaxation:
\begin{align}
\hat{X}_{sdp}(g) = \argmax_{  X \succeq 0 \atop{ X_{ii} = 1, \forall i \in [n]  \atop{  X 1^n = 0 }}} \tr(AX). \label{minbib2}
\end{align}
%A possible approach to handle the constraint $X 1^n = 0 $ is to replace the adjacency matrix $A$ by the matrix $B$ such that $B_{ij}=1$ if there is an edge between vertices $i$ and $j$, and $B_{ij}=-1$ otherwise. Using $-T$ for a large $T$ instead of $-1$ for non-edges would force the clusters to be balanced, and it turns out that $-1$ is already sufficient for our purpose. This gives another SDP:
or the version using the matrix $B$ such that $B_{ij}=1$ if there is an edge between vertices $i$ and $j$, and $B_{ij}=-1$ otherwise, 
\begin{align}
\hat{X}_{SDP}(g) = \argmax_{  X \succeq 0 \atop{ X_{ii} = 1, \forall i \in [n]  }} \tr(BX). \label{minbi2b2}
\end{align}
The dual of this SDP is given by 
\begin{align}
\min_{  Y_{ij}=0 \forall 1\le i \neq j \le n \atop{ Y \succeq B}} \tr(Y). \label{minbi32}
\end{align}
Since the dual minimization gives an upper-bound on the primal maximization, a solution is optimal if it makes the dual minimum match the primal maximum. The Ansatz here consists in taking $Y=2(D_{in}-D_{out})+I_n$ as a candidate for the diagonal matrix $Y$, which gives the primal maxima. It we thus have $Y \succeq B(g)$, this is a feasible solution for the dual, and we obtain a dual certificate. The following is shown in \cite{abh} based on this reasoning. 
\begin{definition}
Define the SBM Laplacian for $G$ drawn under the symmetric SBM with two communities by 
\begin{align}
L_{\mathrm{SBM}}= D(G_{in}) - D(G_{out}) - A,
\end{align}
where $D(G_{in})$ ($D(G_{out})$) are the degree matrices of the subgraphs of $G$ containing only the edges inside (respectively across) the clusters, and $A$ is the adjacency matrix of $G$.
\end{definition}
\begin{theorem}\cite{abh}
The SDP solves exact recovery in the symmetric SBM with 2 communities if $2L_{\mathrm{SBM}}+ 11^t+I_n \succeq 0$ and $\lambda_2(2L_{\mathrm{SBM}}+ 11^t+I_n) >0$. 
\end{theorem}
The second requirement above is used for the uniqueness of the solution.
This condition is verified all the way down to the exact recovery threshold. In \cite{abh}, it is shown that this condition holds in a regime that does not exactly match the threshold, off roughly by a factor of 2 for large degrees.  This gap is closed in \cite{new-xu,afonso_single}, which show that SDPs achieve the exact recovery threshold in the symmetric case. Results for unbalanced communities were also obtained in \cite{wein_sdp}, although it is still open to achieve the general CH threshold with SDPs. Many other works have studied SDPs for the stochastic block model, we refer to \cite{levina,abh,afonso_single,afonso3,new-xu,montanari_sen,wein_sdp} for further references. In particular, \cite{montanari_sen} shows that SDPs can approach the threshold for weak recovery in the two-community SSBM arbitrarily close when the expected degrees diverge. Recently, \cite{rigollet2} also studied an Ising model with block structure, obtaining results for exact recovery and SDPs.

\chapter{Weak recovery for two communities}
The focus on weak recovery, also called detection, was initiated\footnote{The earlier work \cite{Reichardt2008} also considers detection in the SBM.} with \cite{coja-sbm,decelle}. Note that weak recovery is typically investigated in SBMs where vertices have constant expected degree, as otherwise the problem can trivially be resolved by exploiting the degree variations. 

%\section{Fundamental limit and KS threshold}
The following conjecture was stated in \cite{decelle} based on deep but non-rigorous statistical physics arguments, and is responsible in part for the resurgent interest in the SBM: \\

\noindent
{\it 
{\bf Conjecture.} \cite{decelle,Mossel_SBM1}
Let $(X,G)$ be drawn from $\ssbm(n,k,a/n,b/n)$, i.e., the symmetric SBM with $k$ communities, probability $a/n$ inside the communities and $b/n$ across. Define $\snr=\frac{(a-b)^2}{k(a+(k-1)b)}$. Then, 
\begin{enumerate}
\item[(i)] For any $k \ge 2$, it is possible to solve weak recovery efficiently if and only if $\snr>1$ (the Kesten-Stigum (KS) threshold);
%\item If $k \in \{2,3,4\}$, and $\snr \leq 1$, it is impossible to detect communities information-theoretically,
\item[(ii)] If\footnote{The conjecture states that $k=5$ is necessary when imposing the constraint that $a>b$, but $k=4$ is enough in general.} $k \geq 4$, it is possible to solve weak recovery information-theoretically (i.e., not necessarily in polynomial time in $n$) for some $\snr$ strictly below 1.\footnote{\cite{decelle} made in fact a more precise conjecture, stating that there is a second transition below the KS threshold for information-theoretic methods when $k\ge 4$, whereas there is a single threshold when $k=3$. 
%\cite{decelle} provides also a more general conjecture for asymmetric SBMs, which was disproved in \cite{colin3cpam} due to mismatch in the definition of weak recovery.
}
\end{enumerate}
}

It was proved in \cite{massoulie-STOC,Mossel_SBM2} that the KS threshold can be achieved efficiently for $k=2$, and \cite{Mossel_SBM1} shows that it is impossible to detect below the KS threshold for $k=2$. Further, \cite{yash_sbm} extends the results for $k=2$ to the case where $a$ and $b$ diverge while maintaining the SNR finite. So weak recovery is closed for $k=2$ in SSBM. 
\begin{theorem}\label{mass_thm}[Converse in \cite{Mossel_SBM1}, achievability in \cite{massoulie-STOC,Mossel_SBM2}]
Weak recovery is solvable (and efficiently so) in $\ssbm(n,2,a/n,b/n)$ when $a,b=O(1)$ if and only if $(a-b)^2>2(a+b)$.
\end{theorem}
\begin{theorem}\label{yash_thm}\cite{yash_sbm}
Weak recovery is solvable (and efficiently so) in $\ssbm(n,2,a_n/n,b_n/n)$ when $a_n,b_n =\omega(1)$ and $(a_n-b_n)^2/(2(a_n+b_n)) \to \lambda $ if and only if $\lambda>1$.
\end{theorem}
Theorem \ref{yash_thm} is discussed in Section \ref{partial} in the context of partial recovery. Here we discuss Theorem \ref{mass_thm}. An additional result is obtained in \cite{Mossel_SBM1}, showing that when $\snr \le 1$, the symmetric SBM with two communities is in fact contiguous to the Erd\H{o}s-R\'enyi model with edge probability $(a+b)/(2n)$, i.e, distinguishability is not solvable in this case. Contiguity is further discussed in Section \ref{info}. 

For several communities, it was also shown in \cite{bordenave} that for SBMs with multiple communities that are balanced and that satisfy a certain asymmetry condition (i.e., the requirement that $\mu_k$ is a simple eigenvalue in Theorem 5 of \cite{bordenave}), the KS threshold can be achieved efficiently. The achievability  parts of previous conjecture for $k \geq 3$ are resolved in \cite{colin3,colin3cpam}. We discuss these in Section \ref{weak}.

Note that the terminology `KS threshold' comes from the reconstruction problem on trees \cite{ks1,evans,MosselRec,tree-andrea}, referring to the first paper by Kesten-Stigum (KS). A transmitter broadcasts a uniform bit to some relays, which themselves forward the received bits to other relays, etc.\ The goal is to reconstruct the root bit from the leaf bits as the depth of the tree diverges. In particular, for two communities, \cite{Mossel_SBM1} makes a reduction between failing in the reconstruction problem in the tree setting and failing in weak recovery in the SBM. This is discussed in more details in next section. The fact that the reconstruction problem on tree also gives the positive behavior for efficient algorithms requires a more involved argument, as discussed in Section \ref{weak_achieve}.

%The number of relays (or offspring) at each generation may be a constant $c$, or may be random such as Poisson distributed of mean $c$. Each relay is assumed to relay the bit by flipping the bit with probability $\e$ (i.e., binary symmetric channels of parameter $\e$). The receiver gets to see all the bits at the leaves. For what values of $c$ and $\e$ could the receiver reconstruct the original bit when the tree depth diverges? The goal is to recover the root bit weakly, i.e., with probability away from $1/2$, and not tending to 1 as usual in information theory. This problem was first solved in \cite{ks1} for binary symmetric channels and constant offspring, showing that weak recovery is possible if and only if  $c > 1/(1-2\e)^2$, which became the KS threshold. It was later solved for more general offsprings, such as the Poisson case, in \cite{evans}. This implies a converse for weak recovery in the 2-community SBM as shown in \cite{Mossel_SBM1}, using a genie-aided argument and the fact that a node's neighborhood in the sparse SBM is tree-like. In this context, we have $c=(a+b)/2$, $\e=b/(a+b)$, and the KS threshold reads $(a-b)^2 > 2(a+b)$ as stated in Conjecture 1. 

Achieving the KS threshold raises an interesting challenge for community detection algorithms, as standard clustering methods fail to achieve the threshold, as discussed in Section \ref{linBP}.
This includes spectral methods based on the adjacency matrix or standard Laplacians, as well as SDPs. 
For standard spectral methods, a first issue is that the fluctuations in the node degrees produce high-degree nodes that disrupt the eigenvectors from concentrating on the clusters. A classical trick is to suppress such high-degree nodes, by either trimming or shifting the matrix entries \cite{joseph,levina-reg,coja-sbm,Vu-arxiv,sbm-groth,new-vu}, throwing away some information, but this does not suffice to achieve the KS threshold \cite{kawamoto2015limitations}. SDPs are a natural alternative, but they also appear to stumble before the KS threshold \cite{sbm-groth,montanari_sen,ankur_SBM}, focusing on the most likely rather than typical clusterings. As shown in \cite{massoulie-STOC,Mossel_SBM2,bordenave,colin3}, approximate BP algorithms or spectral algorithms on more robust graph operators instead allow us to achieve the KS threshold.

\section{Warm up: broadcasting on trees}\label{bot}
As for exact recovery, we start with a simpler problem that will play a key role in understanding weak recovery. The idea is similar to that of the exact recovery warm up, except that we do not reveal only the direct neigbors of a vertex, but the neighbors at a small depth. In particular, at depth $(1/2-\e)\log(n)/\log((a+b)/2)$, the SBM neigborhood of a vertex can be coupled with a Galton-Watson tree, and so we will consider trees for the warm up. In contrast to exact recovery, we will not be interested in reconstructing the isolated vertex with probability tending to 1, but with probability greater than $1/2$. This is known in the literature as the reconstruction problem for broadcasting on trees. We refer to \cite{MosselRec} for a survey on this topic.
 
The problem consists of broadcasting a bit from the root of a tree down to its leaves, and trying to guess back this bit from the leaf bits at large depth. First consider the case of a deterministic tree with fixed degree $c+1$, i.e., each vertex has exactly $c$ descendants (note that the root has degree $c$). Assume that on each branch of the tree the incoming bit is flipped with probability $\e \in [0,1]$, and that each branch flips independently. Let $X^{(t)}$ be the bits received at depth $t$ in this tree, with $X^{(0)}$ being the root bit, assumed to be drawn uniformly at random in $\{0,1\}$. 

We now define weak recovery in this context. Note that $\E(X^{(0)}|X^{(t)})$ is a random variable that gives the probability that $X^{(0)}=1$ given the leaf bits, as a function of the leaf bits $X^{(t)}$. If this probability is equal to $1/2$, then the leaf bits provide no useful information about the root, and we are interested in understanding whether this takes place in the limit of large $t$ or not.
\begin{definition}
Weak recovery (also called reconstruction) is solvable in broadcasting on a regular tree if $\lim_{t \to \infty} \E |\E (X^{(0)}|X^{(t)}) - 1/2| >0$. Equivalently, weak recovery is solvable if $\lim_{t \to \infty} I(X^{(0)};X^{(t)}) >0$, where $I$ is the mutual information. 
\end{definition}
Note that the above limits exist due to monotonicity arguments. The first result on this model is due to Kesten-Stigum:
\begin{theorem}\label{tree-ks}
In the tree model with constant degree $c$ and flip probability $\e$,  
\begin{itemize}
\item \cite{ks1} weak recovery is solvable if $c(1-2\e)^2 >1$, 
\item \cite{ble,evans} weak recovery is not solvable\footnote{The proof from \cite{evans} appeared first in 1996.} if $c(1-2\e)^2 \le 1$.  
\end{itemize}
\end{theorem}
In fact, one can show a seemingly stronger result where weak recovery fails in the sense that, when $c(1-2\e)^2 \le 1$, 
\begin{align}
\lim_{t \to \infty} \E(X^{(0)}|X^{(t)}) = 1/2 \quad \text{a.s.}
\end{align}
Thus weak recovery in the tree model is solvable if and only if $c(1-2\e)^2 >1$, which gives rise to the so-called Kesten-Stigum (KS) threshold in this tree context. Note that  \cite{MosselRec} further shows that the KS threshold is sharp for ``census reconstruction,'' i.e., deciding about the root-bit by taking majority on the leaf-bits, which is shown to still hold in models such as the multicolor Potts model where the KS threshold is no longer sharp for reconstruction. 

%This result can be guessed by computing the first two moments of the random variable that counts the number of 1's minus 0's at depth $t$. For $t=1$, the expectation is $c(1-2\e)$ and the variance is $2c\e(1-\e)$. Another point of view is obtained from the mutual information. 
To see this result, let us compute the moments of the number of 0-bits minus 1-bits at generation $t$. First, consider $\pm 1$ variables rather than bits, i.e., redefine $X_i^{(t)} \leftarrow (-1)^{X_i^{(t)}}$, and consider the difference variable:
\begin{align}
\Delta^{(t)} = \sum_{i \in [c^t]} X_i^{(t)}.
\end{align}
Note that, if there are $x$ bits of value $1$ and $y$ bits of value $-1$ at generation $t$, then $\Delta^{(t+1)}$ would be the sum of $xc$ Radamacher$(1-\e)$ and $yc$ Radamacher$(\e)$ (all independent), and since the expectation of Radamacher$(1-\e)$ is $1-2\e$, 
%and the variance of Radamacher$(1-\e)$ or Radamacher$(\e)$ is $2\e(1-2\e)$, 
we have 
\begin{align}
\E (\Delta^{(t+1)} | \Delta^{(t)})  &= c (1-2\e) \Delta^{(t)}.%\\
%\Var \Delta^{(t+1)} =  \Var (\Delta^{(t+1)}|X^{(0)})  &= c 2 \e(1-\e) c^t.
\end{align}
Since $X^{(0)} - \Delta^{(t)} - \Delta^{(t+1)}$ forms a Markov chain and $\E (\Delta^{(0)}|X^{(0)})  = X^{(0)}$,
\begin{align}
\E( \Delta^{(t+1)}|X^{(0)}) &= \E( \E ( \Delta^{(t+1)} |\Delta^{(t)}) |X^{(0)}) \\
&=c (1-2\e)\E( \Delta^{(t)}|X^{(0)})  \\
&= c^{t+1} (1-2\e)^{t+1} X^{(0)}.
\end{align}
We now look at the second moment. A direct computation gives 
\begin{align}
\E ((\Delta^{(t+1)})^2 | \Delta^{(t)})  &= c^{t+1} 4 \e (1-\e) + c^2 (1-2\e)^2 (\Delta^{(t)})^2.%\\
%\Var \Delta^{(t+1)} =  \Var (\Delta^{(t+1)}|X^{(0)})  &= c 2 \e(1-\e) c^t.
\end{align}
Defining the signal-to-noise ratio as $\snr= \sqrt{c}(1-2\e)$ and assuming that $\snr>1 $, we have by iteration 
\begin{align}
\E( (\Delta^{(t+1)})^2 |X^{(0)}) &= c^{2(t+1)} (1-2\e)^{2(t+1)} (\frac{4 \e (1-\e)}{ c (1-2\e)^2-1} + 1)(1+o_t(1))
\end{align}
and
\begin{align}
\Var( \Delta^{(t+1)} |X^{(0)}) &= c^{2(t+1)} (1-2\e)^{2(t+1)} \frac{4 \e (1-\e)}{c (1-2\e)^2-1} (1+o_t(1)).
\end{align}
%Therefore, we have 
%\begin{align}
%(\E( \Delta^{(t)}|X^{(0)}= 1), \Var( \Delta^{(t)}|X^{(0)}= 1)) &\asymp (c^t(1-2\e)^t,c^t),\\
%(\E( \Delta^{(t)}|X^{(0)}= -1), \Var( \Delta^{(t)}|X^{(0)}= -1)) &\asymp (-c^t(1-2\e)^t,c^t).
%\end{align}
Denote now by $\mu_+$ the distribution of $\Delta^{(t)}$ given that $X^{(0)}= 1$, and by $\mu_-$ the distribution of $\Delta^{(t)}$ given that $X^{(0)}= -1$. We have by Cauchy-Schwarz, 
\begin{align}
&\E( \Delta^{(t)}|X^{(0)}= 1) - \E( \Delta^{(t)}|X^{(0)}= -1) \\
&= \sum_\delta \delta (\mu_+(\delta) - \mu_-(\delta))\\
&\le \sqrt{ \sum_\delta \frac{(\mu_+(\delta) - \mu_-(\delta))^2}{\mu_+(\delta) + \mu_-(\delta)} } \sqrt{ \sum_\delta \delta^2 (\mu_+(\delta) + \mu_-(\delta)) }\\
&=\sqrt{ \sum_\delta \frac{(\mu_+(\delta) - \mu_-(\delta))^2}{\mu_+(\delta) + \mu_-(\delta)} } \sqrt{ \E( (\Delta^{(t)})^2|X^{(0)}= 1) + \E( (\Delta^{(t)})^2|X^{(0)}= -1) }
\end{align} 
and (trivially) 
\begin{align}
\sum_\delta \frac{(\mu_+(\delta) - \mu_-(\delta))^2}{\mu_+(\delta) + \mu_-(\delta)} \le \| \mu_+ - \mu_- \|_1  .
\end{align}
Therefore, we have from previous expansions that the total variation distance between $\mu_+$ and $\mu_-$ is $\Omega(1)$, which implies that it is possible to distinguish between the hypotheses $X^{(0)}= 1$ and $X^{(0)}= -1$ with a probability of error that is $1/2-\Omega(1)$.
%we have that the statistics $\Delta^{(t)}$ will preserve information about $X^{(0)}$, in that the standard deviation of each posterior will not overflow the magnitude of their respective means, if $(\sqrt{c}(1-2\e))^t$ gets amplified as $t$ grows, i.e., if $\snr>1$.

Conversely, irrespective of the statistics used, one can show that the mutual information $I(X^{(0)};X^{(t)})$ vanishes as $t$ grows if $\snr < 1$. 
In the binary case, it turns out that the mutual information is subadditive among leaves \cite{evans}, i.e.,
\begin{align}
I(X^{(0)};X^{(t)}) \le \sum_{i=1}^{c^t} I(X^{(0)};X_i^{(t)})= c^t I(X^{(0)};X_1^{(t)}). \label{mmi-sa}
\end{align}
Note that this subadditivity holds in greater generality for the first layer, i.e., if we have a Markov chain $Y_1-X-Y_2$, such as happens when a root variable $X$ is broadcast on two independent channels producing $Y_1$ and $Y_2$, then
\begin{align}
&I(X;Y_1) + I(X;Y_2) -I(X;Y_1,Y_2)\\
%&=H(Y_1) - H(Y_1|X) + H(Y_2)- H(Y_2|X) - H(Y_1,Y_2) + H(Y_1,Y_2|X) \\
&=H(Y_1) - H(Y_1|X) + H(Y_2)- H(Y_2|X) \\
&- H(Y_1,Y_2) + H(Y_1|X) + H(Y_2|X) \\
&=H(Y_1) + H(Y_2) - H(Y_1,Y_2)\\
&=I(Y_1;Y_2) \ge 0.
\end{align}
However, going to depth 2 of the tree, there is no simple inequality as the above that shows the subadditivity, and in fact, the subadditivity is not true in general for binary non-symmetric noise or for non-binary labels. 
For binary labels and symmetric channels, it is shown in \cite{evans} that the distribution of labels on the tree can be obtained by degradation from a so-called ``stringy tree'', where branches are ``detached'', which implies the inequality.

Further, the channel between $X^{(0)}$ and a one leaf-bit such as $X_1^{(t)}$ corresponds to $t$ Bernoulli$(\e)$ random variables added, and the mutual information scales as $(1-2\e)^{2t}$, 
\begin{align}
I(X^{(0)};X_1^{(t)}) = O((1-2\e)^{2t})
\end{align}
which implies with the subadditivity
\begin{align}
I(X^{(0)};X^{(t)}) = O(c^t(1-2\e)^{2t}).
\end{align}
Therefore, if $c(1-2\e)^2<1$, 
$$I(X^{(0)};X^{(t)}) = o(1)$$
and the information of the root-bit gets lost in the leaves. 

The subadditivity property in \eqref{mmi-sa} can also be established for the Chi-squared mutual information, i.e., using $I_2(X;Y)=D_{\chi_2}(p_{X,Y}\| p_X p_Y)$ where $D_{\chi_2}$ is the Chi-squared $f$-divergence with $f(z)=(z-1)^2$. We describe here how this can be obtained with an induction on the tree depth, assuming the following two building blocks:
\begin{align}
&I_2(X^{(0)};Y_1,Y_2) \le I_2(X^{(0)};Y_1) + I_2(X^{(0)};Y_2) \label{ineqi1}
\end{align}
if $\,\,\, Y_1-X^{(0)}-Y_2$ (i.e., $Y_1,Y_2$ are independent conditionally on $X^{(0)}$ and $X^{(0)}$ is Bernoulli$(1/2)$ as before) and 
\begin{align}
&I_2(X^{(0)};Y_t) = I_2(X^{(0)};Y_1) I_2(Y_1;Y_t) \label{eqi1} 
\end{align}
if $Y_1$ is a direct descendant of $X^{(0)}$ and $Y_t$ are descendants of $Y_1$ at depth $t$ (in the broadcasting on trees model). One then obtains the subadditivity by applying inequality \eqref{ineqi1} to the sub-trees pending at the root, and equality \eqref{eqi1} to factor out the first edge on each sub-tree, reducing the depth by one and allowing for the induction hypothesis.    
Interestingly, while the inequality also holds for the mutual information, the equality does not hold for the mutual information, and can in fact go in the wrong direction. %Therefore such an induction proof does not work for the mutual information. 
On the other hand, one has 
\begin{align}
I_2(X^{(0)};X_1^{(t)}) =(1-2\e)^{2t},
\end{align} 
and since the Chi-squared mutual information upper-bounds the classical mutual information for binary inputs, 
the same threshold of $c(1-2\e)^2=1$ is obtained with this argument.

We will soon turn to the connection between the reconstruction on trees problem and weak recovery in the SBM. It is easy to guess that the tree will not be a fixed degree tree for us, but the local neighborhood of a vertex in the SBM, which behaves like a Galton-Watson tree of Poisson offspring. We first state the above results for Galton-Watson trees.

\begin{definition}
A Galton-Watson tree with offspring distribution $\mu$ on $\mZ_+$ is a rooted tree where the number of descendants from each vertex is independently drawn under the distribution $\mu$. We denote by $T^{(t)} \sim GW(\mu)$ a Galton-Watson tree with offspring $\mu$ and $t$ generations of descendants, where $T^{(0)}$ is the root vertex. 
\end{definition}
Define as before $X^{(t)}$ as the variables at generation $t$ obtained from broadcasting the root-bit on a Galton-Watson tree $T^{(t)}$.
\begin{definition}
Weak recovery is solvable in broadcasting on a Galton-Watson tree $\{T^{(t)}\}_{t \ge 0}$ if $\lim_{t \to \infty} \E |\E (X^{(0)}|X^{(t)},T^{(t)}) - 1/2| >0$. Equivalently, weak recovery is solvable if $\lim_{t \to \infty} I(X^{(0)};X^{(t)}| T^{(t)}) >0$, where $I$ is the mutual information. 
\end{definition}
In \cite{evans}, it is shown that the threshold $c(1-2\e)^2 >0$ is necessary and sufficient for weak recovery for a large class of offspring distributions, where $c$ is the expectation of $\mu$. The case of $\mu$ being the Poisson$(c)$ distribution is of particular interest to us:
\begin{theorem}\cite{evans}\label{tree-ks-2}
In the broadcasting model with a Galton-Watson tree of Poisson$(c)$ offspring and flip probability $\e$,  
weak recovery is solvable if and only if $$c(1-2\e)^2 >1.$$ 
\end{theorem}
Another important extension is the `robust reconstruction' problem \cite{robust_rec}, where the leaves are not revealed exactly but with the addition of independent noise. It was shown in \cite{robust_rec} that for very noisy leaves, the KS threshold is also tight.

\section{The information-theoretic threshold}
In this Section, we discuss the proof for the threshold of Theorem \ref{mass_thm}, i.e., 
that weak recovery is solvable in $\ssbm(n,2,a/n,b/n)$ (and efficiently so) if and only if $$(a-b)^2 > 2(a+b).$$
%\begin{itemize}
%\item \cite{Mossel_SBM1} weak recovery is not solvable if $(a-b)^2 \le 2(a+b)$.
%\item \cite{massoulie,Mossel_SBM2} weak recovery is solvable (and efficiently so) if $(a-b)^2 > 2(a+b)$.
%\end{itemize}
%\end{theorem}
%This result is obtained by a reduction to the problem of reconstruction on trees, which is next discussed; we refer to \cite{MosselRec} for a survey on this problem. 
We focus in particular on the information-theoretic converse. Interestingly, for the achievability part, there is not currently a simpler proof available in the literature than the proof that directly gives an efficient algorithm. Below we discuss some attempts, and in Section \ref{weak_achieve} we will the efficient achievability.

\subsection{Converse}\label{bot-conv}
We will now prove the converse using two important technical results proved in \cite{Mossel_SBM1}: (1) the coupling of a neighborhood of the SBM with the broadcasting on Galton-Watson tree, (2) the weak effect of non-edges. The second point refers to the fact that the absence of an edge between two vertices does not make their probability of being in the same community exactly half; in fact, $\pp\{X_1=X_2 | E_{1,2}=0\}= \frac{1-a/n}{1-a/n+1-b/n} = 1/2(1 + (b-a)/2n) + o(1/n)$, and there is a slight repulsion towards being in the same community.

The following result is shown in \cite{Mossel_SBM1}:
\begin{theorem}\cite{Mossel_SBM1}
Let $(X,G) \sim \ssbm(n,2,a/n,b/n)$, the SBM with two symmetric communities. If $(a-b)^2 \le 2(a+b)$, 
\begin{align}
\pp \{ X_1=1 | G,X_2=1 \} \to 1/2  \text{ a.a.s.}
\end{align}
\end{theorem}
Here we show the following equivalent result that also implies the impossibility of weak recovery:
\begin{theorem} \label{conv-w}\cite{Mossel_SBM1}
Let $(X,G) \sim \ssbm(n,2,a/n,b/n)$ with $(a-b)^2 \le 2(a+b)$. Then, 
\begin{align}
I(X_1;X_2|G) = o(1).
\end{align}
\end{theorem}

Note that the role of vertex 1 and 2 is arbitrary above, it could be any fixed pair of vertices (not chosen based on the graph). 
%\begin{corollary}\label{conv-w2}\cite{Mossel_SBM1}
%Weak recovery is not solvable in $\ssbm(n,2,a/n,b/n)$ if $(a-b)^2 \le 2(a+b)$.
%%\end{itemize}
%\end{corollary}

The connection between the warm-up problem and the SBM comes from the fact that if one picks a vertex $v$ in the SBM graph, its neighborhood at small enough depth behaves likes a Galton-Watson tree of offspring Poisson$(c)$, $c=((a+b)/2)$, and the labelling on the vertices behaves like the broadcasting process discussed above with a flip probability of $\e=b/(a+b)$. Note that the latter parameter is precisely the probability that two vertices have different labels given that there is an edge between them. More formally, if the depth is $t\le (1/2-\delta)\log(n)/\log(a+b)/2)$ for some $\delta>0$, then the true distribution and the above have a vanishing total variation when $n$ diverges. This depth requirement can be understood from the fact that the expected number of leaves in that case is in expectation $n^{1/2-\delta}$, and by the birthday paradox, no collision will likely occur between two vertices' neighborhoods if $\delta>0$ (hence no loops take place and the neighborhood is a tree with high probability).  

To establish the converse of Theorem \ref{mass_thm}, it is sufficient to argue that, if it is impossible to weakly recover a single vertex when a genie reveals all the leaves at such a depth, it must be impossible to solve weak recovery. In fact, consider $\pp\{X_u = x_u | G=g,X_v=x_v\}$, the posterior distribution given the graph and an arbitrary vertex revealed (here $u$ and $v$ are arbitrary and chosen before the graph is drawn). With high probability, these vertices will not be within small graph-distance of each other (e.g., at distance $\Omega(\log(n))$ with high probability), and one can open a small neighborhood around $u$ of diverging byt small enough depth (e.g., $\e \log(n)$ for a small enough $\e$.) Now reveal not only the value of $X_v$ but in fact all the values at the boundary of this neighborhood. This is an easier problem since the neighborhood is a tree with high probability and since there is approximately a Markov relationship between these boundary vertices and the original $X_v$ (note that `approximately' is used here since there is a negligible effect due to non-edges). We are now back to the broadcasting problem on trees discussed above, and the requirement $c(1-2\e)^2 \le 0$ gives the theorem's bound (since $c=(a+b)/2$ and $\e=b/(a+b)$). 

The reduction extends to more than two communities (i.e., non-binary labels broadcasted on trees) and to asymmetrical communities, but the tightness of the KS bound is no longer present in these cases. For two asymmetrical communities, the result still extends if the communities are roughly symmetrical, using \cite{roughly_sym} and \cite{mossel2,private_mossel}. 
For more than three symmetric communities, new gap phenomena take place \cite{colin3cpam}; see Section \ref{info}. 

We now proceed to proving Theorem \ref{conv-w}.
The first step is to formalize what is meant by the fact that the neighborhoods of the SBM look like a Galton-Watson tree. 
Let $G$
\begin{lemma}\label{coupling} \cite{Mossel_SBM1}
Let $(X,G) \sim \ssbm(n,2,a/n,b/n)$ and $R=R(n)=\lfloor \frac{1}{10} \log(n)/\log(2(a+b)) \rfloor$. 
Let $B_R:=\{v \in [n]: d_G(1,v)\le R  \}$ be the set of vertices at graph distance at most $R$ from vertex $1$, $G_R$ be the restriction of $G$ on $B_R$, and let $X_{R}=\{X_u : u \in B_R\}$.
Let $T_R$ be a Galton-Watson tree with offspring Poisson$(a+b)/2$ and $R$ generations, and let $\tilde{X}^{(t)}$ be the labelling on the vertices at generation $t$ obtained by broadcasting the bit $\tilde{X}^{(0)}:=X_1$ from the root with flip probability $b/(a+b)$. Let $\tilde{X}_{R}=\{\tilde{X}^{(t)}_u : t \le R\}$. Then, there exists a coupling between $(G_R,X_{R})$ and $(T_R,\tilde{X}_{R})$ such that 
\begin{align}
(G_R,X_{R}) = (T_R,\tilde{X}_{R}) \text{ a.a.s.}
\end{align}
\end{lemma}
The second technical lemma that we need regards the negligible effect of non-edges outside from our local neighborhood of a vertex. The difficulty here is that these non-edges are negligible if the vertices labels are more or less balanced; for example, the effect of non-edges would not be negligible if the vertices had all the same labels. 
\begin{lemma}\cite{Mossel_SBM1}\label{indep-l}
Let $(X,G) \sim \ssbm(n,2,a/n,b/n)$, $R=R(n)=\lfloor \frac{1}{10} \log(n)/\log(2(a+b)) \rfloor$ and $X_{\partial R}=\{X_u : d_G(u,1)=R\}$. Then, 
\begin{align}
\pp\{X_1=1 | X_{\partial R}, X_2,G\} = (1+o(1)) \pp\{X_1=1 | X_{\partial R}, G_R\} \label{indep}
\end{align}
for a.a.e.\ $(X,G)$.
\end{lemma}
The statement means that the probability that $(X,G)$ satisfies \eqref{indep} tends to one as $n$ tends to infinity. 
Note that $R$ could actually be $c\log(n)/\log((a+b)/2) \rfloor$ for any $c <1/2$; we keep the same $R$ as in the previous lemma to reduce the number of parameters. 
\begin{corollary}
Using the same definition as in Lemma \ref{indep-l}, 
\begin{align}
H(X_1 | X_{\partial R}, G, X_2 )=H(X_1 | X_{\partial R}, G_R) +o(1).
\end{align}
\end{corollary}
We can now prove Theorem \ref{conv-w}.
\begin{proof}[Proof of Theorem \ref{conv-w}.]
Let $T_R$ and $\{\tilde{X}^{(t)}\}_{t=0}^R$ be the random variables appearing in the coupling of Lemma \ref{coupling}. We have,
\begin{align}
1\ge H(X_1| G, X_2) &\ge  H(X_1| X_{\partial R}, G, X_2) \label{l1}\\
&= H(X_1 | X_{\partial R}, G_R) +o(1) \label{l2} \\
&=H(\tilde{X}^{(1)} | \tilde{X}^{(R)}, T_R) + o(1) \label{l3} \\
&= 1 + o(1) \label{l4}
\end{align}
where \eqref{l1} follows from the fact that conditioning reduces entropy, \eqref{l2} follows from the asymptotic conditional independence of Lemma \ref{indep-l}, \eqref{l3} from the fact that the neighbordhood of a vertex can be coupled with the broadcasting on Galton-Watson tree, i.e., Lemma \ref{coupling}, and \eqref{l4} follows from the fact that below the KS threshold weak recovery is not solvable for the broadcasting on trees problem, i.e., Theorem \ref{tree-ks-2}. Therefore,   
\begin{align}
I(X_1;X_2|G) =1- H(X_1| G, X_2) = o(1).
\end{align}
\end{proof}

%The number of relays (or offspring) at each generation may be a constant $c$, or may be random such as Poisson distributed of mean $c$. Each relay is assumed to relay the bit by flipping the bit with probability $\e$ (i.e., binary symmetric channels of parameter $\e$). The receiver gets to see all the bits at the leaves. For what values of $c$ and $\e$ could the receiver reconstruct the original bit when the tree depth diverges? The goal is to recover the root bit weakly, i.e., with probability away from $1/2$, and not tending to 1 as usual in information theory. This problem was first solved in \cite{ks1} for binary symmetric channels and constant offspring, showing that weak recovery is possible if and only if  $c > 1/(1-2\e)^2$, which became the KS threshold. It was later solved for more general offsprings, such as the Poisson case, in \cite{evans}. This implies a converse for weak recovery in the 2-community SBM as shown in \cite{Mossel_SBM1}, using a genie-aided argument and the fact that a node's neighborhood in the sparse SBM is tree-like. In this context, we have $c=(a+b)/2$, $\e=b/(a+b)$, and the KS threshold reads $(a-b)^2 > 2(a+b)$ as stated in Conjecture 1. 

\subsection{Achievability}\label{weak-it}
Interestingly, the achievability part of Theorem \ref{mass_thm} was directly proved for an efficient algorithm, as discussed in the next section. 
Using an efficient algorithm should a priori require more work than what could be achieved if complexity considerations were put aside, but a short information-theoretic proof has not been given in the literature yet. 
Here we sketch how this could potentially be achieved, although there may be simpler ways.
In Section \ref{info}, we discuss an alternative approach that is simpler than the approach described below; however, it does not provide the right constant for two communities.
%Let us first recall the result:
%\begin{theorem}
%Weak recovery is solvable in $\ssbm(n,2,a/n,b/n)$ if $(a-b)^2 > 2(a+b)$.
%\end{theorem}

Since $(a-b)^2 > 2(a+b)$ is the threshold for weak recovery in the broadcasting on trees problem (when the expected degree is $(a+b)/2$ and the flip probability is $b/(a+b)$), one would hope to find a proof that reduces weak recovery in the SBM to this broadcasting on trees problem. For a converse argument, it is fairly easy to connect to this problem since one can always use a genie that gives further information in a converse argument, in this case, the values at the boundary of a vertex' neighborhood. How would one connect to the broadcasting on trees problem for an achievability result? 

We will next discuss how one can hope to replace the genie by {\it many} random guesses. 
First consider the effect of making a random guess about each vertex. 
Let $(X,G) \sim \ssbm(n,2,a/n,b/n)$ and let $X(\e)=\{X_v + \mathrm{Ber}_v(1/2-\e): v \in [n]\}$ be the noisy genie, i.e., a corruption of each community labels with independent additive Bernoulli noise for each vertex with flip probability $1/2- \e$, with $\e \in [0,1/2]$.
Note that $X(0)$ is pure noise, so we can assume that we have access to $(X(0),G)$ rather than $G$ only (which may seem irrelevant). Next we argue that we can replace $(X(0),G)$ with $(X(\Theta(1/\sqrt{n})),G)$, i.e., not a purely noisy genie but a genie with a very small bias of order $\Theta(1/\sqrt{n})$ towards the truth. 
\begin{lemma}
If weak recovery is solvable by observing $(X(1/\sqrt{n}),G)$, then it is solvable by observing $G$ only.  
\end{lemma}
The proof follows by noting that if weak recovery is solvable using $(X(0),G)$, then it is solvable using $G$ only since $X(0)$ is independent of $(X,G)$. But with high probability $X(0)$ produces a bias of $\Theta(\sqrt{n})$ vertices on either the good or bad side (by the Central Limit Theorem); since both are equivalent in view of weak recovery, we can assume that we have a genie that gives $\Theta(\sqrt{n})$ vertices on the good side.  

{\it Naive plan:} as one can obtain a weak genie `for free' by random guessing, one may hope to connect to the broadcasting problem on trees by amplifying this weak genie for each vertex at tree-like depth. That is, take a vertex $v$ in the graph, open a neighborhood at depth $R(n)$ as in the converse argument of the previous section, and re-decide for the vertex $v$ by solving the broadcasting on trees problem using the noisy vertex labels at the leaves. Do this for each vertex in parallel; assuming that correlations between different vertices are negligible. 

We next explain why this plan is doomed to fail, because the depth $R(n)$ is too small to amplify such a weak genie. In fact, For a vertex $v$ and integer $t$, let $N_t(v)$ be the number of vertices $t$ edges away from $v$, $\Delta_t(v)$ be the difference between the number of vertices $t$ edges away from $v$ that are in community $1$ and the number of vertices $t$ edges away from $v$ that are in community $2$, and $\widetilde{\Delta}_t(v)$ be the difference between the number of vertices $t$ edges away from $v$ that are in $C_1$ and the number of vertices $t$ edges away from $v$ that are in $C_2$. For small $t$,

\[E[N_t(v)] \asymp \left(\frac{a+b}{2}\right)^t\]
and 
\[E[\Delta_t(v)] \asymp \left(\frac{a-b}{2}\right)^t\cdot (-1)^{X_v}.\]

For any fixed values of $N_t(v)$ and $\Delta_t(v)$, the probability distribution of $\widetilde{\Delta}_t(v)$ is essentially a Gaussian distribution with a mean of $\Theta(\Delta_t(v)/\sqrt{n})$ and a variance of $\approx N_t(v)$ because it is the sum of $N_t(v)$ nearly independent variables that are approximately equally likely to be $1$ or $-1$. So, $\widetilde{\Delta}_t(v)$ is positive with a probability of $\frac{1}{2}+\Theta(\Delta_t(v)/\sqrt{|N_t(v)|n})$. In other words, if $v$ is in community $1$ then $\widetilde{\Delta}_t(v)$ is positive with a probability of 
\[\frac{1}{2}-\Theta\left(\left( \frac{a-b}{2}\right)^t\cdot\left(\frac{a+b}{2}\right)^{-t/2}\frac{1}{\sqrt{n}}\right)\] 
and if $v$ is in community $2$ then $\widetilde{\Delta}_t(v)$ is positive with a probability of 
\[\frac{1}{2}+\Theta\left(\left( \frac{a-b}{2}\right)^t\cdot\left(\frac{a+b}{2}\right)^{-t/2}\frac{1}{\sqrt{n}} \right).\] 
If $(a-b)^2\le 2(a+b)$, then this is not improving the accuracy of the classification, so this technique is useless. On the other hand, if $(a-b)^2>2(a+b)$, the classification becomes more accurate as $t$ increases. However, this formula says that to classify vertices with an accuracy of $1/2+\Omega(1)$, we would need to have $t$ such that
\[\left( \frac{a-b}{2}\right)^{2t}=\Omega\left( \left(\frac{a+b}{2}\right)^{t}n\right).\]
However, unless
%\footnote{If $a=0$ and $k=2$ or $b=0$, then classifying vertices based on the sign of $\widetilde{\Delta}_t(v)$ for suitable $t$ is likely to work, but this is pointlessly complicated because every component of the graph consists of all vertices of one community or vertices of alternating communities.} 
$a$ or $b$ is $0$, that would imply that
\begin{align}
\left( \frac{a+b}{2}\right)^{2t}=\omega\left(\left( \frac{a-b}{2}\right)^{2t}\right)=\omega\left( \left(\frac{a+b}{2}\right)^{t}n\right) \label{needed}
\end{align}
which means that $(\frac{a+b}{2})^{t}=\omega(n)$. It is obviously impossible for $N_t(v)$ to be greater than $n$, so this $t$ is too large for the approximation to hold. In any case, this shows that working in the tree-like regime is not going to suffice.

\begin{figure}[H]
\centering
\begin{subfigure}{.5\textwidth}
  \centering
  \includegraphics[width=0.9\linewidth]{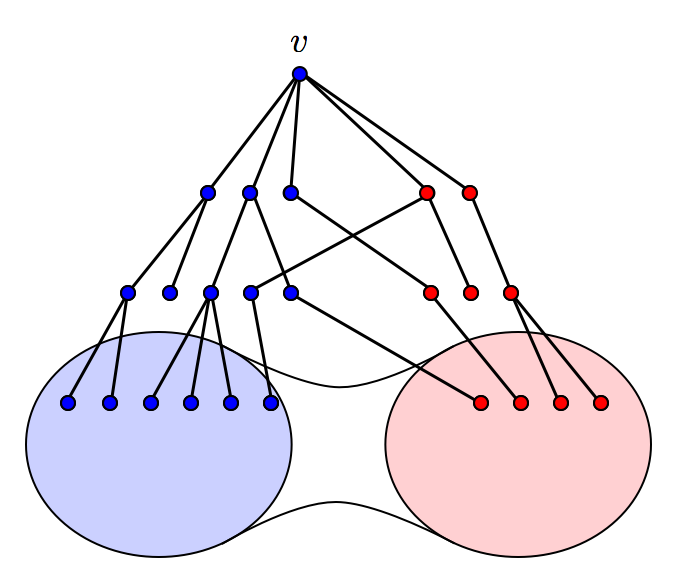}
%  \caption{Random configuration.}
  \label{depth2}
\end{subfigure}%
\begin{subfigure}{.5\textwidth}
  \centering
  \includegraphics[width=0.9\linewidth]{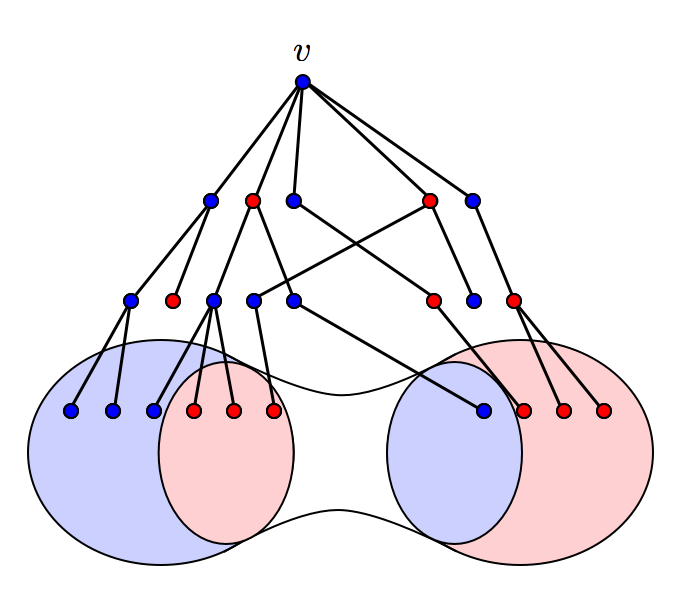}
 % \caption{Clustered configuration.}
  \label{depth2r}
\end{subfigure}
\caption{The left figure shows the neighborhood of vertex $v$ pulled from the SBM graph at depth $c\log_{\lambda_1} n$, $c<1/2$, which is a tree with high probability. If one had an educated guess about each vertex's label, of good enough accuracy, then it would be possible to amplify that guess by considering only such small neighborhoods (deciding with the majority at the leaves). However, we do not have such an educated guess. We thus initialize our labels purely at random, obtaining a small advantage of roughly $\sqrt{n}$ vertices by luck (i.e., the central limit theorem), in either an agreement or disagreement form. This is illustrated in agreement form in the right figure. 
%We next attempt to amplify that lucky guess by exploiting the information of the SBM graph. Unfortunately, the graph is too sparse to let us amplify that guess by considering tree like or even loopy neighborhoods; the vertices would have to be exhausted. This takes us to considering walks.
}       
\label{tree}
\end{figure}

This takes us to two possibilities:
\begin{itemize}
\item {\it Go deeper.} In order to amplify our weak bias of order $1/\sqrt{n}$ to a constant bias, we can go deeper in the neighborhoods, leaving the regime where the neighborhood is tree-like. In fact, according to \eqref{needed}, this requires going beyond the diameter of the graph whcih is $~\log(n)/\log((a+b)/2)$, and having to  repeat vertices (i.e., count walks). The problem is, the above approximation assumes that each vertex at a distance of $t-1$ from $v$ has one edge leading back towards $v$, and that the rest of its edges lead towards new vertices. Once a significant fraction of the vertices are fewer than $t$ edges away from $v$, a significant fraction of the edges incident to vertices $t-1$ edges away from $v$ are part of loops and thus do not lead to new vertices. This obviously creates significant complications. This is the approach that is discussed in the next section, using nonbacktracking walks. In fact, we re-derive this approach in the next section from our principles on weak recovery from Section \ref{linBP}. Note also that this approach is efficient, and it is thus legitimate to ask whether a simpler argument could be obtained information-theoretically, which takes us to the next point;
\item {\it Repeat guessing.} A single random guess gives a bias of order $1/\sqrt{n}$. However, if we keep guessing again and again, eventually one random guess will be atypically correlated with the ground truth. In particular, with enough guesses, one would get a strong enough correlation for the random guess that the {\it naive plan} described above could work at the tree-like depth (connecting us to the (robust) broadcasting on trees problem, as desired). Of course, a new difficulty is now to identify which of these many random guesses leads to a good reconstruction. For this, we propose to use a graph-splitting argument as discussed next. 
\end{itemize} 

\noindent
{\it Information-theoretic procedure:}\\ 
(1) Graph-split $G$ into $G_1$ and $G_2$ such that $G_1$ is above the KS threshold.\\
(2) Take $M=M(n)$ independent random guesses (i.e., partitions of $[n]$) and amplify each at the three-like depth on $G_1$. 
Let $\hat{X}_1,\dots, \hat{X}_M$ be the {\it amplified} guesses, which represent each a partition of $[n]$.\\
(3) Test the edge density of the residue graph $G_2$ on each partition $\hat{X}_i$, and output the first $\hat{X}_i$ that gives a non-trivial edge density (i.e., an edge density above by a constant factor of what a purely random partition gives in expectation). \\

\noindent
We conjecture that there exists an appropriate choice of $M$ such that (i) a ``good guess'' will come up whp, i.e., a guess with enough initial correlation that the naive plan described above amplifies that guess to a weak recovery solution using $G_1$, (ii) the ``good guess'' amplification is tested positive on the residue graph $G_2$ before any bad guess amplification is potentially tested positive. Note that one could use variants for testing the validity of the good guess, for example, using the $\Delta_t$ statistics of previous section to set the validity test. The advantage of this plan is that its analysis would mainly be based on estimates at tree-like depths and moment computations. 

\section{Achieving the threshold}\label{weak_achieve}
In the previous section we mentioned two plans to amplify a random guess to a valid weak recovery reconstruction: (i) one can repeat guessing exponentially many times until one hits an atypically good random guess that can be amplified on shallow neighorhoods to a valid weak recovery solution; (ii) one can take a single random guess and amplify it on deep neighborhoods to directly reach a valid weak recovery construction. We now discuss the latter plan, which can be run efficiently. 

For this, we continue the reasoning from the previous section that explained why the tree-like regime was not sufficient to amplify a random guess. An obvious way to solve the problem caused by running out of vertices would be to simply count the walks of length $t$ from $v$ to vertices in $C_1$ or $C_2$. Recall that a {\em walk} is a series of vertices such that each vertex in the walk is adjacent to the next, and a {\em path} is a walk with no repeated vertices. The last vertex of such a walk will be adjacent to an average of approximately $a/2$ vertices in its community outside the walk and $b/2$ vertices in the other community outside the walk. However, it will also be adjacent to the second to last vertex of the walk, and maybe some of the other vertices in the walk as well. As a result, the number of walks of length $t$ from $v$ to vertices in $C_1$ or $C_2$ cannot be easily predicted in terms of $v$'s community. So, the numbers of such walks are not useful for classifying vertices.

We could deal with this issue by counting paths\footnote{This type of approach is considered in \cite{bhatt-bickel}.} of length $t$ from $v$ to vertices in $C_1$ and $C_2$. The expected number of paths of length $t$ from $v$ is approximately $(\frac{a+b}{2})^t$ and the expected difference between the number that end in vertices in the same community as $v$ and the number that end in the other community is approximately $(\frac{a-b}{2})^t$. The problem with this is that counting all of these paths is inefficient.

A compromise is to count nonbacktracking walks ending at $v$, i.e. walks that never repeat the same edge twice in a row. We can efficiently determine how many nonbacktracking walks of length $t$ there are from vertices in $C_i$ to $v$. Furthermore, most nonbacktracking walks of a given length logarithmic in $n$ are paths, so it seems reasonable to expect that counting nonbacktracking walks instead of paths in our algorithm will have a negligible effect on the accuracy.

\begin{figure}[H]
\centering
  \includegraphics[width=0.52\linewidth]{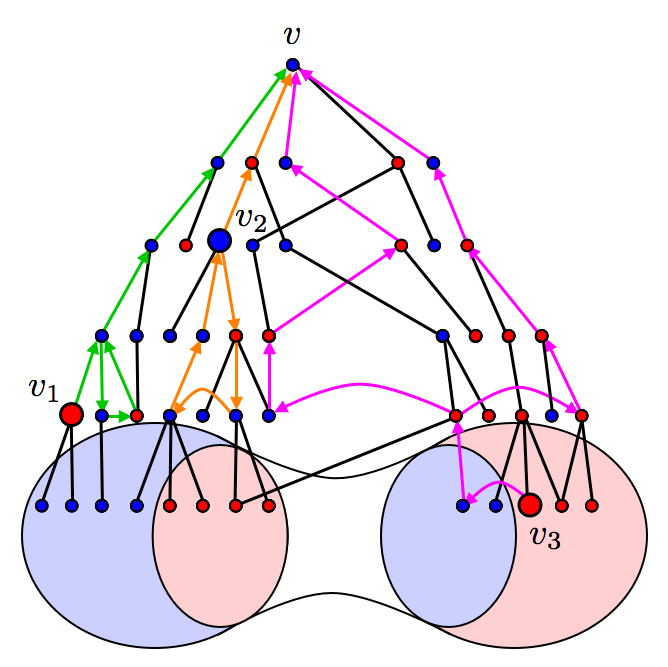}
\caption{This figure extends Figure \ref{tree} to a larger neighborhood. The ABP algorithm amplifies the belief of vertex $v$ by considering all the walks of a given length that end at it. To avoid being disrupted by backtracking or cycling the beliefs on short loops, the algorithm considers only walks that do not repeat the same vertex within $r$ steps, i.e., $r$-nonbacktracking walks. For example, when $r=3$ and when the walks have length $7$, the green walk starting at vertex $v_1$ is discarded, whereas the orange walk starting at the vertex $v_2$ is counted. Note also that the same vertex can lead to multiple walks, as illustrated with the two magenta walks from $v_3$. Since there are approximately equally many such walks between any two vertices, if the majority of the vertices were initially classified as blue, this is likely to classify all of the vertices as blue. Hence we need a compensation step to prevent the classification from becoming biased towards one community.}
\label{loops}
\end{figure}
 
\subsection{Linearized BP and the nonbacktracking matrix}\label{linBP}
To derive the algorithm more formally, we now go back to the weak recovery benchmarks discussed in Section \ref{bit-map}. The idea of using nonbacktracking walks results from a series of papers \cite{redemption,Mossel_SBM2,bordenave,colin3}, as discussed in Section \ref{history}.

Recall that the Belief Propagation algorithm presented in Section \ref{bit-map} as a derivation of the Bayes Optimal estimator. Recall also that we purposely work with a slightly more general SBM to break the symmetry: i.e., a weakly symmetric SBM with community prior $p=(p_1,p_2)$ and connectivity channel $W=Q/n$ such that $\diag(p)Q$ has constant row sums, i.e., the expected degrees in the graph are constant $d$.
As mentioned in Section \ref{bit-map}, the BP algorithm ends with a probability distribution for the community of each vertex, and taking for each vertex the most likely assignment is conjectured to give a weak recovery solution. 
However, this algorithm has several downsides. First of all, it uses a nonlinear formula to calculate each sucessive set of probability distributions (Bayes rule), and its analysis remains challenging to date. From a practical point of view, one needs to know the parameters of the model in order to run the algorithm, which makes it model-dependent.  

We now discuss how both issues can be mitigated by ``linearizing'' the algorithm. First recall that our original guesses of the vertices' communities gives only a very weak bias. 
%If vertices from some communities have higher expected degrees than others, we can quickly gain significant amounts of evidence about the vertice's communities by counting their degrees, but otherwise our beliefs about all of the vertice's communities will be fairly close to their prior probabilities for a while. 
As such, it may be useful to focus on the first order approximation of our formula when our beliefs about the communities of $v_1,...,v_m$ are all close to the prior probabilities for a vertex's community. In this case, every entry of $Qp$ must be equal to $d$. So, we have that
\begin{align*}
P[X_{v_0}=i|G=g]&\approx \frac{p_i\prod_{j=1}^m (Q q_j)_i}{\sum_{i'=1}^k p_{i'}\prod_{j=1}^m (Q q_j)_{i'}}\\ \\
&=\frac{p_i\prod_{j=1}^m (d+(Q (q_j-p))_i)}{\sum_{i'=1}^k p_{i'}\prod_{j=1}^m (d+(Q (q_j-p))_{i'})}\\ \\
&\approx \frac{p_i (1+\sum_{j=1}^m (Q (q_j-p))_i/d )}{\sum_{i'=1}^k p_{i'} (1+\sum_{j=1}^m ((Q (q_j-p))_{i'}/d))}\\ \\
&= \frac{p_i(1+\sum_{j=1}^m (Q (q_j-p))_i/d)}{1+\sum_{j=1}^m p\cdot Q (q_j-p)/d}\\ \\
&=p_i+p_i\sum_{j=1}^m (Q (q_j-p))_i/d.\\
\end{align*}

We can then rewrite the Belief Propagation Algorithm using this approximation in order to get the following algorithm.

\vspace{1 cm}
\noindent
{\em Pseudo Linearized Belief Propagation Algorithm (t, p, Q, q, G):}
\begin{enumerate}
\item Set $\epsilon^{(0)}_{v,v'}=q_{v,v'}-p$ for all $(v,v')\in E(G)$.

\item For each $0<t'< t$, and each $v\in G$, set
\[\epsilon^{(t')}_{v,v'}=\sum_{v'': (v',v'')\in E(G), v''\ne v} PQ\epsilon^{(t'-1)}_{v',v''}/d.\]

\item For each $(v,v')\in E(G)$, set 
\[q^{(t)}_{v}=p+\sum_{v': (v,v')\in E(G)} PQ\epsilon^{(t-1)}_{v,v'}/d.\]

\item Return $q^{(t)}$.
\end{enumerate}
%\vspace{1 cm}

It is time to bring up an issue that was swept under the rug until now, i.e., the effect of non-edges. 
If one has access to a good initial guess and operates a short depth, then the non-edges have negligible effects and the above algorithm can be used. However, in the current context where we will run the iteration at large depth, the non-edges need to be factored in. 
The fundamental problem is that the absence of an edge between $v$ and $v'$ provides slight evidence that these vertices are in different communities, and this algorithm fails to take that into account. Generally, as long as our current estimates of vertices communities assign the right number of vertices to each community, each vertex's nonneigbors are balanced between the communities, so the nonedges provide negligible amounts of evidence. However, if they are not taken into account, then any bias of the estimates towards one community can grow over time, and one may end up classifying all vertices in the community that was initially more represented.

%The obvious solution to this would be to add terms to the formula for $q^{(t')}_{v_0,...,v_{t-t'}}$ that account for all $v_{t-t'+1}$ that are not adjacent to $v_{t-t'}$. However, that would mean that instead of merely having a value of $q_{v_0,...,v_t}$ for every path of length $t$ or less we would need to have a value of $q_{v_0,...,v_t}$ for every $(t+1)$-tuple of vertices. That would increase the memory and time needed to run this algorithm from the already questionable $ne^{O(t)}$ to the completely impractical $O(n^{t+1})$. A reasonable way to make the algorithm more efficient is to assume that  $P[X_{v''}=i'|G\backslash\{v,v'\}=g\backslash\{v,v'\}]\approx P[X_{v''}=i'|G\backslash\{v'\}=g\backslash\{v'\}]$, which should generally hold as long as there is no small cycle containing $v$, $v'$, and $v''$. We can also reasonably assume that if $v$ is not adjacent to $v'$ then $P[X_{v'}=i'|G\backslash\{v\}=g\backslash\{v\}]\approx P[X_{v'}=i'|G=g]$. Using that to simplify the proposed approach, we would need an intial estimate of $P[X_{v'}=i|G\backslash\{v\}=g\backslash\{v\}]$, $(q_{v,v'})_i$, for each community $i$ and each adjacent $v$ and $v'$, and an initial estimate of $P[X_{v}=i|G=g]$, $(q_v)_i$ for each community $i$ and vertex $v$. 
%Then we could develop better guesses by means of the following algorithm.

Re-deriving BP by taking into account the non-edges in Bayes' rule and writing down the proper linearization leads to the following algorithm.\\

\vspace{1 cm}
\noindent
{\em Linearized Belief Propagation Algorithm (t, p, Q, q, G):}
\begin{enumerate}
\item Set $\epsilon^{(0)}_{v,v'}=q_{v,v'}-p$ for all $(v,v')\in E(G)$.

\item Set $\epsilon^{(0)}_v=q_v-p$ for all $v\in G$.

\item For each $0<t'< t$:
\begin{enumerate}

\item For each $(v,v')\in E(G)$, set
\begin{align*}
\epsilon^{(t')}_{v,v'}&=\sum_{v'': (v',v'')\in E(G), v''\ne v} PQ\epsilon^{(t'-1)}_{v',v''}/d\\
& -\sum_{v'':(v',v'')\not\in E(G),v''\ne v'} PQ\epsilon_{v''}^{(t'-1)}/n.
\end{align*}

\item For each $v\in G$, set
\begin{align*}
\epsilon^{(t')}_{v}&=\sum_{v': (v,v')\in E(G)} PQ\epsilon^{(t'-1)}_{v,v'}/d\\
&-\sum_{v': (v,v')\not\in E(G),v'\ne v} PQ\epsilon^{(t'-1)}_{v'}/n.
\end{align*}

\end{enumerate}
\item For each $v\in G$, set 
\[q^{(t)}_{v}=p+\sum_{v': (v,v')\in E(G)} PQ\epsilon^{(t-1)}_{v,v'}/d-\sum_{v': (v,v')\not\in E(G),v'\ne v} PQ\epsilon^{(t-1)}_{v'}/n.\]

\item Return $q^{(t)}$.
\end{enumerate}
\vspace{1 cm}

We will now discuss a spectral implementation of this algorithm, as the above resembles a power iteration method on a linear operator. 
Define the graph's nonbacktracking walk and adjusted nonbacktracking matrix as follows.
%\begin{definition}
%Given a graph, the graph's nonbacktracking walk matrix, $B$, is a matrix over the vector space with an orthonormal basis consisting of a vector for each directed edge in the graph. $B_{(v_1,v_2),(v'_1,v'_2)}$ is $1$ if $v'_2=v_1$ and $v'_1\ne v_2$, otherwise it is $0$. In other words, $B$ maps a directed edge to the sum of all directed edges starting at its end, except the reverse of the original edge.
%\end{definition}
\begin{definition}\label{NB}
Given a graph $(V,E)$, the graph's nonbacktracking walk matrix, $B$, is a matrix of dimension $|E_2| \times |E_2|$, where $E_2$ is the set of directed edges on $E$ (with $|E_2|=2|E|$), such that for two directed edges $e=(i,j)$, $f=(k,l)$,
\begin{align}
B_{e,f} = \1 (l=i,k\ne j).
\end{align}
In other words, $B$ maps a directed edge to the sum of all directed edges starting at its end, except for the reversal edge.
\end{definition}
\begin{definition}
Given a graph $G$, and $d>0$, the graph's adjusted nonbacktracking walk matrix, $\widehat{B}$ is the $(|E_2|+n)\times (|E_2|+n)$ matrix such that for all $w$ in the vector space with a dimension for each directed edge and each vertex, we have that $\widehat{B} w=w'$, where $w'$ is defined such that \[w'_{v,v'}=\sum_{v'': (v',v'')\in E(G), v''\ne v} w_{v',v''}/d-\sum_{v'':(v',v'')\not\in E(G),v''\ne v'} w_{v''}/n\]
for all $(v,v')\in E(G)$ and
\[w'_{v}=\sum_{v': (v,v')\in E(G)} w_{v,v'}/d-\sum_{v': (v,v')\not\in E(G),v'\ne v} w_{v'}/n\] 
for all $v\in G$.
\end{definition}

These definitions allows us to state the following fact:

\begin{theorem}
When the Linearized Belief Propagation Algorithm is run, for every $0<t'<t$, we have that
\[\epsilon^{(t')}=\left (\widehat{B}\otimes PQ\right)^{t'}\epsilon^{(0)}.\]
\end{theorem}

%\begin{proof}
This follows from the definition of $\widehat{B}$ and the fact that the propagation step of the Linearized Belief Propagation Algorithm gives $\epsilon^{(t')}=\left (\widehat{B}\otimes PQ\right)\epsilon^{(t'-1)}$ for all $0<t'<t$. 
%\end{proof}

\begin{comment}
for each $0<t'<t$, 

We procede by induction on $t'$. Obviously, $\epsilon^{(0)}=(\widehat{B}\otimes PQ)^{0}\epsilon^{(0)}$. Now, assume that this holds for $t'-1$. For every $(v,v')\in E(G)$, we have that
\begin{align*}
\epsilon^{(t')}_{v,v'}&=\sum_{v'': (v',v'')\in E(G), v''\ne v} PQ\epsilon^{(t'-1)}_{v',v''}/d-\sum_{v'': (v',v'')\not\in E(G),v''\ne v'} PQ\epsilon^{(t'-1)}_{v',v''}/n\\
&=\sum_{v'': (v',v'')\in E(G), v''\ne v} (1/d+1/n)PQ\epsilon^{(t'-1)}_{v',v''}-\sum_{v''\ne v,v''\ne v'} PQ\epsilon^{(t'-1)}_{v',v''}/n\\
&=\sum_{v'':(v',v'')\in E(G)} B_{(v',v''),(v,v')}  (PQ)\epsilon^{(t'-1)}_{v',v''}/d\\
&=[(e_{(v,v')}B)\otimes (PQ)]\sum_{(v_1,v_2)\in E(G)} e_{(v_1,v2)}\otimes \epsilon^{(t'-1)}_{v_1,v_2}/d\\
\end{align*}
So,
\begin{align*}
\epsilon^{(t')}&=\left(\frac{1}{d} B\otimes PQ\right) \epsilon^{(t'-1)}\\
&=\left(\frac{1}{d} B\otimes PQ\right)^{t'}\epsilon^{(0)}
\end{align*}
\end{comment}

In other words, the Linearized Belief Propagation Algorithm is essentially a power iteration algorithm that finds the eigenvector of $\widehat{B}\otimes PQ$ with the largest eigenvalue. $\widehat{B}\otimes PQ$ has an eigenbasis consisting of tensor products of eigenvectors of $\widehat{B}$ and eigenvectors of $PQ$, with eigenvalues equal to the products of the corresponding eigenvalues of $\widehat{B}$ and $PQ$. As such, this suggests that for large $t'$, $\epsilon^{(t')}$ would be approximately equal to a tensor product of eigenvectors of $\widehat{B}$ and $PQ$ with maximum corresponding eigenvalues. For the sake of concreteness, assume that $w$ and $\rho$ are eigenvectors of $\widehat{B}$ and $PQ$ such that $\epsilon^{(t')}\approx w\otimes \rho$. That corresponds to estimating that
\[P[X_v=i|G=g]\approx p_i+w_v \rho_i\]
for each vertex $v$ and community $i$. If these estimates are any good, they must estimate that there are approximately $p_in$ vertices in community $i$ for each $i$. In other words, it must be the case that the sum over all vertices of the estimated probabilities that they are in community $i$ is approximately $p_i n$. That means that either $\rho$ is small, in which case these estimates are trivial, or $\sum_{v\in G} w_v\approx 0$. Now, let the eigenvalue corresponding to $w$ be $\lambda$. If $\sum_{v\in G} w_v\approx 0$, then for each $(v,v')\in E(G)$, we have that
\begin{align*}
\lambda w_{v,v'} &\approx \sum_{v'':(v',v'')\in E(G),v''\ne v} w_{v',v''}/d\\
&= \sum_{v'':(v',v'')\in E(G)} B_{(v',v''),(v,v')} w_{v',v''}/d
\end{align*}
So, the restriction of $w$ to the space spanned by vectors corresponding to directed edges will be approximately an eigenvector of $B$ with an eigenvalue of approximately $\lambda/d$. Conversely, any eigenvector of $B$ that is balanced in the sense that its entries add up to approximately $0$ should correspond to an eigenvector of $\widehat{B}$. So, we could try to determine what communities vertices were in by finding some of the balanced eigenvectors of $B$ with the largest eignvalues, adding together the entries corresponding to edges ending at each vertex, and thresholding.The eigenvector of $B$ with the largest eigenvalue will have solely nonnegative entries, so it will not be balanced. However, it is reasonable to expect that its next few eigenvectors would be relatively well balanced. 

This approach has a couple of advantages over the full BP algorithm. First of all, one does not need to know anything about the graph's parameters to find the top few eigenvectors of $B$, so this algorithm works on a graph drawn from an SBM with unknown parameters. Secondly, the approximation of the top few eigenvectors of $B$ will tend to be simpler than the analysis of the BP algorithm. Note that balanced eigenvectors of $B$ will be approximately eigenvectors of $B-\frac{\lambda_1}{n}\mathbb{1}$, where $\lambda_1$ is the largest eigenvalue of $B$ and $\mathbb{1}$ is the matrix thats entries are all $1$. Therefore, we could also look for the main eigenvectors of $B-\frac{\lambda_1}{n}\mathbb{1}$ instead of looking for the main balanced eigenvectors of $B$. We give next two variants of the resulting algorithm for the case of the SSBM.\footnote{Note that the symmetry breaking is used to derived the algorithm but we can now apply it equally well to the symmetric SBM.}\\

\noindent
{\bf Nonbacktracking eigenvector extraction algorithm \cite{redemption,bordenave}.}\\
Input: A graph $G$ and a parameter $\tau \in \mR$. \\
(1) Construct the nonbacktracking matrix $B$ of the graph $G$.\\
(2) Extract the eigenvector $\xi_2$ corresponding to the second largest eigenvalue of $B$.\\
(3) Assign vertex $v$ to the first community if $\sum_{e : e_2=v} \xi_2(e)> \tau/\sqrt{|V(G)|}$ and to the second community otherwise. \\

\begin{theorem}\cite{bordenave}
If $(a-b)^2>2(a+b)$, then there exists $\tau \in \mR$ such that previous algorithm solves weak recovery in $\ssbm(n,2,a/n,b/n)$.
\end{theorem}
%We believe that the above algorithm should also work by simply taking the sign of $\sum_{e : e_2=v} \xi_2(e)$ to obtain the clusters.

Extracting the second eigenvector of the nonbacktracking matrix directly may not be the most efficient way to proceed, especially as the graph gets denser. A power iteration method is a natural implementation, which requires additional proofs as done in \cite{colin3cpam}. 
%The approach of \cite{Mossel_SBM2} based on the count of weighted nonbacktracking walks between vertices provides another practical alternative. 
Below is the message-passing implementation.\\

\noindent
{\bf Approximate Belief Propagation (ABP) algorithm. \cite{colin3nips,colin3cpam}}\\
Inputs: A graph $G$ and a parameter $m \in \mZ_+$.\\
(1) For each adjacent $v$ and $v'$ in $G$, randomly draw $y^{(1)}_{v,v'}$ from a Gaussian distribution with mean $0$ and variance $1$. Assign $y^{(t)}_{v,v'}$ to value of $0$ for $t<1$.\\
(2) For each $1<t\le m$, set for all adjacent $v$ and $v'$  
\begin{align*}
z^{(t-1)}_{v,v'}&=y^{(t-1)}_{v,v'}-\frac{1}{2|E(G)|}\sum_{(v'',v''')\in E(G)} y^{(t-1)}_{v'',v'''}, \\
y^{(t)}_{v,v'}&=\sum_{v'':(v',v'')\in E(G),v''\ne v} z^{(t-1)}_{v',v''}.
\end{align*}
%For each adjacent $v,v'$ in $G$ that are not part of a cycle of length $r$ or less, set
%\[y^{(t)}_{v,v'}=\sum_{v'':(v',v'')\in E(G),v''\ne v} z^{(t-1)}_{v',v''},\]
%and for the other adjacent $v,v'$ in $G$, let the other vertex in the cycle that is adjacent to $v$ be $v'''$, the length of the cycle be $r'$, and set \[y_{v,v'}^{(t)}=\sum_{v'':(v',v'')\in E(G),v''\ne v} z_{v',v''}^{(t-1)}-\sum_{v'':(v,v'')\in E(G),v''\ne v',v''\ne v'''} z_{v,v''}^{(t-r')}\]
%unless $t=r'$, in which case, set $y_{v,v'}^{(t)}=\sum_{v'':(v',v'')\in E(G),v''\ne v} z_{v',v''}^{(t-1)}-z^{(1)}_{v''',v}$.
(3) Set for all $v\in G$, $y'_v=\sum_{v':(v',v)\in E(G)} y^{(m)}_{v,v'}$.
 Return $(\{v:y'_v>0\},\{v: y'_v\le 0\})$.\\

In \cite{colin3cpam}, an extension of the above algorithm that prohibits backtrack of higher order (i.e, avoiding short loops rather than just self-loops) is shown to achieve the threshold for weak recovery in the SBM when $m=2\log(n)/\log(\snr)+\omega(1)$. The idea of prohibiting short loops is further discussed in the next section.

\subsection{Algorithms robustness and graph powering}\label{robust}
\begin{figure}[h]
\begin{center}
\includegraphics[scale=.32]{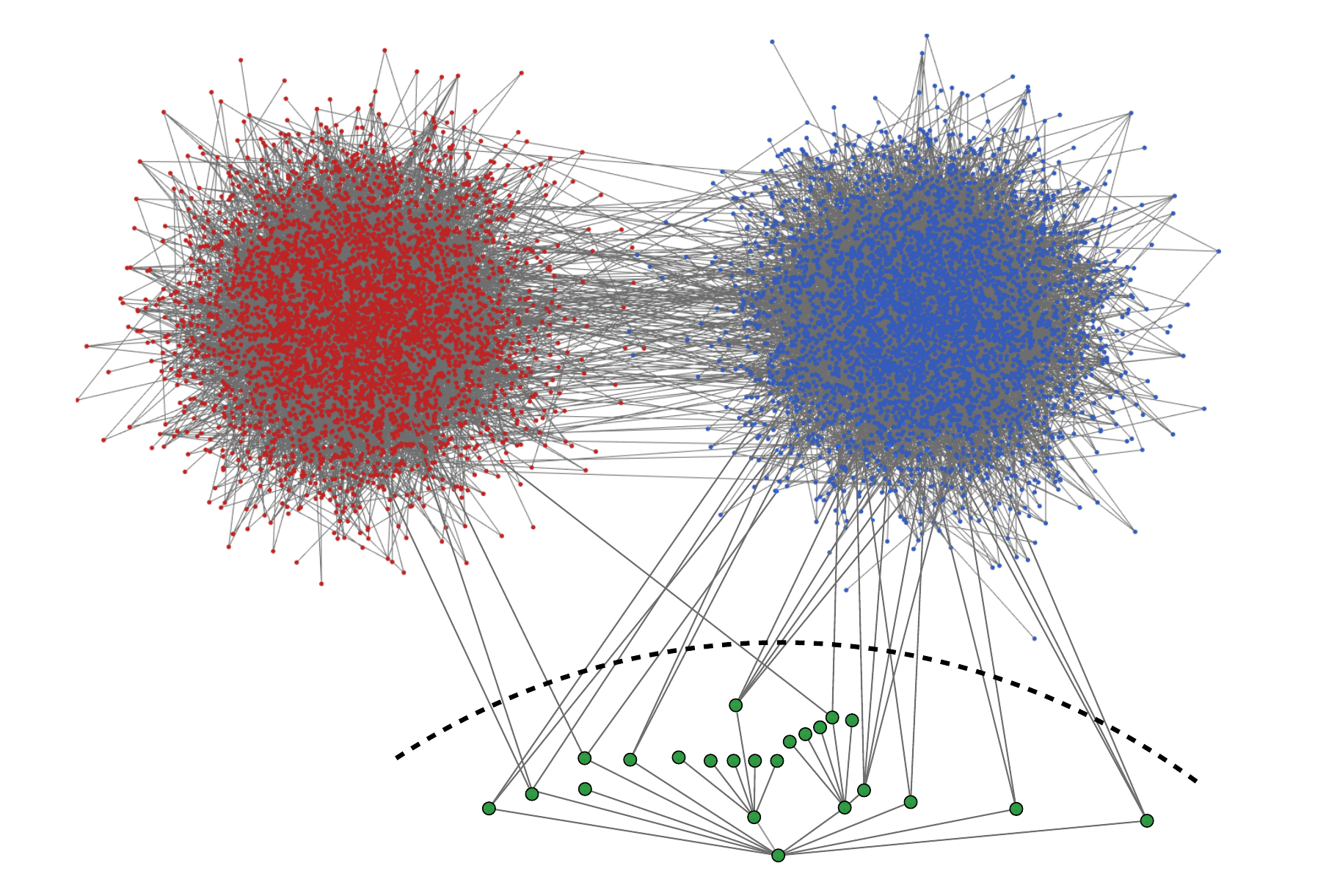}
\caption{The communities obtained with the spectral algorithm on the adjacency matrix in a sparse symmetric SBM above the KS threshold ($n=100000, a=2.2, b=0.06$): one community corresponds to the neighborhood of a high-degree vertex, and all other vertices are put in the second community.}
\label{high-degree}
\end{center}
\end{figure}

The quick intuition on why the nonbacktracking matrix is more amenable to community detection than the adjacency matrix can be seen by taking powers of these matrices. In the case of the adjacency matrix, powers are counting walks from a vertex to another, and these get amplified around high-degree vertices since the walk can come in and out in many ways. This creates large eigenvalues with eigenvectors localized around high-degree vertices. This phenomenon is well documented in the literature; see Figure \ref{high-degree} for a illustration of a real output of the spectral algorithm on a SBM with two symmetric communities (above the KS threshold).\\

Instead, by construction of the nonbacktracking matrix, taking powers forces a directed edge to leave to another directed edge that does not backtrack, preventing such amplifications around high-degree vertices. So nonbacktracking gives a way to mitigate the degree-variations and to avoid localized eigenvectors (recall discussion in Section \ref{tackle}). Note also that one cannot simply remove the high-degree vertices in order to achieve the threshold; one would have to remove too many of them and the graph would lose the information about the communities. This is one of the reasons why the weak recovery regime is interesting. \\

\begin{figure}[h]
\begin{center}
\includegraphics[scale=.4]{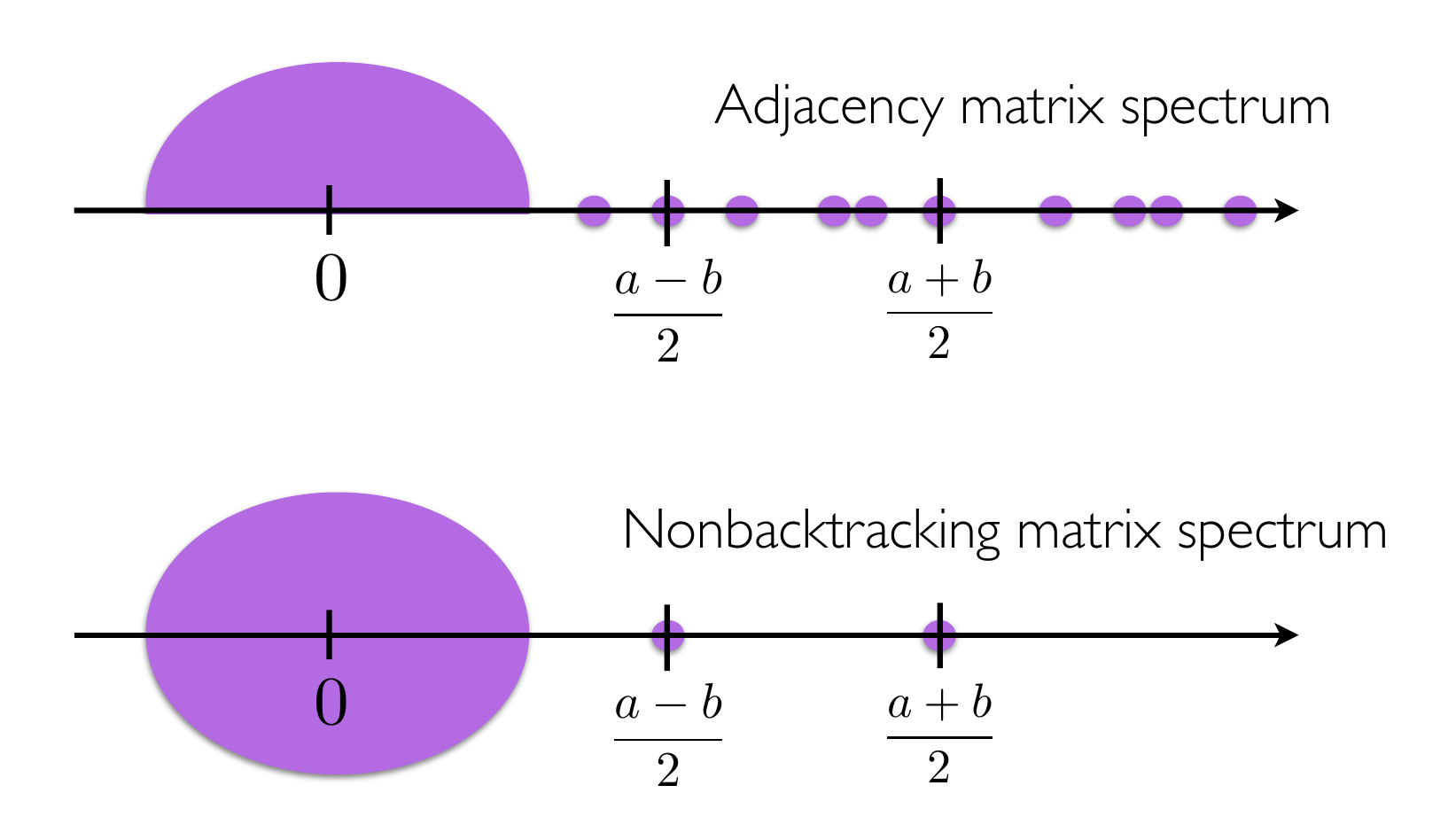}
\caption{Illustration of the spectrum of the adjacency and nonbacktracking matrices for the SBM with two symmetric communities above the KS threshold.}
\label{2spectrums}
\end{center}
\end{figure}

This robustness property of the nonbacktracking matrix is reflected in its spectrum, which has largest magnitude eigenvalue $\lambda_1$ (which is real positive), and second largest magnitude eigenvalue $\lambda_2$ which appears before $\sqrt{\lambda_1}$ above the KS threshold:
\begin{align}
\sqrt{\lambda_1} < |\lambda_2| < \lambda_1.
\end{align}
Then weak recovery can be solved by using the eigenvector corresponding to $\lambda_2$; see the previous section. 
Figure \ref{2spectrums} provides an illustration for the SBM with two symmetric communities.

However the robustness of the NB matrix may not be as strong as desired. It happens that in the SBM, being pushed away from a high-degree vertex makes it unlikely for the walk to go back to a high-degree vertex. Therefore, avoiding direct backtracks suffices. Unfortunately, in many real data graphs, loops are much more frequent than they are in the SBM. Consider for example the geometric block model with two Gaussians discussed in Section \ref{others}; in such a model, being pushed away from a high degree vertex likely brings the walk back to another neighbor of that same high degree vertex, and prohibiting direct backtracks does not help much. In fact, this issue is also present for BP itself (rather than linearized BP), which is originally designed\footnote{Although it also works in some loopy context \cite{loopy}; in addition to the AMP framework that applies to the cases of denser graphs} for locally tree-like models as motived in Section \ref{bit-map}, although BP has the advantage over ABP to pass probability messages that cannot grow out of proportions (being bounded to $[0,1]$). 

\begin{figure}[h]
\begin{center}
\includegraphics[scale=.33]{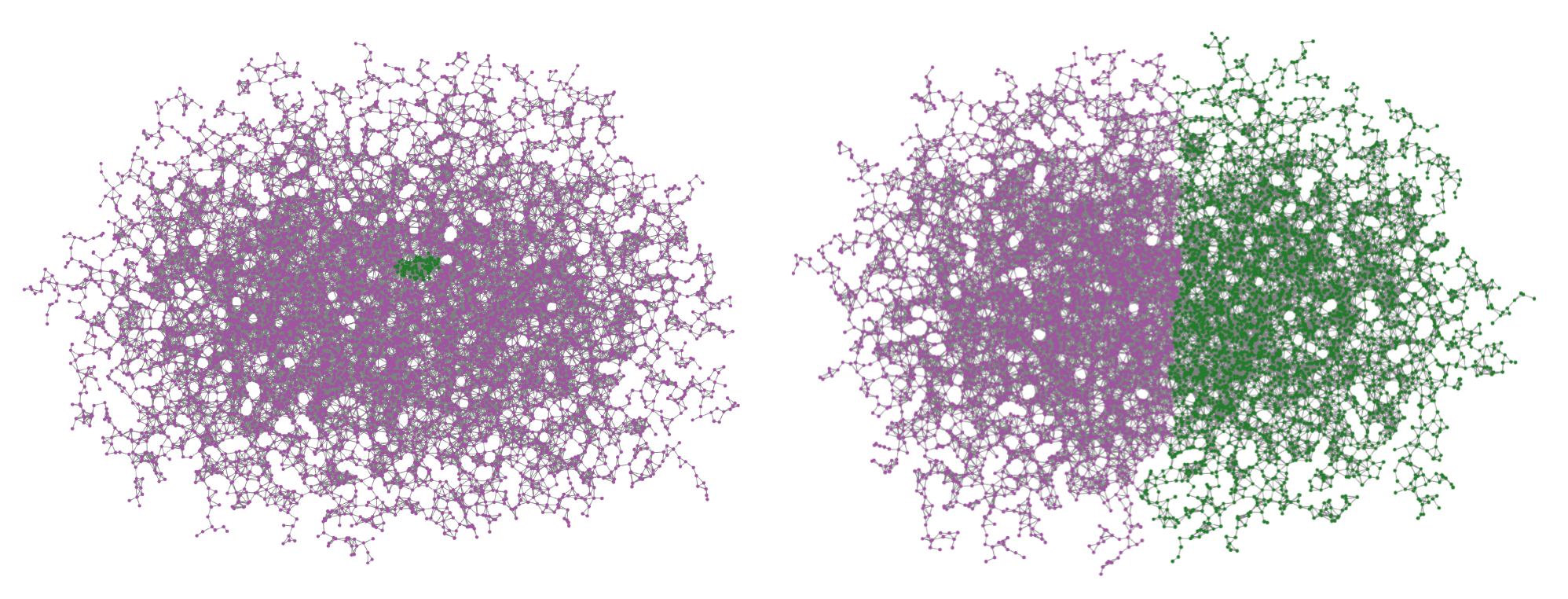}
\caption{A graph drawn from the mixture-GBM$(n,2,T,S)$ defined in Section \ref{others}, where $n/2$ points are sampled i.i.d.\ from an isotropic Gaussian in dimension 2 centered at $(0,0)$ and $n/2$ points are sampled i.i.d.\ from an isotropic Gaussian in dimension 2 centered at $(S,0)$, and any points at distance less than $T$ are connected (here $n=10000$, $S=2$ and $T=10/\sqrt{n}$). The spectral algorithm on the NB matrix gives the right plot, which puts a small fraction of densely connected vertices (a tangle) in a community, and all other vertices in the second community. The right plot is the desired community output, which graph powering produces.}
\label{gbm-comp}
\end{center}
\end{figure}

A natural attempt to improve on this is to then extend the notion of nonbacktracking beyond direct backtracks, prohibiting any repeat of a vertex within $r$ steps of the walk (rather than just 2 steps). In fact, this idea was already used in \cite{colin3cpam} for the SBM, as the increased robustness also helped with simplifying the proofs (even though it is likely unnecessary for the final result to hold). We now formally define the $r$-NB matrix of a graph:

\begin{definition}\label{r-NB}[The $r$-nonbacktracking ($r$-NB) matrix.] Let $G=(V,E)$ be a simple graph and let $\vec{E}_r$ be the set of directed paths of length $r-1$ obtained on $E$. The $r$-nonbacktracking matrix $B^{(r)}$ is a $|\vec{E}_r| \times |\vec{E}_r|$ matrix indexed by the elements of $\vec{E}_r$ such that, for two sequences of $r-1$ direct edges $e=(e_1,\dots,e_{r-1}),f=(f_1,\dots,f_{r-1})$ that form directed paths in $\vec{E}_r$,
\begin{align}
B^{(r)}_{e,f}=\prod_{i=1}^{r-1} \1(f_{i+1}=e_i)\1((f_{1})_1\neq (e_{r-1})_2),
\end{align}
i.e., entry $(e,f)$ of $B^{(r)}$ is 1 if $e$ extends $f$ by one edge without creating a loop, and 0 otherwise. 
\end{definition}

\begin{figure}[H]
\centering
  \includegraphics[width=0.4\linewidth]{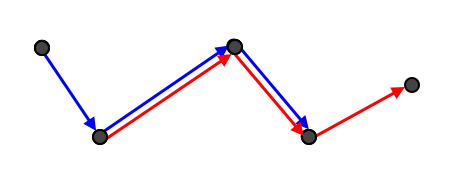}
\caption{Two paths of length $3$ that contribute to an entry of $1$ in $B^{(4)}$.}
\label{b4p}
\end{figure}

\begin{remark}
Note that $B^{(2)}=B$ is the classical nonbacktracking matrix from Definition \ref{NB}. As for $r=2$, we have that $((B^{(r)})^{k-1})_{e,f}$ counts the number of $r$-nonbacktracking walks of length $k$ from $f$ to $e$.  
\end{remark}

While one gains further robustness by using $r$-NB with larger $r$, this may still require $r$ to be impractically large in cases where a large number of cliques and tangles are present in the graph (such as in the two-Gaussian geometric block model mentioned before). We now discuss three alternatives to further increase such robustness.\\

\noindent
{\bf (1) SDPs.} While the SDP in Section \ref{achieve-exact} was motivated as a relaxation of the min-bisection estimator that is optimal for exact recovery but not necessarily for weak recovery, one may still use SDPs for weak recovery as well. In fact, the SDP benefits from a key feature; recall that the SDPs discussed in Section \ref{achieve-exact} takes the form:
\begin{align}
\hat{X}_{SDP}(g) = \argmax_{  X \succeq 0 \atop{ X_{ii} = 1, \forall i \in [n]  }} \tr(B X).
\end{align}
for a matrix $B$ which is a centered version of the adjacency matrix. The key feature is that the constraint 
$$ X_{ii} = 1$$
on the matrix $X$ does not make the above optimization hard, as opposed to the original min-bisection problem that requires
$$ x_{i}^2 = 1$$
on the vector $x$, which makes min-bisection an NP-hard integral optimization problem. 
The advantage is that $X_{ii} = 1$ prohibits the entries of $X$ to grow out of proportion ($X$ needs to be also PSD), and the SDP is less sensitive to producing localized eigenvector. 

Several works have investigated SDPs for SBMs \cite{levina,abh,afonso_single,new-xu,wein_sdp}, with a precise picture obtained for weak recovery in \cite{sbm-groth,montanari_sen,adel_sbm}. We now mention a result from \cite{montanari_sen} that shows that SDPs allow to approach the threshold for weak recovery in the two-community SSBM arbitrarily close when the expected degrees diverge. 

\begin{theorem}\cite{montanari_sen}
There exists $\delta(a,b)$ satisfying $\delta(a,b)\to 0$ as $(a+b) \to \infty$, such that if  
$\frac{(a-b)^2}{2(a+b)} >1 + \delta(a,b)$, then the SDP solves weak recovery in $\ssbm(n,2,a/n,b/n)$.
\end{theorem}

So for large degrees, SDPs are both performing well and allowing for further robustness compared to NB spectral methods. For example, \cite{feige_monotone,ankur_SBM} shows that SDPs are robust to certain monotone adversary, that can add edges within clusters and remove edges across clusters. Such adversary could instead create trouble to NB spectral, e.g., by adding a clique within a community to create a localized eigenvector of large eigenvalue. 

\looseness=-1On the flip side, SDPs have two issues: (1) they do not perform as well in very sparse regimes; take for instance the example covered in Section \ref{achieve-exact} showing that the SDP fails to find the clusters when they are disjoint and which can generalize to more subtle cases where the sparsest cut is at the periphery of the clusters rather than in their middle; (2) most importantly, SDPs are not practical on large graphs. One may use various tricks to initialize or accelerate SDPs, but these instead make the analysis more challenging.%\\  
\medskip

\noindent
{\bf (2) Laplacian and Normalized Laplacian.}
In contrast to SDPs, spectral methods afford much better complexity attributes. 
A possible way to improve their robustness to degree variations would be to simply normalize the matrix by taking into account degree variations, such as 
\begin{align}
L&:=D-A  \label{L1},\\
L_{\mathrm{norm}}&:= I -D^{-1/2} A D^{-1/2}  \quad \Leftrightarrow \quad D^{-1/2} A D^{-1/2} .  \label{L2}
\end{align}
These can also be viewed as relaxations of min-cuts where one does not constrain the number of vertices to be strictly balanced in each community, but where one weighs in the volume of the communities in terms of the number of vertices or degrees:
\begin{align}
\mathrm{normcut}_1&:= \frac{|\partial(S)|}{|S||S^c|} |V| = \frac{|\partial(S)|}{|S|} + \frac{|\partial(S)|}{|S^c|}\\
\mathrm{normcut}_2&:= \frac{|\partial(S)|}{d(S)d(S^c)} d(V) = \frac{|\partial(S)|}{d(S)} + \frac{|\partial(S)|}{d(S^c)}
\end{align}
where $d(S)=\sum_{v \in S}\mathrm{degree}(v)$, $V=[n]$.

One can now look for a $\{0,1\}$-valued vector, i.e., the indicator vector on a subset $S$, that minimizes these normalized cuts. This is still an NP-hard problem, but its spectral relaxation obtained by removing the integral constraints leads to the smallest eigenvector of the matrices in \eqref{L1} and \eqref{L2} (ignoring the 0 eigenvalue), which now have some balanceness properties embedded. 

\begin{figure}[h]
\begin{center}
\includegraphics[scale=.32]{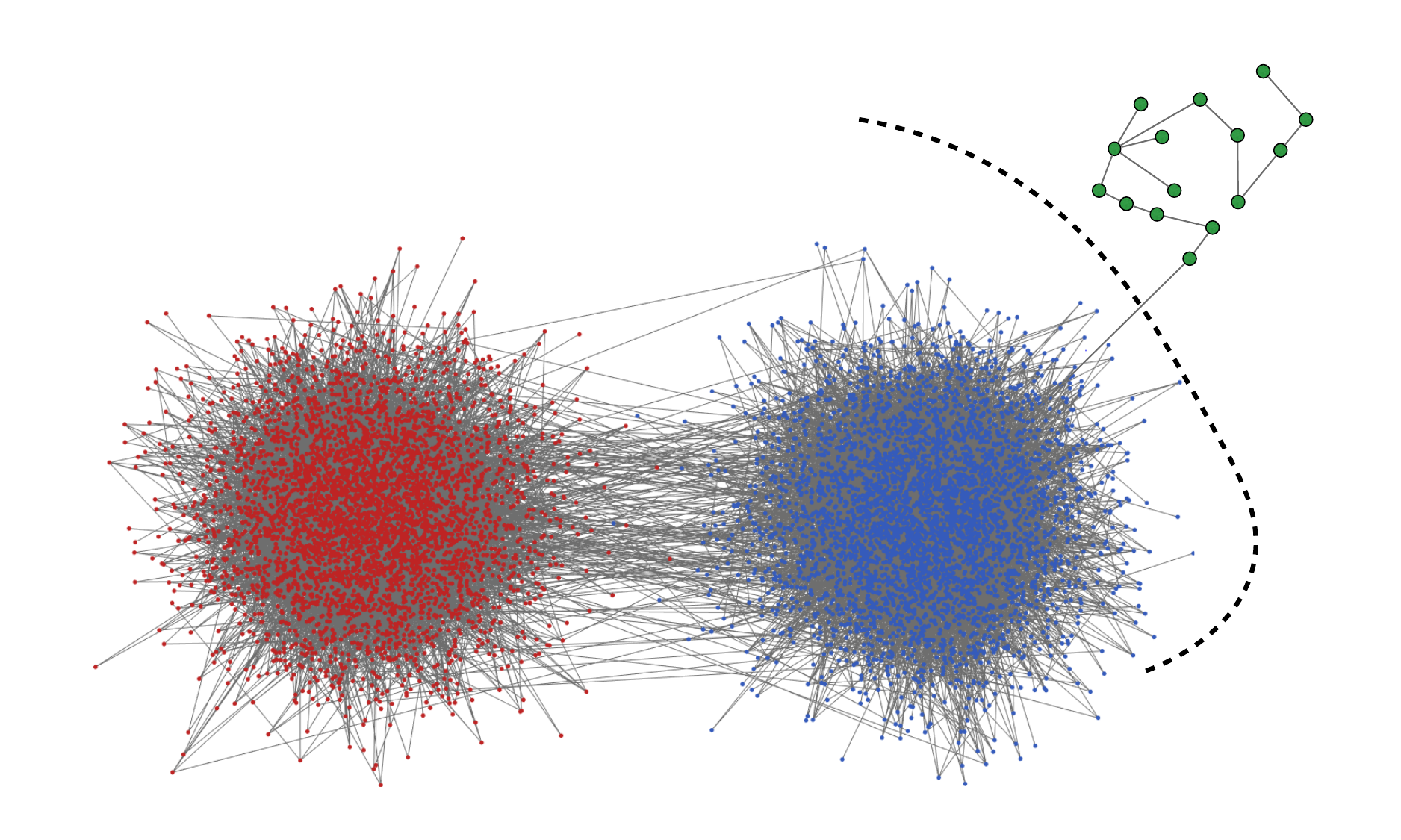}
\caption{The communities obtained with the spectral algorithm on the Laplacian matrix in a sparse symmetric SBM above the KS threshold ($n=100000, a=2.2, b=0.06$): one community corresponds to a ``tail'' of the graph (i.e., a small region connected by a single edge to the giant component), and all other vertices are put in the second community. The same outcome takes place for the normalized Laplacian.}
\label{lap-cut}
\end{center}
\end{figure}

In fact these afford better robustness to high-degree vertices. However, they tend to overdo the degree correction and can have trouble with low-degree regions of the graph in models such as the SBM. Attached are two examples of SBMs that are above the weak recovery threshold, but where these two normalized spectral methods produce communities that are peripheral, i.e., cutting a small tail of the graph that has a sparse normalized cut --- see Figure \ref{lap-cut}.

\looseness=-1In fact, we conjecture that neither of these two operators achieve the weak recovery threshold in general. But, more importantly, one should be reminded of the principled approach pursued here: The normalized Laplacians are motivated by combinatorial benchmarks, i.e., normalized cuts, which do not have a clear connection to the Bayes optimal estimator. \medskip

\noindent
{\bf (3) Graph powering.} We conclude with a recent proposal to bridge the advantage of spectral methods with robustness attributes, while keeping a Bayes-inspired construction. The method is developed in \cite{powering1} and relies on the following operator. 

\begin{definition}[Graph powering]
We give two equivalent definitions: 
\begin{itemize}
\item Given a graph $G$ and a positive integer $r$, define the $r$-th graph power of $G$ as the graph $G^{(r)}$ with the same vertex set as $G$ and where two vertices are connected if there exists a path of length $\le r$ between them in $G$. 

Note: $G^{(r)}$ contains the edges of $G$ and adds edges between any two vertices at distance $\le r$ in $G$. Note also that one can equivalently ask for a walk of length $\le r$, rather than a path. 

\item If $A$ denotes the adjacency matrix of $G$ with 1 on the diagonal (i.e., put self-loops to every vertex of $G$), then the adjacency matrix $A^{(r)}$ of $G^{(r)}$ is defined by 
\begin{align}
A^{(r)} = \mathbb{1} (A^r \ge 1).
\end{align}
Note: $A^r$ has the same spectrum as $A$ (up to rescaling), but the action of the non-linearity $\mathbb{1} (\cdot \ge 1)$ gives a key modification to the spectrum. 
\end{itemize}
\end{definition}

\begin{definition}[Deep cuts] For a graph $G$, a $r$-deepcut in $G$ corresponds to a cut in $G^{(r)}$, i.e., 
\begin{align}
\partial_r(S) = \{ u \in S, v \in S^c: (A^{(r)})_{u,v}=1\}, \quad S \subseteq V(G).
\end{align}
\end{definition}
We now discuss two key attributes of graph powering: 
\begin{itemize}
\item {\it Deep cuts as Bayes-like cuts.}
The cut-based algorithms discussed previously for $A$, $L$ or $L_{\mathrm{norm}}$ can be viewed as relaxations of the MAP estimator, i.e., min-bisection. As said multiple times, this is not the right objective for weak recovery. Let us again illustrate the distinction on a toy example, illustrated in Figure \ref{bad-tree}.

\begin{figure}[h]
\begin{center}
\includegraphics[scale=.35]{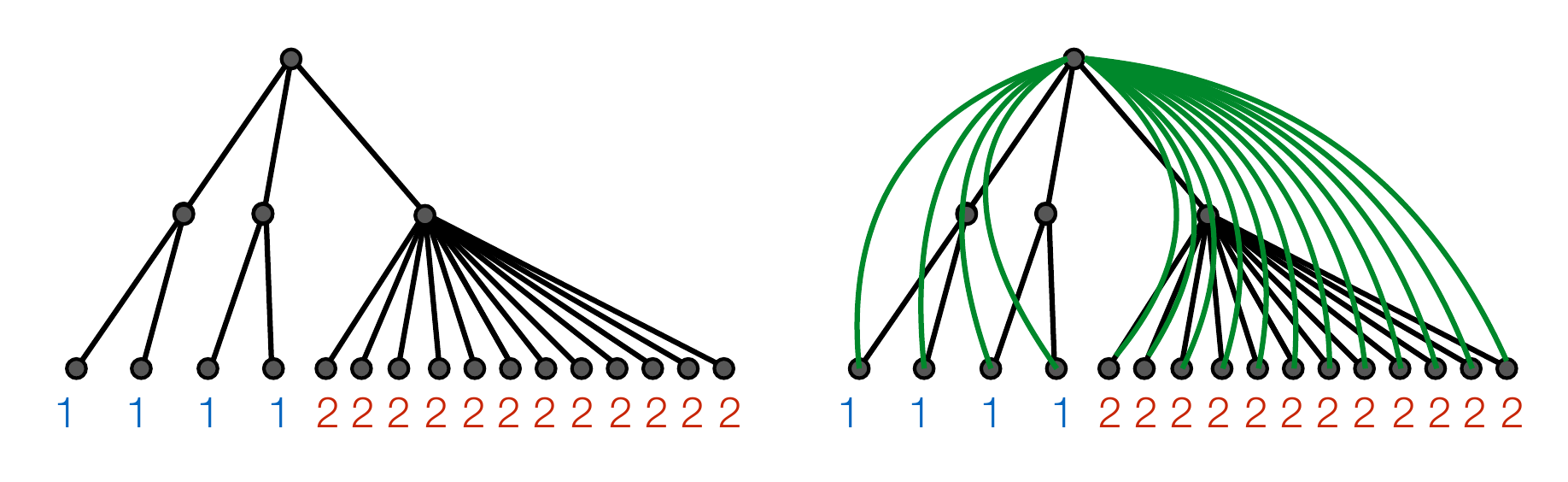}
\caption{In the left graph, assumed to come from $\ssbm(n,2,3/n,2/n)$, the root vertex is labelled community 1 from the ML estimator given the leaf labels, which corresponds to the min-cut around that vertex. 
In contrast, the Bayes optimal estimator puts the root vertex in community 2, as the belief of its right descendent towards community 2 is much stronger than the belief of its two left descendents towards community 1. This corresponds in fact to the min-deep-cut obtained from the right graph, where 2-deep edges are added by a graph-power.}
\label{bad-tree}
\end{center}
\end{figure}

Imagine that a graph drawn from $\ssbm(n,2,3/n,2/n)$ contained the following induced subgraph. $v_0$ is adjacent to $v_1$, $v_2$, and $v_3$. $v_1$ and $v_2$ are each adjacent to two outside vertices that are known to be in community $1$, and $v_3$ is adjacent to a large number of vertices that are known to be in community $2$. $v_1$ and $v_2$ are more likely to be in community $1$ than they are to be in community $2$, and $v_3$ is more likely to be in community $2$ than it is to be in community $1$. So, the single most likely scenario is that $v_0$, $v_1$, and $v_2$ are in community $1$ while $v_3$ is in community $2$. In particular, this puts $v_0$ in the community that produces the sparsest cut (1 edge in the cut vs 2 edges in the other case). However, $v_3$ is almost certain to be in community $2$, while if we disregard any evidence provided by their adjacency to $v_0$, we would conclude that $v_1$ and $v_2$ are each only about $69\%$ likely to be in community $1$. As a result, $v_0$ is actually slightly more likely to be in community $2$ than it is to be in community $1$.

The reason why powering and deepcuts matter is that it helps getting feedback from vertices that are further away, producing a form of combinatorial likelihood that is measured by the number of vertices that are not just directly connected to a vertex, but also neighbors at a deeper depth. The deepcuts are thus more ``Bayes-like'' and less ``MAP-like,'' as seen in the previous example where vertex $1$ is now assigned to community $2$ using $2$-deepcuts rather than community $1$ using standard cuts (i.e., $1$-deepcuts). 

\item {\it Powering to homogenize the graph.} Powering further helps mitigate the degree variations, and more generally density variations in the graph, both with respect to  high and low densities. Since the degree of all vertices is raised with powering, both high and low density regions do not contrast as much under powering. Large degree vertices (as in Figure \ref{high-degree}) do not stick out as much and tails (as in Figure \ref{lap-cut}) are thickened, and the more macroscopic properties of the graph can prevail.
\end{itemize}
Of course, powering is useful only if the power $r$ is not too low or not too large. If it is too low, say $r=2$, powering may not help. If it is too large, say $r\ge \mathrm{diameter}(G)$, then powering turns any graph to a complete graph, which destroys all the useful information. However, powering with $r$ below the diameter and larger than $\log\log(n)$, such as $r=\lfloor \sqrt{\log(n)} \rfloor$, allows to regularize the SBM graph to achieve the weak recovery threshold with the vanilla spectral algorithm. 

The key property is captured by the following pictorial representation of the spectrum of $A^{(r)}$ (say for $r=\lfloor \sqrt{\log(n)} \rfloor$):

\begin{figure}[h]
\begin{center}
\includegraphics[scale=.3]{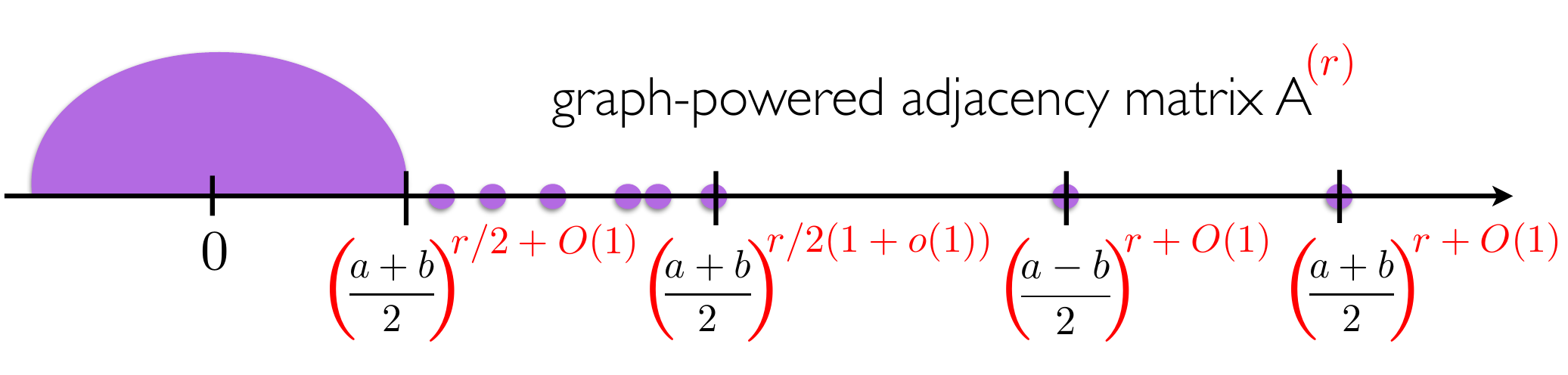}
\caption{Illustration of the spectrum of the adjacency matrix of the powered graph for a two community symmetric SBM above the KS threshold. The power must be below half the diameter and larger than $\Omega(\log\log(n))$, such as $r=(\log\log(n))^2$ or $r=\e \log n$, $\e$ small enough. The bulk is delimited by the eigenvalues of random vectors, while the localized eigenvectors on high-degree vertices mark the next transition (note that the previous two regions may not appear separated as in the figure on a real plot), followed by the isolated eigenvalue containing the community information, and at last the Perron-Frobenius eigenvalue close to the average degree. 
A similar picture takes place for the operator of \cite{massoulie-STOC}, which does not use the non-linearity of graph-powering, and which is more sensitive to tangles as for the NB operator in Figure \ref{gbm-comp}.
}
\label{powering-spectrum}
\end{center}
\end{figure}

Note also that a similar picture to Figure \ref{powering-spectrum} holds in the SBM when taking a different operator, namely $W$ where $W_{ij}$ counts the number of paths of length $r$ between $i$ and $j$, as shown by Massouli\'e in the first proof of the KS threshold achievability \cite{massoulie-STOC}. The key difference is that the operator $W$ suffers from the same issues as discussed above for the nonbacktracking operator: it allows us to mitigate high-degree vertices, but not denser regions such as tangles. In particular, planting a moderately small clique in a community of the SBM or taking the mixture GBM as in Figure  \ref{gbm-comp} can drive the spectral algorithm on $W$ to produce localized eigenvector, while the powering operator is more robust due to the non-linearity $\mathbb{1} (\cdot \ge 1)$ that flattens out the entries of large magnitude. 

Note also that one could look for a non-linearity function that is `optimal' (say for the agreement) rather than $\mathbb{1} (\cdot \ge 1)$; however this is likely to be model dependent, while the previous choice seems to be natural and generic. A downside of powering approaches is that they densify the graph, so one would ideally combine graph powering with degree normalizations to reduce the number of powers (since powering raises the density of the graph, normalization may no longer have the issues mentioned previously with tails) or some form of graph sparsification (such as \cite{sparsify}). Note that powering and sparsifying do not cancel each other: powering adds edges to ``complete the graph'' where edges should be present, while sparsifying prunes down the graph by adding weights on representative edges. Finally, one may also peel out leaves (with a few recursions) and paths to further reduce the powering order.

\chapter{Partial recovery for two communities}\label{partial}
In this section, we discuss various results for partial recovery.
\section{Almost Exact Recovery}\label{almost}
Almost exact recovery, also called weak consistency in the statistics literature, or strong recovery, has been investigated in various papers such as \cite{prout_colt,levina,harrison,mossel-consist,prout,colin1focs}. 

\begin{theorem}\label{almost_exact}
Almost exact recovery is solvable in $\ssbm(n,2,a_n/n,$ $b_n/n)$ (and efficiently so) if and only if
\begin{align}
\frac{(a_n-b_n)^2}{2(a_n+b_n)}= \omega(1).
\end{align} 
\end{theorem}
Note that the constant $2$ above is not necessary but it makes the connection to the $\snr$ of previous section more explicit.  
This result appeared in several papers; a first appearance is in \cite{prout_colt} where it results from the case of non-adaptive random samples, and it is also shown in \cite{mossel-consist,colin1focs}. 

We point out two approaches to obtain this result:
\begin{itemize}
\item {\it Boosting weak recovery with graph-splitting.} One can use graph-splitting repeatedly to turn the results of the previous section on weak recovery into Theorem \ref{almost_exact}. However this requires having an algorithm to solve weak recovery first, which represents more work than needed to obtain Theorem \ref{almost_exact}. We mention nonetheless how one can take a shortcut assuming such an algorithm. 

The idea is to graph-split $G$ into $k$ subgraphs with equal split-probabilities into $G_1,\dots,G_k$ with $k=\lfloor \sqrt{\log(n)} \rfloor$. Note that each $G_i$ is marginally an SBM with parameters $a_n=a \sqrt{\log(n)} $, $b_n=b \sqrt{\log(n)}$, so largely above the KS threshold. Now apply the algorithm that solves weak recovery in each of these graphs, to obtain a collection of $k$ clusterings $\hat{X}_1, \dots, \hat{X}_k$. One can now boost the accuracy by doing a vote for each pair of vertices over these different clusterings. E.g., take vertices $1$ and $2$, and define 
\begin{align}
V_{1,2} = \sum_{i=1}^k \1(\hat{X}_i(1)=\hat{X}_i(2)),
\end{align}
which counts how many times vertices $1$ and $2$ have been classified as being in the same community over the $k$ trials. If these graphs were independently drawn, since $\pp(\hat{X}_i(1)=\hat{X}_i(2)) = 1/2 + \e$ from the weak recovery algorithm, and as $V_{1,2}$ concentrates towards its mean, one can decide with probability $1-o(1)$ whether $1$ and $2$ are in the same community or not. The conclusion stays the same with graph-splitting as one has an approximate independence. The reason is that the graphs are in a sparse enough regime. In particular, the probability that vertex $1$ and $2$ would receive multiple edges from $k$ truly independent such SBMs is only $1-(1-p_{+})^2-2(1-p_{+})p_{+} =O(p_+) = O(\sqrt{\log(n)}/n)$, where $p_+$ is the probability of placing an edge given by $p_+=(a+b)\sqrt{\log(n)}/(2n)$. 
\item {\it Sphere comparison.} While previous argument cuts shorter if one has a weak recovery algorithm, one can also obtain almost exact recovery directly. One possibility is to count the common neighbors at large enough depth between each pair of vertices. This uses ``two poles'' for comparison rather than a ``single pole'' as used in previous section for weak recovery (where we decided for each vertex by looking at the neighbors at large depths, rather than comparing two neighborhoods).
We found that it is simpler to use a single pole when having to work at very large depth as needed for the sparse regime of weak recovery, whereas if one has the advantage of having diverging degrees, and thus the possibility of working at shorter depth, then using two poles allows for simplifications. The name ``sphere comparision'' used in \cite{colin1focs} refers to the fact that one compares the ``spheres'' of two vertices, i.e., the neighbors at a given depth from each vertex. This goes for example with the general intuition that the social spheres of two like-minded people should be more similar. In particular, the count of common neighbors is a natural benchmark of comparison. 

The depth at which spheres need to be compared needs to be above half the graph diameter, so that spheres can overlap. However, in contrast to the constant degree regime, diverging degrees allow us to compare spheres at depths below the diameter, circumventing the use of walks.   
In \cite{colin1focs}, graph splitting is also used to inject independence in the comparisons of the sphere, as the direct count of common neighbors is a challenging quantity to analyze due to dependencies. Instead, \cite{colin1focs} graph-splits the original graph into a work-graph and a bridge-graph, counting how many edges from the bridge-graph connect two spheres in the work-graph. We next provide more formal statements about this approach.
\end{itemize}
\begin{definition}\label{def-n1}
For any vertex $v$, let $N_{r[G]}(v)$ be the set of all vertices with shortest path in $G$ to $v$ of length $r$. We often drop the subscript $G$ if the graph in question is the original SBM. 
%We also refer to $\bar{N}_{r}(v)$ as the vector whose $i$-th entry is the number of vertices in $N_{r}(v)$ that are in community $i$.
\end{definition}

For an arbitrary vertex $v$ and reasonably small $r$, there will typically be about $d^r$ vertices in $N_r(v)$ (recall $d=(a+b)/2$), and about $(\frac{a-b}{2})^r$ more of them will be in $v$'s community than in each other community. Of course, this only holds when $r< \log n/\log d$ because there are not enough vertices in the graph otherwise. The obvious way to try to determine whether or not two vertices $v$ and $v'$ are in the same community is to guess that they are in the same community if $|N_r(v)\cap N_{r}(v')|>d^{2r}/n$ and different communities otherwise. Unfortunately, whether or not a vertex is in $N_r(v)$ is not independent of whether or not it is in $N_{r}(v')$, which compromises this plan. This is why we use the graph-splitting step: 
Randomly assign every edge in $G$ to some set $E$ with a fixed probability $c$ and then count the number of edges in $E$ that connect $N_{r[G\backslash E]}$ and $N_{r'[G\backslash E]}$:

\begin{definition}
For any $v, v'\in G$, $r,r'\in \mathbb{Z}$, and subset of $G$'s edges $E$, let $N_{r,r'[E]}(v\cdot v')$ be the number of pairs $(v_1,v_2)$ such that $v_1\in N_{r[G\backslash E]}(v)$, $v_2\in N_{r'[G\backslash E]}(v')$, and $(v_1,v_2)\in E$.
\end{definition}

\begin{figure}[h]
\begin{center}
\includegraphics[scale=.42]{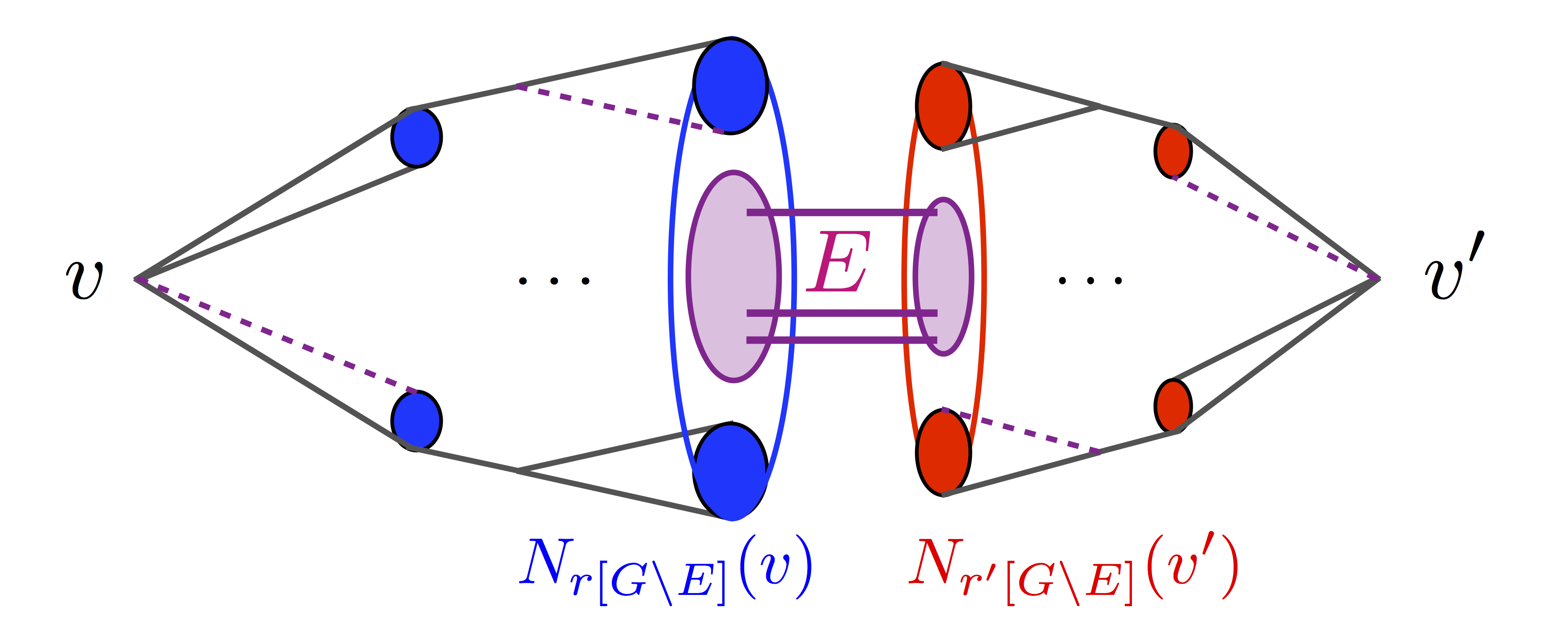}
\caption{Sphere comparison: The algorithm takes a graph-splitting of the graph with a constant probability, and decides whether two vertices are in the same community or not based on the number of crossing edges (in the first graph of the graph-split) between the two neighborhoods' spheres at a given depth of each vertices (in the second graph of the graph-split). A careful (unbalanced) choice of $r,r'$ allows us to reduce the complexity of the algorithm, but in general, $r=r'=\frac{3}{4}\log n/\log d$ suffices for the algorithm to succeed (where $d$ is the average degree).}
\label{sphere_plot}
\end{center}
\end{figure}

Figure \ref{sphere_plot} provides an illustration of the statistics used. 
Further, \cite{colin2nips} develops an invariant statistics that does not require the knowledge of the model parameters in order to compare the spheres:
\begin{definition} 
Let $G$ be a graph and let $E$ be the edge set obtained by sampling each edge with probability $c$ (i.e., the graph-split). For two vertices $v,v'$ in $G$, define the {\bf sign-invariant statistics} as 
\begin{align}
&I_{r,r'[E]} (v\cdot v'):= N_{r+2,r'[E]}(v\cdot v')\cdot N_{r,r'[E]}(v\cdot v')-N^2_{r+1,r'[E]}(v\cdot v').
\end{align}
\end{definition}
The key property is that this statistics $I_{r,r'[E]}$ with high probability scales as 
\begin{align}
\frac{c^2(1-c)^{2r+2r'+2}}{n^2}\cdot\left(d-\frac{a-b}{2}\right)^2\cdot d^{r+r'+1}\left(\frac{a-b}{2}\right)^{r+r'+1}(2\delta_{X_v,X_{v'}}-1).
\end{align}
In particular, for $r+r'$ odd, $I_{r,r'[E]} (v\cdot v')$ will tend to be positive if $v$ and $v'$ are in the same community and negative otherwise, irrespective of the specific values of $a,b$. That suggests the following agnostic algorithm for partial recovery:
\begin{algorithm}
{\it Agnostic-sphere-comparison.}
%Assume knowledge of $\delta>0$ such that $\min_{i\in [k]} p_i  \ge \delta$ and 
Input: an $n$-vertex graph and a parameter $\tau\le 0$. 
Let $d$ be the average degree of the graph:
\begin{enumerate}
\item Set $r=r'=\frac{3}{4}\log n/\log d$ and put each of the graph's edges in $E$ with probability $c=1/10$.
%\item Set $k_{\max}=1/\delta$ and select $k_{\max}\log(4k_{\max})$ random vertices, $v_1,...,v_{k_{\max}\log(4k_{\max})}$.
\item Among distinct vertices $u_1,u_2$ such that $I_{r,r'[E]} (u_1\cdot u_2)\le \tau$, pick two of maximal sum-degree, and assign each of these to a different community. Assign each vertex $u \notin \{u_1,u_2\}$ to the community $i \in \{1,2\}$ maximizing $I_{r,r'[E]} (u\cdot u_i)$. 
\end{enumerate}
\end{algorithm}

A variant of this algorithm (that applies to the general SBM) is shown in \cite{colin2nips} to solve almost exact recovery efficiently under the conditions of Theorem \ref{almost_exact}. 
We conclude this section by noting that one can also study more specific almost exact recovery requirements, allowing for a specified number of misclassified vertices $s(n)$. This is investigated in \cite{prout2} when $s(n)$ is moderately small (at most logarithmic), with an extension of Theorem \ref{almost_exact} that applies to this more general setting. The case where $s(n)$ is linear, i.e., a constant fraction of errors, is more challenging and is discussed in the next sections. 
%As long as $s(n)$ decays moderately slow, techniques above apply with minor modifications. 

\section{Partial recovery at finite SNR}\label{partial2}
Recall that partial recovery refers to the case that a fraction of misclassified vertices is constant, whereas the previous section investigates the case that a fraction of misclassified vertices is vanishing.

In the symmetric $\ssbm(n,2,a/n,b/n)$, the regime for partial recovery takes place when the following notion of SNR is finite: 
\begin{align}
\snr:=\frac{(a-b)^2}{2(a+b)}= O(1).
\end{align} 
This takes place under two circumstances:
\begin{itemize}
\item[I.] If $a,b$ are constant, i.e., the constant degree regime, 
\item[II.] If $a,b$ are functions of $n$ that diverge such that the numerator and denominator in $\snr$ scale proportionally.
\end{itemize}
Our main goal is to identify the optimal tradeoff between $\snr$ and the fraction of misclassified vertices, or between $\snr$ and the MMSE or the mutual information of the reconstruction. The latter has particular application to the compression of graphs \cite{abbe_allerton16,asadi1}. We first mention some bounds. 

Upper bounds on the fraction of incorrectly recovered vertices were demonstrated,
among others, in \cite{colin1focs,prout,new-vu,harrison}, taking form $C\exp(-c\snr)$ when $\snr$ is large. 
A bound that applies to the general SBM with arbitrary connectivity matrix $W=Q/n$ is also provided in \cite{colin1focs}.
In \cite{prout}, a spectral algorithm is shown to reach an upper-bound
of $C\exp\{-\snr/2\}$ for the two-community symmetric case and
in a suitable asymptotic sense.  An upper bound of the form 
$C\exp(-\snr/4.1)$---again for a spectral algorithm---was obtained earlier in \cite{new-vu}. 
Further, \cite{harrison} also establishes minimax optimal rate of $C\exp\{-\snr/2\}$ in the case of large $\snr$ and for certain types of SBMs, further handling a growing number of communities (to the expense of looser bounds). 

The optimal fraction of nodes that can be recovered was obtained in \cite{mossel2} for two symmetric communities when the degrees are constant but the SNR is sufficiently large, connecting to the broadcasting problem on tree problem 
\cite{evans}. 
%The reconstruction problem on tree consist in broadcasting a bit on a Galton-Watson tree with Poisson($(a+b)/2$) offspring and with binary symmetric channels of bias $b/(a+b)$ on each branch.
%That is, the probability of recovering the bit correctly from the leaves at large depth allows to determine the fraction of nodes that can be correctly labeled in the SBM in this regime. 
This result is further discussed below. 
It remains open to establish such a result at arbitrary finite SNR.  %This gives an exact characterization of the optimal fraction of misclassified vertices when the SNR is large enough. 

We next describe a result that gives the optimal tradeoff  between the SNR and the MMSE (or the mutual information) for the two-symmetric SBM in the second regime, where SNR is finite (and arbitrary) but where degrees diverge. After that, we discuss results for constant degrees but large enough SNR. 

\section{Mutual Information-SNR tradeoff}\label{wigner}
In this section, we study the finite SNR regime with diverging degrees and show that the SBM is essentially equivalent to a spiked Wigner model (i.e., a low-rank matrix perturbed by a Wigner random matrix), where the spiked signal has a block structure (rather than a sparse structure as in sparse PCA \cite{rigollet,yash_pca}).
To compare the two models, we use the mutual information.  
 
In this section we use $p_n,q_n$ in lieu of $q_{\mathrm{in}},q_{\mathrm{out}}$, i.e., the inner- and outer-cluster probabilities. 
For $(X,G) \sim \ssbm(n,2,p_n,q_n)$, the mutual information of the SBM is $I(X;G)$, where 
\begin{align*}
I(X;G)=H(G)-H(G|X)=H(X)-H(X|G),
\end{align*}
and $H$ denotes the entropy. We next introduce the normalized MMSE of the SBM:
\begin{align}
\MMSE_n(\snr)&\equiv
\frac{1}{n(n-1)}\E\Big\{\big\| XX^{\sT}-\E\{XX^{\sT}|G\}\big\|_F^2\Big\}\, \label{eq:MMSEdef}\\
&= \min_{\hx_{12}: \cG_n\to \reals}
\E\big\{\big[X_1 X_2-\hx_{12}(G)\big]^2\big\}\,  \label{eq:MMSE_2vert_opt}
\end{align}
where $\| \cdot \|_F$ denotes the Frobenius norm. 

To state the result that provides a single-letter
characterization of the per-vertex MMSE (or mutual information), we need to introduce the
\emph{effective Gaussian
scalar channel}.
Namely, define the Gaussian channel
\begin{align}
  Y_0 = Y_0(\gamma) = \sqrt{\gamma} \, X_0 + Z_0, \label{eq:scalarprob_1}
\end{align}
where $X_0\sim \Unif(\{+1,-1\})$ independent of $Z_0\sim\normal(0,1)$.
We denote by $\mmse(\gamma)$ and $\Info(\gamma)$ the corresponding
minimum mean square error and
mutual information:
\begin{align}
\Info(\gamma) & = \E\,\log \Big\{\frac{\de p_{Y|X}(Y_0(\gamma)|X_0)}{\de p_{Y}(Y_0(\gamma))}\Big\}\, ,\label{eq:InfoDef}\\
  \mmse(\gamma ) &= \E\left\{ ( X_0- \E\left\{ X_0| Y_0(\gamma)\right\})^2
                   \right\}\, . \label{eq:mmseDef}
\end{align}
Note that these quantities can be written explicitly as
Gaussian integrals of elementary functions:
\begin{align}
\Info(\gamma) & = \gamma-\E\log\cosh\big(\gamma+\sqrt{\gamma}\,
Z_0\big)\, ,\label{eq:InfoFormula}\\
\mmse(\gamma ) &= 1- \E\big\{\tanh(\gamma+\sqrt{\gamma}\, Z_0)^2\big\}\, . \label{eq:mmseFormula}
\end{align}

We are now in position to state the result.

\begin{theorem}\label{thm:main}\cite{yash_sbm}
  For any $\lambda>0$, let $\gamma_* = \gamma_*(\lambda)$ be the
  largest non-negative solution of the equation
  \begin{align}
      \gamma &= \lambda\big(1- \mmse(\gamma)\big)\,  \label{eq:MainEquation}
  \end{align}
and
%We refer to $\gamma_*(\lambda)$ as to the \emph{effective
%  signal-to-noise ratio}.
%  Further, define $\Psi(\gamma, \lambda)$ by:
%
  \begin{align}
    \Psi(\gamma, \lambda ) &=
    \frac{\lambda}{4}+\frac{\gamma^2}{4\lambda}-\frac{\gamma}{2}+\Info(\gamma)\, .
\end{align}
%

%Recall that for $\ssbm(n,2,p_n,q_n)$ we have\footnote{In fact, if one wants to also consider denser regimes where $p+q$ tends to 1, then one should use the general definition of $\snr:=\frac{n(p-q)^2}{4(p+q)/2(1-(p+q)/2)}$.} 
%\begin{align}
%\snr:=\frac{n(p-q)^2}{2(p+q)}.
%\end{align}  
\noindent
Let $(X,G)\sim \ssbm(n,2,p_n,q_n)$ and define\footnote{Note that this is asymptotically the same notion of SNR as defined earlier when $p_n,q_n$ vanish.} $\snr:= n\,
(p_n-q_n)^2/(2 (p_n+q_n) (1-(p_n+q_n)/2)).$    
Assume that, as $n\to\infty$, $(i)$ $\snr\to\lambda$ and $(ii)$  $n(p_n+q_n)/2(1-(p_n+q_n)/2)\to\infty$.
Then, 
\begin{align}
\lim_{n\to\infty} \MMSE_n(\snr) &=
1-\frac{\gamma_*(\lambda)^2}{\lambda^2}\, \\
    \lim_{n\to\infty}  \frac{1}{n}\, I(X;G)&= \Psi(\gamma_*(\lambda),
    \lambda)\, .
\end{align}

Further, this implies $\lim_{n\to\infty} \MMSE_n(\snr) = 1$ for
$\lambda\le 1$ (i.e., weak recovery unsolvable) and $\lim_{n\to\infty} \MMSE_n(\snr) < 1$ for
$\lambda>1$ (i.e., weak recovery solvable). 

\end{theorem}

\begin{corollary}\label{corr:sparsedegbound}\cite{yash_sbm}
    When $p_n = a/n, q_n = b/n$, where $a, b$ are bounded as $n$ diverges, 
    there exists an absolute constant $C$ such that
\begin{align}
  \limsup_{n\to\infty}
  \Big|\frac{1}{n}I(X;G)-\Psi(\gamma_*(\lambda),\lambda) \Big| \le \frac{C\lambda^{3/2} }{\sqrt{a+b}}\,. 
\end{align}
Here $\lambda$, $\psi(\gamma, \lambda)$ and $\gamma_*(\lambda)$ are as in Theorem 
\ref{thm:main}.
\end{corollary}
\noindent
A few remarks about the previous theorem and corollary:
\begin{itemize}
\item Theorem \ref{thm:main} shows that the normalized MMSE (or mutual information) is non-trivial if and only if $\lambda>1$. This extends the results on weak recovery \cite{massoulie-STOC,Mossel_SBM1} discussed in Section \ref{weak} for the regime of finite SNR with diverging degrees, completing weak recovery in the SSBM with two communities and any choice of parameters;
\item Theorem \ref{thm:main} also gives upper and lower bounds for the optimal agreement. Let $$\overlap_n(\snr)= \frac{1}{n}\sup_{\hbs: \cG_n\to
  \{+1,-1\}^n}\E\big\{|\<X,\hbs(G)\>|\big\}.$$ Then, 
  %Note that by returning $\hs_i(G) \in\{+1,-1\}$ uniformly at random, we have $\E\big\{|\<X,\hbs(G)\>|\big\}/n = O(n^{-1/2})\to 0$. 
\begin{align}
    &1-\MMSE_n(\snr) + O(n^{-1}) \le \overlap_n(\snr) \\
    &\le \sqrt{1-\MMSE_n(\snr)} + O(n^{-1/2}).
\end{align}
\item In \cite{mossel-xu}, tight expressions similar to those obtained in Theorem \ref{thm:main} for the MMSE are obtained for the optimal expected agreement with additional scaling requirements. Namely, it is shown that for $\ssbm(n,2,a/n,b/n)$ with $a=b + \mu \sqrt{b}$ and $b=o(\log n)$, the least fraction of misclassified vertices is in expectation given by $Q(\sqrt{v^*})$ where $v^*$ is the unique fixed point of the equation $v=\frac{\mu^2}{4} \E \tanh(v + v\sqrt{Z})$, $Z$ is normal distributed, and $Q$ is the Q-function for the normal distribution. Similar expressions were also derived in \cite{zhang_moore} for the overlap metric, and \cite{lenka_rank} for the MMSE. 
%\item The above results are information-theoretic, and use AMP as a proof technique to obtain tight bounds on the MMSE rather than showing optimality of the algorithm. It is however believed that AMP would also achieve the optimal tradeoff. 
\item Note that Theorem \ref{thm:main} requires merely diverging degrees (arbitrarily slowly), in contrast to general results from random matrix theory such as \cite{peche} that would require poly-logarithmic degrees to extract communities from the spiked Wigner model point of view. We refer to \cite{wein1,wein_tensor,banks3} and references therein for generalizations of the spiked Wigner model discussed here with more general input signals.  
\end{itemize}

\begin{figure}[h]
\begin{center}
\includegraphics[scale=.4]{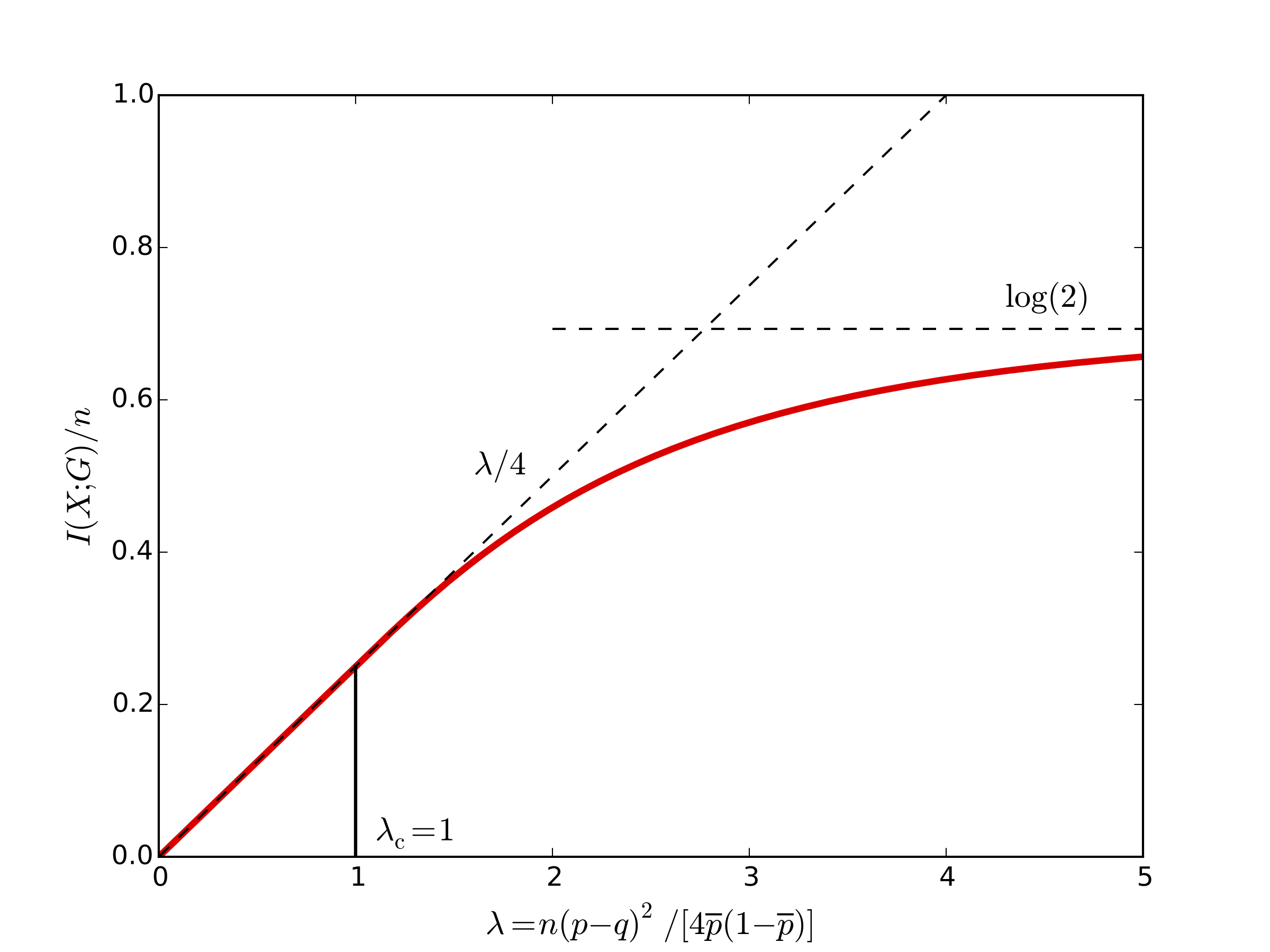}
\caption{Asymptotic mutual information per vertex of the
   symmetric stochastic block model with two communities, as a function of the
   signal-to-noise ratio $\lambda$. The dashed lines are simple upper
   bounds: $\lim_{n\to\infty} I(X;G)/n\le \lambda/4$ and $I(X;G)/n\le \log 2$. }
\label{i-mmse-plot}
\end{center}
\end{figure}

\section{Proof technique and connections to spiked Wigner models}
Theorem \ref{thm:main} gives an exact expression for the normalized MMSE and mutual information in terms of an effective Gaussian noise channel. The Gaussian distribution emerges due to a universality result established in the proof: in the regime of the theorem, the SBM model is equivalent to a spiked Wigner model given by  
$$Y= \sqrt{\lambda/n} XX^t + Z$$
where $Z$ is a Wigner random matrix (i.e., symmetric with i.i.d.\ Normal entries), and where we recall that $\lambda$ corresponds to the limit of $\snr$.

The formal statement of the equivalence is as follows:

\begin{theorem}[Equivalence of SBM and a spiked Wigner model]\label{prop:gausseq}
Let $I(X;G)$ be the mutual information of $\ssbm(n,2,p_n,q_n)$ with $\snr\to\lambda$ and $n(p_n+q_n)/2(1-(p_n+q_n)/2)\to\infty$, and $I(X;Y)$ be the mutual information for spiked Wigner model $Y= \sqrt{\lambda/n} XX^t + Z$. Then, there is a constant $C$ independent of $n$ such that
  \begin{align}
    &\frac{1}{n} \big|I(X;G) - I(X; Y)\big| \notag \\
    &\le C\left(
    \frac{\lambda^{3/2}}{\sqrt{n(p_n+q_n)/2(1-(p_n+q_n)/2)}} + |\snr-\lambda|
    \right).
  \end{align}
\end{theorem}

To obtain the limiting expression for the normalized mutual information in Theorem \ref{thm:main}, first notice that for $Y(\lambda)= \sqrt{\lambda/n} XX^t + Z$,
\begin{align*}
\frac{1}{n}I(X;Y(0))=0 \qquad \frac{1}{n}I(X;Y(\infty))=\log(2).
\end{align*}
Next, (i) use the fundamental theorem of calculus to express these boundary conditions as an integral of the derivative of the mutual information, (ii) use the I-MMSE identity \cite{i-mmse} to express this derivative in terms of the MMSE, (iii) upper-bound the MMSE using a specific estimator obtained from the AMP algorithm \cite{amp} (or any algorithm performing optimally in this regime), (iv) evaluate the asymptotic performance of the AMP estimate using the density evolution technique \cite{bayati,yash_pca}, and (v) note that the obtained bound matches the original value of $\log(2)$ in the limit of $n$ tending to infinty: 
\begin{align}
\log(2)&\stackrel{(i)}{=} \frac{1}{n}\int_{0}^{\infty} \frac{\partial}{\partial \lambda} I(XX^t;Y(\lambda))\, d \lambda \\
&\stackrel{(ii)}{=} \frac{1}{4n^2}   \int_{0}^{\infty} \mathrm{MMSE}(XX^t|Y(\lambda)) \, d \lambda \\
& \stackrel{(iii)}{\le } \frac{1}{4n^2}   \int_{0}^{\infty} \mathbb{E}(XX^t -\hat{x}_{\mathrm{AMP},\lambda}(\infty)\hat{x}^t_{\mathrm{AMP},\lambda}(\infty) )^2 \, d \lambda \\
& \stackrel{(iv)}{=} \Psi(\gamma_*(\infty), \infty)-\Psi(\gamma_*(0), 0) + o_n(1)\\
& \stackrel{(v)}{=}\log(2) + o_n(1).
\end{align}
This implies that (iii) is in fact an equality asymptotically, and using monotonicity and continuity properties of the integrand, the identity must hold for all SNRs as stated in the theorem. The only caveat not discussed here is the fact that AMP needs an initialization that is not fully symmetric to converge to the right solution, which causes the insertion in the proof of a noisy genie on the true labels $X$ at the channel output to break the symmetry for AMP. The genie is then removed by taking noise parameters that are arbitrarily large.

%In contrast, Theorem \ref{thm:MainEstimation} provides the exact asymptotic ($n\to\infty$) expression for the estimation error for any fixed $\lambda$. Equivalently, for any pre-assigned estimation error $\delta\in (0,1)$,  Theorem \ref{thm:MainEstimation} characterizes a sharp threshold in signal-to-noise ratio $\lambda_*(\delta)$ above which $\lim_{n\to\infty}\MMSE_n(\lambda)<\delta$. Further, the expression in Theorem \ref{thm:MainEstimation} (together with some straightforward analysis of the function $\gamma_*(\lambda)$) implies $\lim_{n\to\infty}\MMSE_n(\lambda) = \exp\{-\lambda^2/2+o(\lambda^2)\}$ for large $\lambda$, and therefore  the minimax optimality of the rate $C\exp\{-\lambda^2_n/2\}$ in this regime In independent work, Zhang and Zhou \cite{zhang2015minimax} also establish that this rate is minimax optimal albeit in the  slightly different regime $\snr\to\infty$ with $n$, and using fairly different techniques.

\section{Partial recovery for constant degrees}\label{boost}
Obtaining the expression for the optimal agreement at finite and arbitrary SNR when the degrees are constant remains an open problem (see also Sections \ref{almost} and \ref{open}). The problem is settled for high enough SNR in \cite{mossel2}, with the following expression relying on reconstruction error for the broadcasting on trees problem. 

Define the optimal agreement fraction as 
\begin{align}
P_{G_n} (a,b):=\frac{1}{2} + \sup_{f} \E |\frac{1}{n} \sum_{v} \1(f(v,G_n)=X_v) - \frac{1}{2}| .
\end{align}
Note that the above expression takes into account the symmetry of the problem and can also been interpreted both as a normalized agreement and probability. Let $P_G(a,b) := \limsup_g P_{G_n} (a,b)$.
Define now the counterpart for the broadcasting problem on tree: Back to the notation of Section \ref{bot},
define $T^{(t)}$ as the Galton-Watson tree with Poisson$((a+b)/2)$ offspring, flip probability $b/(a+b)$ and depth $t$, and define the optimal inference probability of the root as  
\begin{align}
P_{T} (a,b):=\frac{1}{2} + \lim_{t \to \infty} \E | \E(X^{(0)}| X^{(t)}) - 1/2| .
\end{align}
The reduction from \cite{Mossel_SBM1} discussed in Section \ref{bot} allows to deduce that $P_{G}(a,b)\le P_T(a,b)$, and this is shown to be an equality for large enough SNR:
\begin{theorem}\cite{mossel2}
There exists $C$ large enough such that if $\snr >C$ then $P_{G}(a,b) = P_T(a,b)$, and this normalized agreement is efficiently achievable.
\end{theorem}
The theorem in \cite{mossel2} has a weaker requirement of $\snr >C \log (a+b)$, but later developments on weak recovery imply for free the version stated above. 
Note that $P_T(a,b)$ gives an implicit expression for the optimal fraction, though it admits a variational representation due to \cite{tree-andrea}. The efficient algorithm is a variant of belief propagation.

In \cite{decelle}, it is conjectured that BP gives the optimal agreement at all SNR. However, as discussed in Section \ref{bit-map}, BP is hard to analyze in the context of loopy graphs with a random initialization. Another strategy is to proceed with a two-round procedure, which is used to establish the above results in \cite{mossel2} for two communities. 
The idea is to use a simpler algorithm to obtain a non-trivial reconstruction when $\snr >1$, see Section \ref{weak}, and then to improve the accuracy using full BP at shorter depth. To show that the accuracy achieved is optimal, one also has to show that a noisy version of the reconstruction on tree problem \cite{robust_rec}, where leaves do not have exact labels but noisy labels, leads to the same probability of error at the root. This is expected to take place for two communities at all SNR above the KS threshold, and it was shown in \cite{mossel2} for the case of large enough SNR. 
This type of claim is not expected to hold for general $k$. For more than two communities, one needs to first convert the output of the algorithm discussed in Section \ref{linBP}, which gives two sets that correlated with their communities, into a nontrivial assignment of a belief to each vertex; this is discussed in \cite{colin3cpam}.  
Then one can use these beliefs as starting probabilities for a belief propagation algorithm of depth $\log(n)/3\log(\lambda_1)$, which now runs on a tree-like graph. %Again, this part preserves optimal accuracy if the reconstruction problem on a tree with noisy leaves admits the same fixed point distribution than the noiseless case, such as for two or three communities.  

\chapter{The general SBM}\label{general}
In this section we discuss results for the general SBM, where communities can take arbitrary relative sizes and where connectivity rates among communities are arbitrary. 
\section{Exact recovery and CH-divergence}\label{exact}
We provide the fundamental limit for exact recovery in the general SBM, in the regime of the phase transition where $W$ scales like $\log(n)Q/n$ for a matrix $Q$ with positive entries. 
\begin{theorem}\cite{colin1focs}\label{exact_thm}
Exact recovery in $\mathrm{SBM}(n,p,\log(n)Q/n)$ is solvable and efficiently so if
\begin{align*}
&I_+(p,Q):=\min_{1 \leq i<j \leq k}  D_+( (\diag(p)Q)_i \| (\diag(p)Q)_j) > 1
\end{align*}
and is not solvable if $I_+(p,Q)<1$, where $D_+$ is defined by 
\begin{align}
D_+(\mu\|\nu):=\max_{t \in [0,1]} \sum_x \nu(x) f_t (\mu(x)/\nu(x)), \quad f_t(y):=1-t+ty-y^t.
\end{align}
\end{theorem}
\begin{remark}
Regarding the behavior at the threshold: 
If all the entries of $Q$ are non-zero, then exact recovery is solvable (and efficiently so) if and only if $I_+(p,Q) \ge 1$. 
In general, exact recovery is solvable at the threshold, i.e., when $I_+(p,Q) = 1$, if and only if any two columns of $\diag(p)Q$ have a component that is non-zero and different in both columns.  
\end{remark}
\begin{remark}
 In the symmetric case $\ssbm(n,k,a \log (n)/n,b \log (n)/n)$, the CH-divergence is maximized at the value of $t=1/2$, and it reduces in this case to the Hellinger divergence between any two columns of $Q$; the theorem's inequality becomes 
$$\frac{1}{k} (\sqrt{a} - \sqrt{b})^2 >1,$$
matching the expression obtained in Theorem \ref{exact1} for two symmetric communities.
\end{remark}

We now discuss some properties of the functional $D_+$ governing the fundamental limit for exact recovery in Theorem \ref{exact_thm}. For $t \in [0,1]$, let
\begin{align}
D_t(\mu\|\nu):= \sum_x \nu(x) f_t (\mu(x)/\nu(x)), \quad f_t(y)=1-t+ty-y^t,
\end{align}
and note that $D_+=\max_{t\in [0,1]} D_t $. Since the function $f_t$ satisfies
\begin{itemize}
\item $f_t(1)=0$
\item $f_t$ is convex on $\mR_+$,
\end{itemize}
the functional $D_t$ is a so-called $f$-divergence \cite{csiszar-f}, like the KL-divergence ($f(y)=y \log y$), the Hellinger divergence, or the Chernoff divergence. Such functionals have a list of common properties described in \cite{csiszar-f}. For example, if two distributions are perturbed by additive noise (i.e., convolving with a distribution), then the divergence always increases, or if some of the elements of the distributions' support are merged, then the divergence always decreases. Each of these properties can be interpreted in terms of community detection (e.g., it is easier to recovery merged communities, etc.). Since $D_{t}$ collapses to the Hellinger divergence when $t=1/2$ and since it matches the Chernoff divergence for probability measures, we call $D_t$ (and $D_+)$ the Chernoff-Hellinger (CH) divergence in \cite{colin1focs}.

Theorem \ref{exact_thm} thus gives a new operational meaning to an $f$-divergence, showing that the fundamental limit for data clustering in SBMs is  governed by the CH-divergence, similarly to the fundamental limit for data transmission in DMCs being governed by the KL-divergence. If the columns of $\diag(p)Q$ are ``different'' enough, where difference is measured in CH-divergence, then one can separate the communities. This is analogous to the channel coding theorem that says that when the output's distributions are ``different'' enough, where difference is measured in KL-divergence, then one can separate the codewords.

\subsection{Converse}\label{genie} 
%\subsection{The Maximum A Posteriori (MAP) estimator}
Let $(X,G)\sim \sbm(n,p,W)$. 
Recall that to solve exact recovery, we need to find the partition of the vertices, but not necessarily the actual labels. Equivalently, the goal is to find the community partition $\Omega=\Omega(X)$ as defined in Section \ref{model}.  
Recall also that the optimal estimator (see Section \ref{exact-2}) is the MAP estimator $\hat{\Omega}_{\map}(\cdot)$ that maximizes the posterior distribution 
\begin{align}
\pp \{ \Omega= s| G=g \}, \label{ap}
\end{align}
or equivalently 
\begin{align}
\sum_{x \in [k]^n: \Omega(x)=s} \pp \{ G= g | X=x \}  \prod_{i=1}^k p_{i}^{|\Omega_i(x)|}, \label{posterior}
\end{align}
and any such maximizer can be chosen arbitrarily. 
%It is worth simplifying these expressions for the the two-symmetric SBM, i.e., $\ssbm(n,2,A,B)$. In this case, denoting by $N_{in}$ and $N_{out}$ the number of edges inside and across clusters respectively, and $w_1=|\Omega_1(x)|$, $w_2=|\Omega_2(x)|$, we have  
%\begin{align}
%\pp\{  G=g |X=x\} &=  A^{N_{in}} (1-A)^{{w_1 \choose 2} +{w_2 \choose 2} -N_{in}} B^{N_{out}} (1-B)^{w_1w_2-N_{out}} \\
%&\propto  \left( \frac{B(1-A)}{A(1-B)} \right)^{N_{out}} \cdot \1(x \text{ is balanced})
%\end{align}
%where $N_{in}$ is the number of edges that $G$ has inside the clusters defined by $x$, and $N_{out}$ is the number of crossing edges. 
If MAP fails in solving exact recovery, no other algorithm can succeed. 
%Note that to succeed for exact recovery, the partition shall be typical in order to make the last factor in \eqref{posterior} non-vanishing (i.e., communities of relative size $p_i+o(1)$ for all $i \in [k]$). Of course, resolving exactly the maximization in \eqref{ap} requires comparing exponentially many terms, so the MAP estimator may not always reveal the computational threshold for exact recovery.

%Recall also that $X$ has with probability exponentially close to one communities of relative size $p_i+o(1)$ for all $i \in [k]$. So the MAP estimator can only succeed in solving exact recovery if it outputs communities of such relative sizes, and the search can be restricted to such typical vectors $x$.

%\subsection{MAP for SSBM}

%\subsection{Converse: the genie-aided approach}\label{genie}
We proceed similarly to the symmetric case to obtain the impossibility part of Theorem \ref{exact_thm}, i.e., we reduce the problem to a genie hypothesis test for recovering a single vertex given the other vertices. However, we now work in the Bernoulli community prior model, for slight conveniences.

Imagine that in addition to observing $G$, a genie provides the observation of $X_{\sim u}=\{X_v: v \in [n] \setminus \{u\}\}$. 
Define now $\hat{X}_v=X_v$ for $v \in [n] \setminus \{u\}$ and 
\begin{align}
\hat{X}_{u,\map}(g,x_{\sim u})=\arg \max_{i \in [k]} \pp\{ X_u=i | G=g,X_{\sim u}=x_{\sim u}  \},
\end{align}  
where ties can be broken arbitrarily if they occur (we assume that an error is declared in case of ties to simplify the analysis). 
%Recall that 
%\begin{align}
%P_e:=\pp\{A(\hat{X}_{\map}(G), X) \neq 1 \} &= \pp\{\forall \pi \in S_k, \exists u \in [n] : (\hat{X}_{\map}(G))_u \neq \pi(X_u) \} ,
%\end{align}
%that is, we fail at exact recovery if no matter how we label the communities, there is one component that is misclassified, hence
If we fail at recovering a single component when all others are revealed, we must fail at solving exact recovery all at once, thus 
\begin{align}
\pp\{ \hat{\Omega}_{\map}(G) \neq \Omega \} &\ge  \pp\{ \exists u \in [n] : \hat{X}_{u,\map}(G,X_{\sim u}) \neq X_u \} . \label{geniebound}
\end{align}
This lower bound may appear to be loose at first, as recovering the entire communities from the graph $G$ seems much more difficult than classifying each vertex by having all others revealed (we call the latter component-MAP). As shown for the two symmetric case in Section \ref{exact-2}, the obtained bound is, however, tight. %In any case, studying when this lower bound is not vanishing always provides a necessary condition for exact recovery. 

Let $E_u := \{ \hat{X}_{u,\map}(G,X_{\sim u}) \neq X_u \}$. If the events $E_u$ were independent, we could write $\pp\{ \cup_u E_u\} = 1-\pp\{\cap_u E_u^c \} = 1- (1-\pp\{E_1\})^{n} \geq 1- e^{-n \pp\{E_1\}} $ and if $\pp\{E_1\}= \omega(1/n)$, this would drive $\pp\{ \cup_u E_u\}$, and thus $P_e$, to 1. The events $E_u$ are not independent, but their dependencies are weak enough  that previous reasoning still applies, and $P_e$ is driven to 1 when $\pp\{E_1\}=\omega(1/n)$. In Section \ref{exact-2}, we used the second moment method to obtain this statement, showing that the events are asymptotically independent, which we also pursue below. In \cite{colin1focs}, a variant of this method is used to obtain the conclusion. 

Recall that for the second moment method, one defines
$$Z=\sum_{u \in [n]} \1( \hat{X}_{u,\map}(G,X_{\sim u}) \neq X_u) ,$$
which counts the number of components where component-MAP fails.  
Note that the right hand side of \eqref{geniebound} corresponds to $\pp\{ Z \ge  1\}$ as desired. Our goal is to show that $\frac{\Var Z}{(\E Z)^2}$ stays strictly below 1 in the limit, or equivalently, $\frac{\E Z^2}{(\E Z)^2}$ stays strictly below 2 in the limit. In fact, the latter tends to $1$ in the converse of Theorem \ref{exact_thm}.

Note that $Z= \sum_{u \in [n]} Z_u$ where $Z_u:= \1( \hat{X}_{u,\map}(G,X_{\sim u}) \neq X_u)$ are binary random variables with $\E Z_u=\E Z_v$ for all $u,v$. 
Hence, 
%\begin{align}
%\E Z&=  n \pp\{Z_1 =1\}\\
%\E Z^2 &= \sum_{u,v \in [n]} \E (Z_u Z_v)  = \sum_{u,v \in [n]} \pp\{Z_u =Z_v=1\} \\
%&= n \pp\{Z_1 =1\} + n(n-1) \pp\{Z_1 =1\} \pp\{Z_2=1 | Z_1 =1\} 
%\end{align}
%and 
$\frac{\E Z^2}{(\E Z)^2}$ tends to 1 if 
\begin{align}
%\frac{n \pp\{Z_1 =1\} + n(n-1) \pp\{Z_1 =1\} \pp\{Z_2=1 | Z_1 =1\} }{n^2 \pp\{Z_1 =1\}^2} = 1+o(1)
%\end{align}
%or
%\begin{align}
\frac{1}{n \pp\{Z_1 =1\}} + \frac{\pp\{Z_2=1 | Z_1 =1\} }{\pp\{Z_1 =1\}} &= 1+o(1)
\end{align}
which takes place if $n \pp\{Z_1 =1\}=\omega(1)$ and $\frac{\pp\{Z_2=1 | Z_1 =1\} }{\pp\{Z_2 =1\}} = 1 +o(1)$.

%The asymptotic independence takes place due to the regime that we consider for the block model in the theorem. To given a related example, in the context of the Erd\H{o}s-R\'enyi model $ER(n,p)$, if $W_1$ is 1 when vertex $u$ is isolated and 0 otherwise, then $\pp\{W_1=1|W_2=1\}= (1-p)^{n-2}$ and $\pp\{W_1=1\}= (1-p)^{n-1}$, and thus $\frac{\pp\{W_1=1 | W_2 =1\} }{\pp\{W_1 =1\}} =(1-p)^{-1}$ tends to 1 as long as $p$ tends to 0. That is, the property of a vetex being isolated is asymptotically independent as long as the edge probability is vanishing. A similar outcomes takes place for the property of MAP-component failing when edge probabilities are vanishing in the block model.

\begin{figure}[h]
\centering
  \includegraphics[width=.93\linewidth]{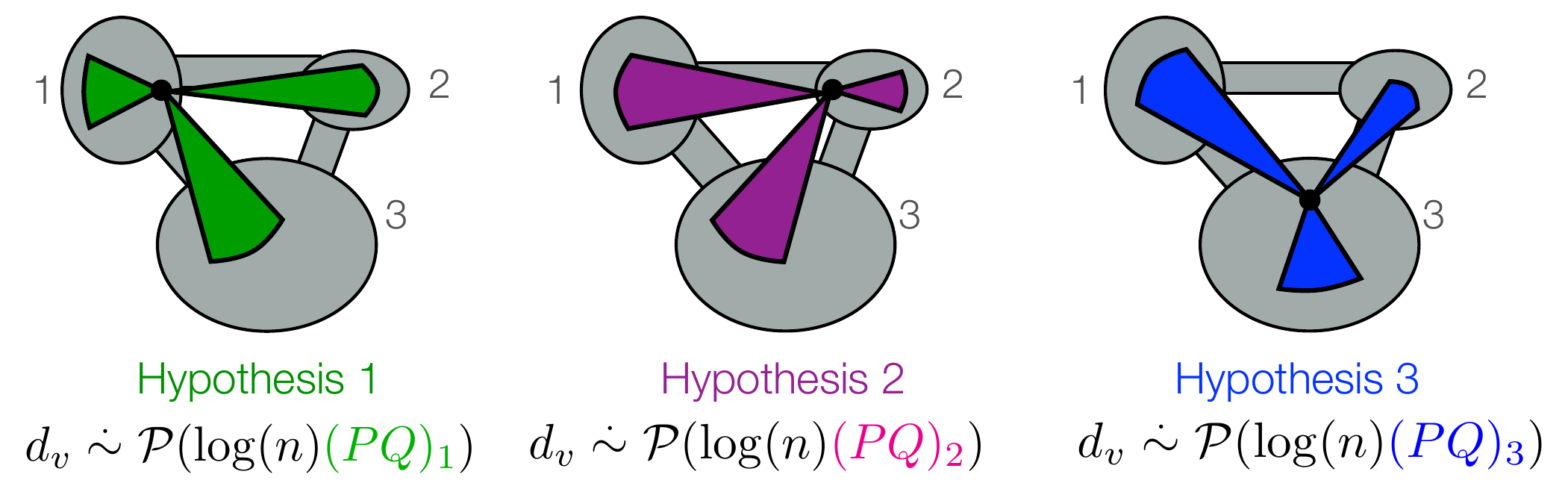}
  \caption{The genie-aided hypothesis test (degree-profiling) to classify a vertex given the labels of all other vertices consists in a multi-hypotheses test with multivariate Poisson distributions of means corresponding to the different community profiles. The probability of error of that test scales as $n^{-I_+(p,Q)}$ where $I_+(p,Q)$ is given by the CH-divergence $D_+$ between the community profiles as in Lemma \ref{overlap}.}\label{degree-pro}
\end{figure}

The location of the threshold is then dictated by requirement that $n \pp\{Z_1 =1\}$ diverges, and this where the CH-divergence threshold emerges from a moderate deviation analysis. We next summarize what we obtained with the above reasoning, and then specialized to the regime of Theorem \ref{exact_thm}.
\begin{theorem}
Let $(X,G)\sim \sbm(n,p,W)$, $E_u:=\1(\hat{X}_{u,\map}(G,X_{\sim u}) \neq X_u)$, $u \in [n]$. If the events $E_1$ and $E_2$ are asymptotically independent, then exact recovery is not solvable if 
\begin{align}
\pp\{\hat{X}_{u,\map}(G,X_{\sim u}) \neq X_u\} = \omega(1/n).
\end{align}
\end{theorem}
The next lemma gives the behavior of $\pp\{Z_1=1\}$ in the logarithmic degree regime.
%\begin{lemma}
%Let $(X,G)\sim \sbm(n,p,W)$ where $p$ is constant and $W$ satisfies . Then $Z_1$ and $Z_2$ are asymptotically independent.
%\end{lemma}
\begin{lemma}\label{overlap}\cite{colin1focs}
Consider the hypothesis test where $H=i$ has prior probability $p_i$ for $i \in [k]$, and where observable $Y$ is distributed $\bin(np,W_{i})$ under hypothesis $H=i$. This is called degree-profiling in \cite{colin1focs}, and is illustrated in Figure \ref{degree-pro}. 
Then the probability of error $P_e(p,W)$ of MAP decoding for this test satisfies 
$\frac{1}{k-1} \mathrm{Over}(n,p,W) \le P_e(p,W) \le \mathrm{Over}(n,p,W)$ where 
\begin{align*}
&\mathrm{Over}(n,p,W) = \\  
&\sum_{i < j} \sum_{z \in \mZ_+^k} \min(\pp \{\bin(np,W_{i})=z\} p_{i} ,  \pp \{\bin(np,W_{j})=z \} p_{j}), 
\end{align*}
and for a symmetric $Q \in \mR_+^{k \times k}$,
\begin{align}
\mathrm{Over}(n,p,\log(n)Q/n) =  n^{-I_+(p,Q) -O(\log \log (n)/\log n)},
\end{align}
where $I_+(p,Q) =\min_{i<j} D_+((\diag(p)Q)_i,(\diag(p)Q)_j)$.
\end{lemma}
\begin{corollary}
Let $(X,G)\sim \sbm(n,p,W)$ where $p$ is constant and $W=Q \frac{\log n}{n}$. Then
\begin{align}
\pp\{\hat{X}_{u,\map}(G,X_{\sim u}) \neq X_u\} = n^{-I_+(p,Q) + o(1)}.
\end{align}
\end{corollary}
A robust extension of this Lemma is proved in \cite{colin1focs} that allows for a slight perturbation of the binomial distributions. 

\subsection{Achievability}\label{split}
Two-round algorithms have proved to be powerful in the context of exact recovery. The general idea consists in using a first algorithm to obtain a good but not necessarily exact clustering, solving a joint assignment of all vertices, and then to switch to a local algorithm that ``cleans up'' the good clustering into an exact one by reclassifying each vertex. This approach has a few advantages:
\begin{itemize}
\item If the clustering of the first round is accurate enough, the second round becomes approximately the genie-aided hypothesis test discussed in previous section, and the approach is built in to achieve the threshold;
\item if the clustering of the first round is efficient, then the overall method is efficient since the second round only performs computations for each single node separately and has thus linear complexity.
\end{itemize}
Some difficulties need to be overome for this program to be carried out:
\begin{itemize}
\item One needs to obtain a good clustering in the first round, which is typically non-trivial;
\item One needs to be able to analyze the probability of success of the second round, as the graph is no longer independent of the obtained clusters.
\end{itemize}
To resolve the latter point, we rely in \cite{abh} on the technique which we call ``graph-splitting'' and which takes again advantage of the sparsity of the graph.
\begin{definition}[Graph-splitting]
Let $G$ be an $n$-vertex graph and $\gamma \in [0,1]$. The graph-splitting of $G$ with split-probability $\gamma$ produces two random graphs $G_1,G_2$ on the same vertex set as $G$. The graph $G_1$ is obtained by sampling each edge of $G$ independently with probability $\gamma$, and $G_2= G \setminus G_1$ (i.e., $G_2$ contains the edges from $G$ that have not been subsampled in $G_1$). 
\end{definition} 
Graph splitting is convenient in part due to the following fact. 
\begin{lemma}\label{split_ind}
Let $(X,G)\sim \sbm(n,p,\log n Q/n)$, $(G_1,G_2)$ be a graph splitting of $G$ with split-probability $\gamma$, and $(X,\tilde{G}_2) \sim \sbm(n,p,(1-\gamma)\log n Q/n)$ with $\tilde{G}_2$ independent of $G_1$. Let $\hat{X}=\hat{X}(G_1)$ be valued in $[k]^n$ such that $\pp\{A(X, \hat{X}) \ge 1-o(n) \}=1-o(1)$. For any $v \in [n]$, $d \in \mZ_+^k$, 
\begin{align}
\pp\{D_v(\hat{X},G_2) = d\} \le (1+o(1)) \pp\{D_v(\hat{X},\tilde{G}_2) = d\} + n^{-\omega(1)},
\end{align}
where $D_v(\hat{X},G_2)$ is the degree profile of vertex $v$, i.e., the $k$-dimensional vector counting the number of neighbors of vertex $v$ in each community using the clustered graph $(\hat{X},G_2)$.
\end{lemma}
The meaning of this lemma is as follows. We can consider $G_1$ and $G_2$ to be approximately independent, and export the output of an algorithm run on $G_1$ to the graph $G_2$ without worrying about dependencies to proceed with component-MAP. Further, if $\gamma$ is to chosen as $\gamma=\tau(n)/\log(n)$ where $\tau(n)=o(\log(n))$, then $G_1$ is distributed as $\sbm(n,p,\tau(n) Q/n)$ and $G_2$ remains approximately as $\sbm(n,p,\log n Q/n)$. This means that from our original SBM graph, we produce essentially `for free' a preliminary graph $G_1$ with $\tau(n)$ expected degrees that can be used to get a preliminary clustering, and we can then improve that clustering on the graph $G_2$ which has still logarithmic expected degree. 

Our goal is to obtain on $G_1$ a clustering that is almost exact, i.e., with only a vanishing fraction of misclassified vertices. If this can be achieved for some $\tau(n)$ that is $o(\log(n))$, then a robust version of the genie-aided hypothesis test described in Section \ref{genie} can be run to re-classify each node successfully when $I_+(p,Q)>1$. Luckily, as we shall see in Section \ref{almost}, almost exact recovery can be solved with the mere requirement that $\tau(n)=\omega(1)$ (i.e., $\tau(n)$ diverges). In particular, setting $\tau(n)=\log \log(n)$ does the job.
We next describe more formally the previous reasoning.

\begin{theorem}\label{reduction}
If almost exact recovery is solvable in $\sbm(n,p,\omega(1) Q /n)$, then exact recovery is solvable in $\sbm(n,p,\log(n) Q /n)$ if 
\begin{align}
I_+(p,Q)>1.
\end{align}
\end{theorem}
To see this, let $(X,G)\sim \sbm(n,p,\tau(n) Q/n)$, and $(G_1,G_2)$ be a graph splitting of $G$ with split-probability $\gamma=\log \log n/\log n$. Let $(X,\tilde{G}_2) \sim \sbm(n,p,(1-\gamma)\tau(n)Q/n)$ with $\tilde{G}_2$ independent of $G_1$ (note that the same $X$ appears twice). Let $\hat{X}=\hat{X}(G_1)$ be valued in $[k]^n$ such that $\pp\{A(X, \hat{X}) \ge 1-o(1) \}=1-o(1)$; note that such an $\hat{X}$ exists from the theorem's hypothesis. 
Since $A(X,\hat{X})=1-o(1)$ with high probability, $(G_2,\hat{X})$ are functions of $G$ and using a union bound, we have 
\begin{align}
&\pp\{ \hat{\Omega}_{\map}(G) \ne \Omega\} \\
&\le \pp \{  \hat{\Omega}_{\map}(G) \ne \Omega | A(X,\hat{X})=1-o(1) \} + o(1)\\
&\le \pp \{  \hat{\Omega}_{\map}(G_2,\hat{X}) \ne \Omega | A(X,\hat{X})=1-o(1) \} + o(1)\\
&\le n \pp \{ X_{1,\map}(G_2,\hat{X}_{\sim 1}) \neq X_1 | A(X,\hat{X})=1-o(1) \} + o(1).
\end{align}
We next replace $G_2$ by $\tilde{G}_2$. Note that $\tilde{G}_2$ has already the same marginal as $G_2$, the only issue is that $G_2$ is not independent from $G_1$ since the two graphs are disjoint, and since $\hat{X}$ is derived from $G_2$, some dependencies are carried along with $G_1$. However, $\tilde{G}_2$ and $G_2$ are `essentially independent' as stated in Lemma \ref{split_ind}, because the probability that $\tilde{G}_2$ samples an edge that is already present in $G_1$ is $O(\log^2 n/n^2)$, and the expected degrees in each graph is $O(\log n)$. This takes us to  
\begin{align}
&\pp\{ \hat{\Omega}_{\map}(G) \ne \Omega\} \le\\
& n \pp \{ X_{1,\map}(\tilde{G}_2,\hat{X}_{\sim 1}) \neq X_1 | A(X,\hat{X})=1-o(1) \} (1+o(1)) + o(1).
\end{align}
We can now replace $\hat{X}_{\sim 1}$ with $X_{\sim 1}$ to the expense that we may blow up this the probability by a factor $n^{o(1)}$ since $A(X,\hat{X})=1-o(1)$, again using the fact that expected degrees are logarithmic.
Thus we have 
\begin{align}
&\pp\{ \hat{\Omega}_{\map}(G) \ne \Omega\} \le \\
& n^{1+o(1)} \pp \{ X_{1,\map}(\tilde{G}_2,X_{\sim 1}) \neq X_1 | A(X,\hat{X})=1-o(1) \} + o(1) \label{tens1}
\end{align}
and the conditioning on $A(X,\hat{X})=1-o(1) $ can now be removed due to independence, so that 
\begin{align}
\pp\{ \hat{\Omega}_{\map}(G) \ne \Omega\} 
&\le n^{1+o(1)} \pp \{ X_{1,\map}(\tilde{G}_2,X_{\sim 1}) \neq X_1  \} + o(1).
\end{align}
The last step consists in closing the loop and replacing $\tilde{G}_2$ by $G$, since $1-\gamma=1-o(1)$, which uses the same type of argument as for the replacement of $G_2$ by $\tilde{G}_2$, with a blow up that is at most $n^{o(1)}$. As a result,
\begin{align}
\pp\{ \hat{\Omega}_{\map}(G) \ne \Omega\} 
&\le n^{1+o(1)} \pp \{ X_{1,\map}(G ,X_{\sim 1}) \neq X_1  \} + o(1), 
\end{align}
%A more delicate analysis behind step \eqref{tens1} can be made, by improving first the classification of $\hat{X}$ into one whose accuracy is not merely $o(1)$ but sufficiently small to prevent the blow up factor $n^{o(1)}$ to appear. Further, by working with a value of $\gamma$ that is itself dependent on the bit error probability $\pp \{ X_{1,\map}(G ,X_{\sim 1}) \neq X_1  \}$, the same $n^{o(1)}$ blow up can also be avoided in the last step above. In other words, one can obtain a tighter bound $\pp\{ \hat{S}_{map}(G) \ne S\} \le n  \pp \{ X_{1,\map}(G ,X_{\sim 1}) \neq X_1  \} + o(1)$. 
and if 
\begin{align}
\pp \{ X_{1,\map}(G ,X_{\sim 1}) \neq X_1  \} = n^{-1-\e}
\end{align}
for $\e>0$, then $\pp\{ \hat{\Omega}_{\map}(G) \ne \Omega\} $ is vanishing as stated in the theorem. 

Therefore, in view of Theorem \ref{reduction}, the achievability part of Theorem \ref{exact_thm} reduces to the following result.
\begin{theorem}\cite{colin1focs}
Almost exact recovery is efficiently solvable in $\sbm(n,p,\omega(1) Q /n)$. 
\end{theorem}
This follows from Theorem \ref{almost_exact2} discussed below, based on the Sphere-comparison algorithm discussed in Section \ref{partial}. %Note that to prove that almost exact recovery is solvable in this regime without worrying about efficiency, the Typicality Sampling Algorithm discussed in Section \ref{info} is already sufficient. 

In conclusion, in the regime of Theorem \ref{exact_thm}, exact recovery can be shown by using graph-splitting and combining almost exact recovery with degrees that grow sub-logarithmically and an additional clean-up phase. The behavior of the component-MAP error (i.e., the probability of misclassifying a single node when others have been revealed) pings down the behavior of the threshold: if this probability is $\omega(1/n)$, exact recovery is not possible, and if it is $o(1/n)$, exact recovery is possible. Decoding for the latter is then resolved by obtaining the exponent of the component-MAP error, which brings in the CH-divergence.\\

\noindent
{\bf Local to global amplification.}
The previous two sections give a lower bound and an upper bound on the probability that MAP fails at recovering the entire clusters, in terms of the probability that MAP fails at recovering a single vertex when others are revealed. Denoting by $P_{\mathrm{global}}$ and $P_{\mathrm{local}}$ these two probabilities of error, we essentially\footnote{The upper bound discussed in Section \ref{split} gives $n^{1+o(1)} P_{\mathrm{local}} + o(1)$, but the analysis can be tighten to yield a factor $n$ instead of $n^{1+o(1)}$.} have
\begin{align}
1-\frac{1}{nP_{\mathrm{local}}} +o(1) \le P_{\mathrm{global}} \le n P_{\mathrm{local}} + o(1).
\end{align}
This implies that $P_{\mathrm{global}}$ has a threshold phenomena as $P_{\mathrm{local}}$ varies:
\begin{align}
P_{\mathrm{global}} \to 	\begin{cases}
	0 & \text{ if } P_{\mathrm{local}} \ll 1/n, \\
	1 & \text{ if } P_{\mathrm{local}} \gg 1/n.
	\end{cases} 
\end{align}
Moreover, deriving this relies mainly on the regime of the model, rather than the specific structure of the SBM. In particular, it mainly relies on the exchangeability of the model (i.e., vertex labels have no relevance) and the fact that the vertex degrees do not grow rapidly. This suggests that this `local to global' phenomenon takes place in a more general class of models. The expression of the threshold for exact recovery in $\sbm(n,p,\log n Q /n)$ as a function of the parameters $p,Q$ is instead specific to the model, and relies on the CH-divergence in the case of the SBM, but the moderate/large deviation analysis of $P_{\mathrm{local}}$ for other models may reveal a different functional or $f$-divergence. 

The local to global approach also has an important implication at the computational level. The achievability proof described in the previous section directly gives an algorithm: use graph-splitting to produce two graphs; solve almost exact recovery on the first graph and locally improve the obtained clusters with the second graph. Since the second round is efficient by construction (it corresponds to $n$ parallel local computations), it is sufficient to solve almost exact recovery efficiently (in the regime of diverging degrees) to obtain for free an efficient algorithm for exact recovery down to the threshold. Thus this gives a computational reduction. In fact, the process can be iterated to further reduce almost exact recovery to a weaker recovery requirements, until a `bottle-neck' recovery problem is attained.

\section{Weak recovery and generalized KS threshold}\label{weak}

We recall the conjecture stated in \cite{decelle}:

\begin{conjecture}\label{c1}\cite{decelle,Mossel_SBM1}
Let $(X,G)$ be drawn from $\ssbm(n,k,a/n,b/n)$, i.e., the symmetric SBM with $k$ communities, probability $a/n$ inside the communities and $b/n$ across. Define $\snr=\frac{(a-b)^2}{k(a+(k-1)b)}$. Then, 
\begin{enumerate}
\item[(i)] For any $k \ge 2$, it is possible to solve weak recovery efficiently if and only if $\snr>1$ (the Kesten-Stigum (KS) threshold);
%\item If $k \in \{2,3,4\}$, and $\snr \leq 1$, it is impossible to detect communities information-theoretically,
\item[(ii)] If\footnote{The conjecture states that $k=5$ is necessary when imposing the constraint that $a>b$, but $k=4$ is enough in general.} $k \geq 4$, it is possible to solve weak recovery information-theoretically (i.e., not necessarily in polynomial time in $n$) for some $\snr$ strictly below 1.\footnote{\cite{decelle} made in fact a more precise conjecture, stating that below the KS threshold, there is a second transition for information-theoretic methods when $k\ge 4$, whereas there is a single threshold (for both efficient or non-efficient algorithms) when $k=3$. 
%\cite{decelle} provides also a more general conjecture for asymmetric SBMs, which was disproved in \cite{colin3cpam} due to mismatch in the definition of weak recovery.
}
\end{enumerate}
\end{conjecture}

%It was proved in \cite{massoulie-STOC,Mossel_SBM2} that the KS threshold can be achieved efficiently for $k=2$, with an alternative proof later given in \cite{bordenave}, and \cite{Mossel_SBM1} shows that it is impossible to detect below the KS threshold for $k=2$. Further, \cite{yash_sbm} extends the results for $k=2$ to the case where $a$ and $b$ diverge while maintaining the SNR finite. So detection is closed for $k=2$ in SSBM. 
It was also shown in \cite{bordenave} that for SBMs with communities that are balanced and for parameters that satisfy a certain asymmetry condition, i.e., the requirement that $\mu_k$ is a simple eigenvalue in Theorem 5 of \cite{bordenave}, the KS threshold can be achieved efficiently. The conditions of \cite{bordenave} do not cover Conjecture \ref{c1} for $k \geq 3$. In \cite{colin3,colin3cpam}, the two positive parts of the above conjecture are proved, with an extended result applying to the general SBM. We next discuss these various results.

Given parameters $p$ and $Q$ in the general model $\sbm(n,p,Q/n)$, let $P$ be the diagonal matrix such that $P_{i,i}=p_i$ for each $i \in [k]$. Also, let $\lambda_1,...,\lambda_h$ be the distinct eigenvalues of $PQ$ in order of nonincreasing magnitude. 
\begin{definition}
Define the signal-to-noise ratio of $\sbm(n,p,Q/n)$ by $$\snr=\lambda_2^2/\lambda_1.$$ 
\end{definition}
In the $k$ community symmetric case where vertices in the same community are connected with probability $a/n$ and vertices in different communities are connected with probability $b/n$, we have $\snr=(\frac{a-b}{k})^2/(\frac{a+(k-1)b}{k})=(a-b)^2/(k(a+(k-1)b))$, which matches the quantity in Conjecture 1 and in all previous sections concerned with two communities. 
%The following is our main result.  

\begin{theorem}\cite{bordenave}
Let $G\sim \sbm(n,p,Q/n)$ such that $p=(\frac{1}{k},\dots,\frac{1}{k})$ and such that $PQ$ has an eigenvector with corresponding eigenvalue in $(\sqrt{\lambda_1},\lambda_1)$ of single multiplicity. If $\snr>1$, then weak recovery is efficiently solvable.
\end{theorem}
\begin{theorem}\label{main}\cite{colin3,colin3cpam}
%Let $G$ be drawn under $\sbm(n,p,Q/n)$. 
Let $G \sim \sbm(n,p,Q/n)$ for $p,Q$ arbitrary. If $\snr>1$, then weak recovery is efficiently solvable.
\end{theorem}
Theorem \ref{main} implies the achievability part of Conjecture 1 part (i).\\

\noindent
The algorithm used in \cite{bordenave} is a spectral algorithm using the second eigenvector of the NB matrix discussed in Section \ref{linBP}. 
The algorithm used in \cite{colin3cpam} is an approximate acyclic belief propagation (ABP) algorithm, which corresponds to a power iteration method to extract the second eigenvector of the $r$-NB matrix.
% then there exist $r\in\mathbb{Z}^+$, $c>0$, and $m:\mathbb{Z}^+\rightarrow \mathbb{Z}^+$ such that $\mathrm{ABP}(G,m(n),r,c,(\lambda_1,...,\lambda_h))$ solves weak recovery and runs in $O(n\log n)$ time.

\begin{remark}
Also note that it is important to use the notion of weak recovery defined in Section \ref{all-defs}, where the agreement is normalized by the sizes of the communities. Without this normalization, using naively ``beating a random guess'' or the definition of \cite{decelle}, the conjecture that weak recovery is efficiently solvable if and only if $\snr >1$ is not true in general; a counter-example is given in \cite{colin3cpam} where the max-detection criteria of \cite{decelle} is not solvable when $\snr>1$.
\end{remark}

%For the symmetric SBM, the above reduces to part (1) of Theorem \ref{colin3_thm}, which proves the first part of Conjecture 1. 

%In particular, the conjecture from \cite{decelle} that max-detection is efficiently solvable if and only if $\snr >1$ for the general SBM is not true
%
%An example is obtained for example by taking an SSBM with two symmetric communities of intra-connectivity $3/n$, small enough extra-connectivity, and breaking one community into 4 identical sub-communities. Then max-detection is not solvable if $\snr>1$, but it is solvable for our notion of weak recovery. 

%If one uses  that definition, then we also conjecture that such a statement holds, i.e., detection is efficiently solvable if and only if $\lambda_2^2 > \lambda_1$ (besides for the case where $\lambda_1$ has multiplicity more than one, for which it is sufficient to have $\lambda_1>1$, and always focusing on the case of constant expected degrees).

%The full version of ABP is described in \cite{colin3}, but as for the symmetric case, the version ABP$^*$ described in previous section applies to the general setting, replacing $d$ with the largest eigenvalue of $PQ$. In \cite{bordenave}, a similar result to Theorem \ref{main} is obtained for the case where $p$ is uniform and $PQ$ has an eigenvalue $\lambda_k$ of multiplicity 1 such that $\lambda_k^2>\lambda_1$.

%Theorem \ref{main} provides the most general condition for solving efficiently detection in the SBM with linear size communities. 
We conjecture that Theorem \ref{main} is tight, i.e., if $\snr <1$, then efficient weak recovery is not solvable. However, establishing formally such a converse argument seems out of reach at the moment: as we shall see in the next section, except for a few possible cases with low values of $k$ (e.g., symmetric SBMs with $k=2,3$), it is possible to detect information-theoretically when $\snr<1$, and thus one cannot get a converse for efficient algorithms by considering all algorithms. 

%On the other hand, \cite{decelle} provides non-formal arguments based on statistical physics arguments that such a converse hold. It would be interesting to further connect the computational barriers occurring here with those from other problems such as planted clique \cite{clique}, or as in \cite{rigollet}. 

\section{Partial recovery}
The {\it Sphere-comparison} algorithm discussed in Section \ref{almost} gives the following result:
\begin{theorem}\cite{colin1}\label{almost_exact2}
For any $k\in \mathbb{Z}$, $p\in (0,1)^k$ and $Q$ with no two rows equal, there exist $\epsilon(c)=O(1/\log(c))$ such that for all sufficiently large $c$, {\it Sphere-comparison} detects communities in SBM$(n,p,c Q/n)$ with accuracy $1-e^{-\Omega(c)}$ and complexity $O_n(n^{1+\epsilon(c)})$. 
\end{theorem}
Tight expressions for partial recovery in the general SBM and at a finite SNR is open. 
We note the exception of a result obtained recently in \cite{coja2}, which requires the assortative regime. 
We also refer to \cite{harrison,prout}.

\chapter{The information-computation gap}\label{info}
In this section we discuss SBM regimes where weak recovery can be solved information-theoretically. As stated in Conjecture 1 and proved in Theorem \ref{main}, the  information-computation gap---defined as the gap between the KS and IT thresholds---takes place when the number of communities $k$ is larger than 4. We provide an information-theoretic (IT) bound for $\ssbm(n,k,a/n,b/n)$ that confirms this, showing further that the gap can grow fast with the number of communities. 

The information-theoretic bound described below is obtained by using a non-efficient algorithm that samples uniformly at random a clustering that is typical, i.e., that has the right proportions of edges inside and across the clusters. We describe below how this gives a tight expression in various regimes. Note that to capture the exact information-theoretic threshold in all regimes, one would have to rely on tighter estimates on the posterior distribution of the clusters given the graph. A possibility is to estimate the limit of the normalized mutual information between the clusters and the graph, i.e., $\frac{1}{n}I(X;G)$, as done in \cite{yash_sbm} for the regime of finite SNR with diverging degrees\footnote{Similar results  were also obtained recently in a more general context in \cite{lelarge_ks,lelarge_wigner}.}---see Section \ref{wigner}. Recent results also obtained the expression for the finite degree regime in the disassortative case \cite{coja2}. Another possibility is to estimate the limiting total variation or KL-divergence between the graph distribution in the SBM vs.\ Erd\H{o}s-R\'enyi model of matching expected degree. The limiting total variation is positive if and only if an hypothesis test can distinguish between the two models with a chance better than half. The easy implication of this is that if the total variation is vanishing, the weak recovery is not solvable (otherwise we would detect virtual clusters in the Erd\H{o}s-R\'enyi model). This used in \cite{banks} to obtain a lower-bound on the information-theoretic threshold, using a contiguity argument, see further details at the end of this section. 

\section{Crossing KS: typicality}

To obtain our information-theoretic upper-bound, we rely on the following sampling algorithm:\\ 

\noindent
{\bf Typicality Sampling Algorithm.} Given an $n$-vertex graph $G$ and $\delta>0$, the algorithm draws $\hat{X}_{\mathrm{typ}}(G)$ uniformly at random in  
\begin{align*}
&T_\delta(G)=\{   x \in \mathrm{Balanced}(n,k):\\
&  \sum_{i=1}^k |\{  G_{u,v} : (u,v) \in {[n] \choose 2} \text{ s.t. } x_u=i, x_v=i \}|  \geq \frac{an}{2k} (1-\delta),\\
& \sum_{i,j \in [k], i<j} |\{  G_{u,v} : (u,v) \in {[n] \choose 2} \text{ s.t. } x_u=i, x_v=j \}|  \leq \frac{bn(k-1)}{2k} (1+\delta) \},
\end{align*}
where the above assumes that $a\ge b$; flip the above two inequalities in the case $a<b$.\\

We present below a bound that we claim is tight at the extremal regimes of $a$ and $b$ (see discussions below). Note that for $b=0$, $\ssbm(n,k,a/n,0)$ is simply a patching of disjoint Erd\H{o}s-R\'enyi random graphs, and thus the information-theoretic threshold corresponds to the giant component threshold, i.e., $a>k$, achieved by separating the giants. This breaks down for $b$ positive, but we expect that the bound derived below remains tight in the scaling of small $b$. For $a=0$, the problem corresponds to planted coloring, which is already challenging \cite{kahale}. The bound obtained below gives that, in this case, weak recovery is information-theoretically solvable if $b >c k \log k + o_k(1)$, $c \in [1,2]$. This scaling is further shown to be tight in \cite{banks}, which also provides a simple upper-bound that scales as $k \log k$ for $a=0$.  Overall, the bound below shows that the KS threshold gives a much more restrictive regime than what is possible information-theoretically, as the latter reads $b > k(k-1)$ for $a=0$. 

%The behaviour of the IT bound is also characterized for $b$ close to 0. Similar though weaker results were also recently posted in \cite{banks}. Note also that the information-computation gap concerns the gap between the KS threshold and what is achieved information-theoretically, which is the gap between the information-theoretic and computational thresholds only under non-formal evidences \cite{decelle}.  Showing formally that no algorithms can succeed below the KS threshold would naturally require novel techniques and major progress on deep complexity theory questions. 

\begin{theorem}\label{main2}
Let $d:=\frac{a+(k-1)b}{k}$, assume $d>1$, and let $\tau=\tau_d$ be the unique solution in $(0,1)$ of $\tau e^{-\tau}=de^{-d}$, i.e., $\tau = \sum_{j=1}^{+ \infty} \frac{j^{j-1}}{j!} (de^{-d})^j$. The Typicality Sampling Algorithm solves weak recovery\footnote{Setting $\delta>0$ small enough gives the existence of $\e>0$ for weak recovery.} communities in $\ssbm(n,k,a/n,b/n)$ if
\begin{align}
&  \frac{a \log a + (k-1) b \log b}{k}  - \frac{a+(k-1)b}{k} \log \frac{a+(k-1)b}{k} \\
& >  \frac{1-\tau}{1-\tau k/(a+(k-1)b)} 2 \log (k) \\
& \quad \wedge ( 2\log (k)- 2\log(2)e^{-a/k}(1-(1-e^{-b/k})^{k-1})). \label{bound1}
\end{align}
%\begin{align}
%&\frac{1}{2 \log k} \left( \frac{a \log a + (k-1) b \log b}{k}  - \frac{a+(k-1)b}{k} \log \frac{a+(k-1)b}{k} \right) > \frac{1-\tau}{1-\tau k/(a+(k-1)b)}. \label{bound1}
%\end{align}
%or
%\[\frac{1}{k} \left(  \frac{a+b(k-1)}{2} \log \frac{k}{(a+(k-1)b)} +  \frac{a}{2} \log a + \frac{b(k-1)}{2} \log b \right)+e^{-a/k}(1-(1-e^{-b/k})^{k-1})\log(2)>\log(k)\]
%\begin{align}
%&\frac{1}{2 \log k} \left( \frac{a \log a + (k-1) b \log b}{k}  - d \log d \right) \left(1-\frac{\tau}{d} \right)^2 + \frac{\tau}{d} (1-\tau)  + \tau >1.
%\end{align}
%\begin{align}
%&\frac{1}{k \log (k)} \left(  \frac{a+b(k-1)}{2} \log \frac{k}{(a+(k-1)b)} +  \frac{a}{2} \log a + \frac{b(k-1)}{2} \log b \right)+ \frac{\tau}{d}\left(1-\frac{\tau}{2}\right) \\
%&+ \frac{1}{\log(k)} \left( \frac{a}{a+(k-1)b} \log\left(\frac{a+(k-1)b}{a}\right) +\frac{ (k-1)b}{a+(k-1)b} \log\left(\frac{a+(k-1)b}{b}\right) \right) \\
%& \phantom{- \frac{1}{\log(k)}} \cdot \left( \frac{\tau^2}{2d} + (d -\tau)  e^{-(d-\tau)}\right)\\
%&>1.
%\end{align}
\end{theorem}
This bound strictly improves on the KS threshold for $k\ge 4$. See \cite{colin3cpam} for a numerical example.
\begin{corollary}
Conjecture 1 part (ii) holds. 
\end{corollary}
%This bound strictly improves on the KS threshold for $k\ge 4$, settling part (ii) of Conjecture 1. 
%see Figure \ref{cross}.
%\Cnote{We should probably get a new figure for $k=4$.}
%\begin{figure}[H]
%\centering
%  \includegraphics[width=0.5\linewidth]{cross.png}
%\caption{The horizontal axis has $b$ varying, while $a=1/10$ and $k=5$. The red (dashed) curve is the KS threshold, i.e., the values of $b$ for which $(a-b)^2/(k(a+(k-1)b))=1$, and the blue (plain) curve is our IT threshold. When the curves are positive, detection is solvable. In particular, detection is solvable below the KS threshold. }
%\label{cross}
%\end{figure}

Note that \eqref{bound1} simplifies to
\begin{align}
&\frac{1}{2 \log k} \left( \frac{a \log a + (k-1) b \log b}{k}  - d\log d \right) > \frac{1-\tau}{1-\tau /d} =:f(\tau,d),
\end{align}
and since $f(\tau,d)<1$ when $d>1$ (which is needed for the presence of the giant), weak recovery is already solvable in $\sbm(n,k,a,b)$ if 
\begin{align}
\frac{1}{2 \log k} \left( \frac{a \log a + (k-1) b \log b}{k}  - d \log d \right) >1. \label{weakb}
\end{align}
The above bound corresponds to the regime where there is no bad clustering that is typical with high probability. The analog of this bound in the unbalanced case already provides examples to crossing KS for two communities, such as for $p=(1/10,9/10)$ and $Q=(0,81;81,72)$. 
However, the above bound is not tight in the extreme regime of $b=0$, since it reads $a>2k$ as opposed to $a>k$, and it only crosses the KS threshold at $k=5$. Before explaining how to obtain tight interpolations, we provide further insight on the bound of Theorem \ref{main2}.
%An intermediate bound is also obtained by dropping the term \eqref{drop} in the proof of the theorem, corresponding to ignoring planted trees in the giant (see Section \ref{proof_tech2}), which gives 
%\begin{align}
%\frac{1}{2 \log k} \left( \frac{a \log a + (k-1) b \log b}{k}  - d \log d \right) >1 - \frac{\tau}{d}\left(1-\frac{\tau}{2}\right).
%\end{align}
%\end{remark}

Defining $a_k(b)$ as the unique solution of 
\begin{align}
&\frac{1}{2 \log k} \left( \frac{a \log a + (k-1) b \log b}{k}  - d\log d \right) \\
&= \min\left(f(\tau,d),1-\frac{e^{-a/k}(1-(1-e^{-b/k})^{k-1})\log(2)}{\log(k)}\right)
\end{align}
and simplifying the bound in Theorem \ref{main2} gives the following. 
\begin{corollary}\label{corol_ext}
Weak recovery is solvable  
\begin{align}
&\text{in}\quad \sbm(n,k,0,b) \quad \text{if} \quad b >   \frac{2k \log k}{(k-1) \log \frac{k}{k-1}} f(\tau,b(k-1)/k), \label{a0} \\ 
%&\sbm(n,k,0,b) \quad \text{if} \quad b >   k \log k \cdot 2f(\tau,1/b) + o_k(1), \label{a0} \\
%&\sbm(n,k,\e,b) \quad \text{if} \quad b > 2+ \Delta_2(\e), \quad \text{where} \quad  \lim_{\e \to 0} \Delta_2(\e)=0 \label{a00}, \\
&\text{in}\quad \sbm(n,k,a,b) \quad \text{if} \quad a > a_k(b), \quad \text{where } a_k(0)=k. \label{b0}
\end{align}
\end{corollary}
\begin{remark}
Note that \eqref{b0} approaches the optimal bound given by the presence of the giant at $b=0$, and we further conjecture that $a_k(b)$ gives the correct first order approximation of the information-theoretic bound for small $b$. 
\end{remark}
\begin{remark}
%Note that the $k \log k$ scaling in \eqref{a0} improves significantly on the KS threshold given by $b>k (k-1)$ at $a=0$. 
Note that the $k$-colorability threshold for Erd\H{o}s-R\'enyi graphs grows as $2k\log k$ \cite{naor_col}. This may be used to obtain an information-theoretic bound, which would however be looser than the one obtained above. 
\end{remark}
%We can cross the KS threshold for some positive values of $a$ at $k=5$. 
%\begin{lemma}\label{cross}
%Detection is solvable in $\mathrm{SBM}(n,k,a,b)$ if 
%\begin{align}
%&  (a+(k-1)b) \log \frac{k}{a+(k-1)b}+ a\log a + (k-1)b \log b  > 2k\log k.
%\end{align}
%In particular, detection is solvable in $\mathrm{SSBM}(n,k,0,b/n)$ if 
%\begin{align}
%\frac{b}{k(k-1)} >  \frac{2 \log k}{(k-1)^2 \log \frac{k}{k-1}}.
%\end{align}
%\end{lemma}
It is possible to see that this gives also the correct scaling in $k$ for $a=0$, i.e., that for $b< (1-\e)k \log (k) + o_k(1)$, $\e>0$, weak recovery is information-theoretically impossible. To see this, consider $v\in G$, $b=(1-\epsilon)k \log(k)$, and assume that we know the communities of all vertices more than $r=\log(\log(n))$ edges away from $v$. For each vertex $r$ edges away from $v$, there will be approximately $k^{\epsilon}$ communities that it has no neighbors in. Then vertices $r-1$ edges away from $v$ have approximately $k^{\epsilon}\log(k)$ neighbors that are potentially in each community, with approximately $\log(k)$ fewer neighbors suspected of being in its community than in the average other community.  At that point, the noise has mostly drowned out the signal and our confidence that we know anything about the vertices' communities continues to degrade with each successive step towards $v$.

A different approach is developed in \cite{banks} to prove that the scaling in $k$ is in fact optimal, obtaining both upper and lower bounds on the information-theoretic threshold that match in the regime of large $k$ when $(a-b)/d = O(1)$. In terms of the expected degree, the threshold reads as follows. 
\begin{theorem}\cite{banks,banks2}
When $(a-b)/d = O(1)$, the critical value of $d$ satisfies $d = \Theta \left( \frac{d^2 k \log k}{(a-b)^2} \right)$, i.e., the critical SNR satisfies $\snr = \Theta (\log(k)/k)$.
\end{theorem}
 The upper-bound in \cite{banks} corresponds essentially to \eqref{weakb}, the regime in which the first moment bound is vanishing. The lower-bound is based on a contiguity argument and second moment estimates from \cite{naor_col}. The idea is to compare the distribution of graphs drawn from the SBM, i.e., 
\begin{align}
\mu_{\sbm} (g) := \sum_{x \in [k]^n} \pp\{ G=g | X=x \} \pp\{X=x \}
\end{align}
with the distribution of graphs drawn from the Erd\H{o}s-R\'enyi model with matching expected degree, call it $\mu_{\mathrm{ER}}$. If one can show that 
\begin{align}
\| \mu_{\sbm} -  \mu_{\mathrm{ER}} \|_1 \to 0, \label{tot}
\end{align}
then upon observing a graph drawn from either of the two models, say with probability half for each, it is impossible to decide from which ensemble the graph is drawn with high probability. Thus it is not possible to solve weak recovery (otherwise one would detect clusters in the Erd\H{o}s-R\'enyi model). A sufficient condition to imply \eqref{tot} is to show that $\mu_{\sbm} \trianglelefteq \mu_{\mathrm{ER}}$, i.e., for any sequence of event $E_n$ such that $\mu_{\mathrm{ER}}(E_n) \to 0$, it must be that $\mu_{\sbm} \to 0$. In particular, $\mu_{\sbm}$ and $\mu_{\mathrm{ER}}$ are called contiguous if $\mu_{\sbm} \trianglelefteq \mu_{\mathrm{ER}}$ and $\mu_{\mathrm{ER}} \trianglelefteq \mu_{\sbm}$, but only the first of these conditions is needed here. Further, this is implied from Cauchy-Schwarz if the ratio function $$ \rho (G):=\mu_{\sbm}(G)/\mu_{\mathrm{ER}}(G)$$ has a bounded second moment, i.e., $$\E_{G \sim \mathrm{ER}}\rho^2(G) =O(1),$$ which is shown in \cite{banks}; see \cite{moore_review} for more details.

\section{Nature of the gap}
The nature of such gap phenomena can be seen from different perspectives. One interpretation comes from the behavior of belief propagation. See also \cite{moore_review} for further discussions. 

Above the Kesten-Stigum threshold, the uniform fixed point is unstable and BP does not get attracted to it and reaches a non-trivial solution on most initializations. In particular, the ABP algorithm discussed in Section \ref{linBP}, which starts with a random initialization with order $\sqrt{n}$ vertices towards the true partition (due to the Central Limit Theorem), is enough to make linearized BP reach a non-trivial fixed point.  
Below the information-theoretic threshold, the non-trivial fixed points are no longer present, and BP settles in a solution that represents a noisy clustering, i.e., one that would also take place in the Er\H{o}s-R\'enyi model due to the noise fluctuations in the model. In the gap region, non-trivial fixed points are still present,  but the trivial fixed points are locally stable and attract most initializations. One could try multiple initializations until a non-trivial fixed point is reached, using for example the graph-splitting technique discussed in Section \ref{weak-it} to test such solutions. However, it is believed that an exponential number of initializations is needed to reach a good solution. 

This connects to the ``energy landscape'' of the possible clusterings: in this gap region, the non-trivial fixed points have a very small basin of attraction, and they can only attract an exponentially small fraction of initializations. To connect to the results from Section \ref{exact}, the success of the two-round procedure can also be related to the energy landscape, i.e., the objective function of MAP in this case. Above the CH threshold, an almost exact solution having $n-o(n)$ correctly labeled vertices can be converted to an exact solution by the degree-profiling hypothesis test. This is essentially saying that BP at depth 1, i.e., computing the likelihood of a vertex based on its direct neighbors, allows us to reach the global maximum of the likelihood function with such a strong initialization. In other words, the BP view, or more precisely understanding the question of how accurate our initial beliefs need to be in order to amplify these to non-trivial levels based on neighbors at a given depth, is related to the landscape of the objective function. 

The gap phenomenon also admits a local manifestation in the context of ABP, having to do with the approximation discussed in Section \ref{linBP}, where the non-linear terms behave differently from $k=3$ to $k=4$ due to the loss of a diminishing return property. Better understanding such gap phenomena is an active research area.

\section{Proof technique for crossing KS} 
\begin{figure}[h]
\centering
  \includegraphics[width=.8\linewidth]{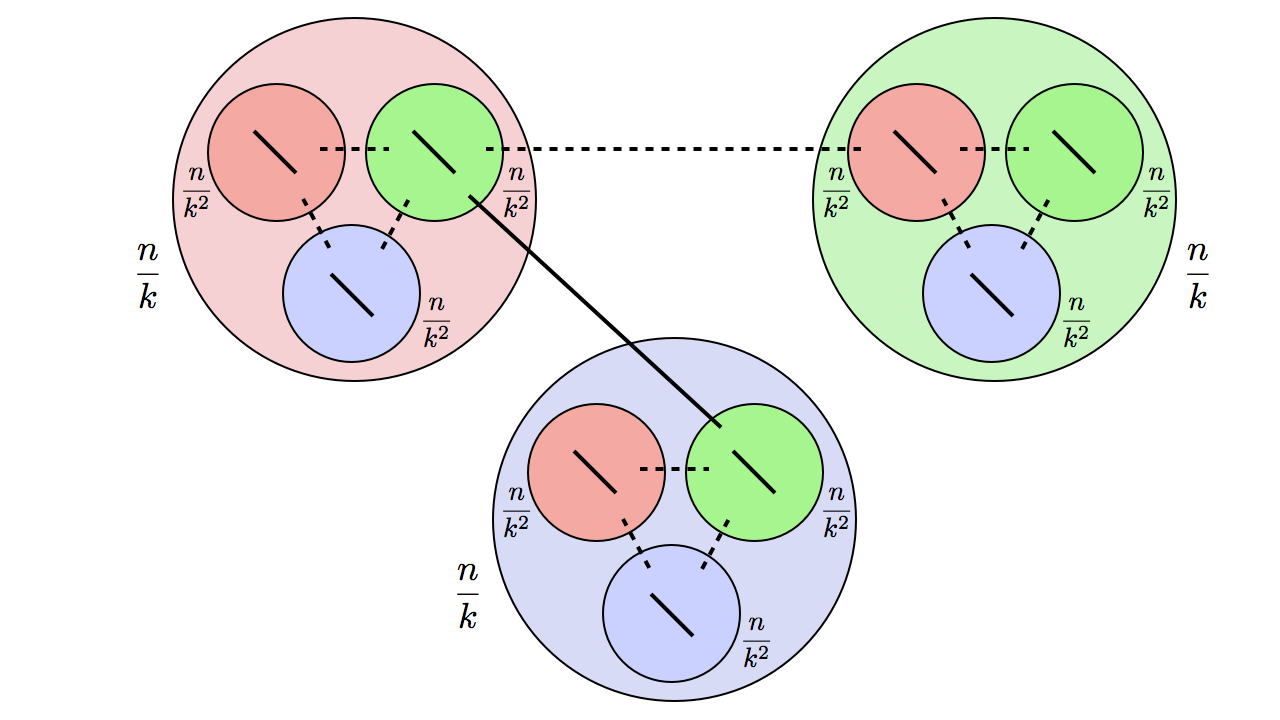}
  \caption{A bad clustering roughly splits each community equally among the $k$ communities. Each pair of nodes connects with probability $a/n$ among vertices of same communities (i.e., same color groups, plain line connections), and $b/n$ across communities (i.e., different color groups, dashed line connections). Only some connections are displayed in the figure to ease the visualization.}
  \label{bad_cluster}
\end{figure}

We explain in this section how to obtain the bound in Theorem \ref{main2}.
A first problem is to estimate the likelihood that a bad clustering, i.e., one that has an overlap close to $1/k$ with the true clustering, belongs to the typical set. As clusters sampled from the TS algorithm are balanced, a bad clustering must split each cluster roughly into $k$ balanced subgroups that belong to each community, see Figure \ref{bad_cluster}. Thus it is unlikely to keep the right proportions of edges inside and across the clusters, but depending on the exponent of this rare event, and since there are exponentially many bad clusterings, there may exist one bad clustering that looks typical.

As illustrated in Figure \ref{bad_cluster}, the number of edges that are contained in the clusters of a bad clustering is roughly distributed as the sum of two Binomial random variables, 
\begin{align*}
E_\mathrm{in} \stackrel{\cdot}{\sim} \bin \left(\frac{n^2}{2k^2},\frac{a}{n}\right) + \bin \left(\frac{(k-1)n^2}{2k^2},\frac{b}{n}\right),
\end{align*}
where we use $\stackrel{\cdot}{\sim}$ to emphasize that this is an approximation that ignores the fact that the clustering is not perfect bad and perfectly balanced. Note that the expectation of the above distribution is $\frac{n}{2k} \frac{a+(k-1)b}{k}$. In contrast, the true clustering would have a distribution given by $\bin(\frac{n^2}{2k},\frac{a}{n})$, which would give an expectation of $\frac{an}{2k}$.
In turn, the number of edges that are crossing the clusters of a bad clustering is roughly distributed as
\begin{align*}
E_\mathrm{out} \stackrel{\cdot}{\sim}  \bin \left(\frac{n^2(k-1)}{2k^2} ,\frac{a}{n}\right) + \bin \left(\frac{n^2(k-1)^2}{2k^2} ,\frac{b}{n}\right),
\end{align*}
which has an expectation of $\frac{n(k-1)}{2k} \frac{a+(k-1)b}{k}$.
In contrast, the true clustering would have the above replaced by $\bin(\frac{n^2(k-1)}{2k},\frac{b}{n})$, and an expectation of $\frac{bn(k-1)}{2k}$. 

Thus, we need to estimate the rare event that the Binomial sum deviates from its expectation. While there is a large list of bounds on Binomial tail events, the number of trials here is quadratic in $n$ and the success bias decays linearly in $n$, which require particular care to ensure tight bounds. We derive these in \cite{colin3}, obtaining that $\mathbb{P} \{ x_{\mathrm{bad}} \in T_\delta(G) | x_{\mathrm{bad}} \in B_\epsilon  \}$ behaves as
\begin{align*}
&\exp \left(- \frac{n}{k} A \right) 
\end{align*}
when $\e,\delta$ are arbitrarily small, where  $A := \frac{a+b(k-1)}{2} \log \frac{k}{a+(k-1)b} +  \frac{a}{2} \log a + \frac{b(k-1)}{2} \log b$. One can then use the fact that $|T_\delta(G)| \ge 1$ with high probability, since the planted clustering is typical with high probability, and using a union bound and the fact that there are at most $k^n$ bad clusterings:

\begin{align}
P\{ \hat{X}(G) \in B_\epsilon \} &= E_G \frac{|T_\delta(G) \cap B_\epsilon|}{|T_\delta(G)|} \label{tighter} \\
&\leq E_G |T_\delta(G) \cap B_\epsilon| + o(1) \label{loose} \\
&\leq k^n  \cdot \, \mathbb{P} \{ x_{\mathrm{bad}} \in T_\delta(G) | x_{\mathrm{bad}} \in B_\epsilon  \} +o(1) . \notag
\end{align}

Checking when the above upper-bound vanishes already gives a regime that crosses the KS threshold when $k\geq 5$ for symmetric communities, when $k \ge 2$ for asymmetric communities (deriving the version of the bound for asymmetric communities), and scales properly in $k$ when $a=0$. However, it does not interpolate the correct behavior of the information-theoretic bound in the extreme regime of $b=0$ and does not cross at $k=4$. In fact, for $b=0$, the union bound requires $a>2k$ to imply no bad typical clustering with high probability, whereas as soon as $a>k$, an algorithm that simply separates the two giants in $\sbm(n,k,a,0)$ and assigns communities uniformly at random for the other vertices solves weak recovery. Thus when $a \in (k,2k]$, the union bound is loose. To remediate this, we next take into account the topology of the SBM graph to tighten our bound on $|T_\delta(G)|$.

Since the algorithm samples a typical clustering, we only need the number of bad and typical clusterings to be small compared to the total number of typical clusterings, in expectation. Namely, we can get a tighter bound on the probability of error of the TS algorithm by obtaining a tighter bound on the typical set size than simply 1, i.e., estimating \eqref{tighter} without relying on the loose bound from \eqref{loose}. We proceed here with three levels of refinements to bound the typical set size. In each level, we construct a random labelling of the vertices that remains typical, and then use entropic estimates to count the number of such typical labellings. 

First we exploit the large fraction of nodes that are in tree-like components outside of the giant, and the labels are distributed on such trees as in the broadcasting on trees problem \ref{bot}. Specifically, for a uniformly drawn root node $X$, each edge in the tree acts as a $k$-ary symmetric channel.
%\begin{align*}
%&Z \sim \\
%&\nu:=\left(\frac{a}{a+(k-1)b}, \frac{b}{a+(k-1)b}, \dots, \frac{b}{a+(k-1)b} \right)
%\end{align*}
Thus, labelling the nodes in the trees according to the above distribution and freezing the giant to the correct labels leads to a typical clustering with high probability. The resulting bound matches the giant component bound at $b=0$, but is unlikely to scale properly for small $b$. To improve on this, we next take into account the vertices in the giant that belong to planted trees, and follow the same program as above, except that the root node (in the giant) is now frozen to the correct label rather than being uniformly drawn. This gives a bound that we claim is tight at the first order approximation when $b$ is small.  Finally, we also take into account vertices that are not saturated, i.e., whose neighbors do not cover all communities and who can thus be swapped without affecting typicality. The final bound allows to cross at $k=4$.

%The combined estimated on the typical set are in the following bound, which leads to the bound in Theorem \ref{}. 

%\vspace{-.8cm}
%\begin{figure}[h]
%\centering
%  \includegraphics[width=1\linewidth]{topo.png}
%  \caption{Illustration of the topology of $\sbm(n,k,a,b)$ for $k=2$. A giant component covering the two communities takes place when $d=\frac{a+(k-1)b}{k}>1$; a linear fraction of vertices belong to isolated trees (including isolate vertices). To estimate the size of the typical set: sample a bit uniformly at random in each isolated tree (green vertices) and propagate the bit according to the symmetric channel with flip probability $b/(a+(k-1)b)$ (plain edges do not flip whereas dashed edges flip). }
%  \label{topo}
%\end{figure}
%\vspace{-.5cm}

\begin{figure}[h]
\centering
  \includegraphics[width=0.9\linewidth]{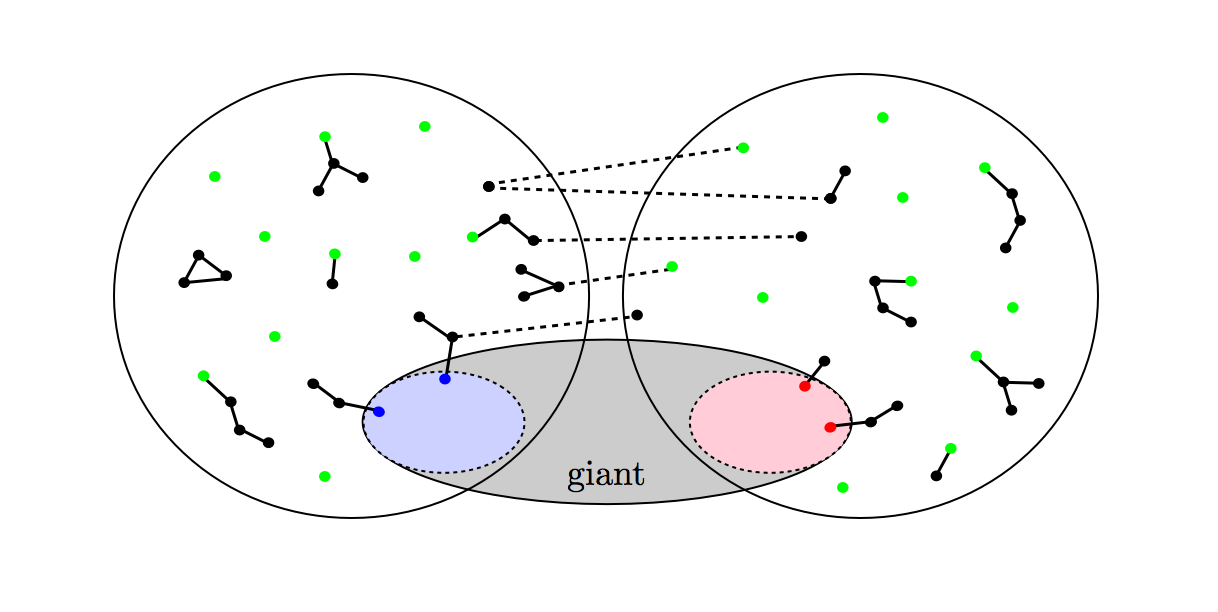}
  \caption{Illustration of the topology of $\sbm(n,k,a,b)$ for $k=2$. A giant component covering the two communities takes place when $d=\frac{a+(k-1)b}{k}>1$; a linear fraction of vertices belong to isolated trees (including isolated vertices), and a linear fraction of vertices in the giant are on planted trees. The following is used to estimate the size of the typical set in \cite{colin3cpam}. For isolated trees, sample a bit uniformly at random for a vertex (green vertices) and propagate the bit according to the symmetric channel with flip probability $b/(a+(k-1)b)$ (plain edges do not flip whereas dashed edges flip). For planted trees, do the same but freeze the root bit to its true value.}
  \label{topofig}
\end{figure}

\chapter{Other block models}\label{others}

%Finally, the extension to models discussed in Section \ref{extensions} can also be understood in the lens of weak recovery. The case of edge labels or hyperedges need a separate treatment, since the reductions described in Section \ref{extensions} are specific to exact recovery. The converse for weak recovery in the labelled SBM is covered in \cite{massoulie_label2}. The bipartite case can instead be treated as a special case of the threshold $\lambda_2^2/\lambda_1>1$ when the matrix $Q$ has 0 entries.  

%Concerning the separation of specific communities, the results from \cite{colin3} imply that communities $i$ and $j$ can be separated if there exists an eigenvector $w$ of $PQ$ such that $w_i\ne w_j$ and the eigenvalue of $PQ$ corresponding to $w$ has a magnitude of at least $\sqrt{\lambda_1}$. 

%\begin{theorem}
%Let $k \in \mZ_+$, $p=\{1/k\}^k$, $a,b>0$, $Q$ be the $k\times k$ matrix with $Q_{i,i}=a$ for all $i$ and $Q_{i,j}=b$ for all $i\ne j$ and $G$ drawn under $\sbm(n,p,Q/n)$. If \eqref{ks1} holds then there exist $r\in\mathbb{Z}^+$ and $m:\mathbb{Z}^+\rightarrow \mathbb{Z}^+$ such that $\mathrm{ABP}(G,m(n),r,((a-b)/k,(a+(k-1)b)/k))$ solves detection (according to our definition and a fortiori to Decelle's definition) and runs in $O(n\log n )$ time.
%\end{theorem}

%\section{Model variants}\label{other_models}
There are various extensions of the basic SBM discussed in previous section.
The variations increase yearly, and we mention here a few basic variants:
\begin{itemize}
\item {\bf Labelled SBMs:} allowing for edges to carry a label, which can model intensities of similarity functions between vertices; see for example \cite{massoulie_label2,jiaming,jog,prout2}; see also \cite{abbe_jmlr} for a reduction from labelled edges to unlabelled edges for certain recovery requirements;
\item {\bf Degree-corrected SBMs:} allowing for a degree parameter for each vertex that scales the edge probabilities in order to make expected degrees match the observed degrees; see for example \cite{newman2} and \cite{massoulie_degree1,massoulie_degree2} for sharp results on such models;
\item {\bf Overlapping SBMs:} allowing for the communities to overlap, such as in the mixed-membership SBM \cite{mixed-core}, where each vertex has a profile of community memberships or a continuous label---see also \cite{fortunato,newman2015generalized,ball,Peixoto2015,palla,prem,colin1focs}; also \cite{abbe_jmlr} for reductions to non-overlapping community models in some cases and recent results that sharp for MMSBM in \cite{david1};
\item {\bf Metric and geometric SBMs:} allowing for labels at the vertices that leave in metric spaces, e.g., a grid, a Euclidean space or a sphere, and where connectivity depends on the distance between the vertices' labels as further discussed below; see  \cite{our_grid,abishek2,arya,powering1};  
\end{itemize}

\begin{definition}[Geometric block models] 
We define here two geometric block models with two communities, the sphere-GBM and the mixture-GBM.
For each model, $X^n=(X_1,\dots,X_n)$ has i.i.d.\ Bernoulli$(1/2)$ components, which represents the {\it abstract} community labels for each vertex. We next add a {\it geometric} label for each vertex, and draw the graph depending on both the abstract and geometric labels:
\begin{itemize}
\item In the sphere-GBM$(n,d,\tau,a,b)$, $U^n =(U_1,\dots,U_n)$ has i.i.d.\ components drawn uniformly at random on a sphere of dimension $d$; the graph $G=([n],E)$ is drawn with edges independent conditionally on $X^n,U^n$, such that for $1 \le i < j \le n$, 
\begin{align}
&\pp\{E_{ij}=1 | X^n=x^n, U^n=u^n \}  \notag \\
&= 
\begin{cases}
a & \text{ if } \| u_i - u_j\| \le \tau \text{ and } x_i=x_j \\
b & \text{ if } \| u_i - u_j\| \le \tau \text{ and } x_i\ne x_j\\
0 & \text{ if } \| u_i - u_j\| > \tau
\end{cases}
\end{align} 
and $\pp\{E_{ij}=0 | X^n=x^n, U^n=u^n \}=1-\pp\{E_{ij}=1 | X^n=x^n, U^n=u^n \}$. 
One can have a variant with a special symbol $\star$ that indicates if $ \| u_i - u_j\| > \tau$.
\item In the mixture-GBM$(n,d,s,\tau)$, $U^n =(U_1,\dots,U_n)$ has independent components conditionally on $X^n$ with $U_i$ drawn from $\mathcal{N}(0^{d},I_{d})$ if $X_i=0$ and from $\mathcal{N}((s,0^{d-1}),I_{d})$ if $X_i=1$ (two isotropic Gaussians in dimension $d$ at distance $s$); the graph $G=([n],E)$ is drawn with edges independent conditionally on $U^n$, such that for $1 \le i < j \le n$, 
\begin{align}
\pp\{E_{ij}=1 | U^n=u^n \} = 
\begin{cases}
1 & \text{ if } \| u_i - u_j\| \le \tau \\
0 & \text{ if } \| u_i - u_j\| > \tau
\end{cases}
\end{align} 
and $\pp\{E_{ij}=0 | U^n=u^n \}=1-\pp\{E_{ij}=1 | U^n=u^n \}$.
\end{itemize}
\end{definition}
In dimension $d$, we want $\tau$ to be at least of order $n^{-1/d}$ with a large enough constant in order for the graphs to have giant components, and $(\log(n)/n)^{1/d}$ with a large enough constant for connectivity. 

Previous models have a much larger number of short loops than the SBM does, which captures a feature of various real world graphs having transitive attributes (``friends of friends are more likely to be friends'').  
On the flip side, these models do not have `abstract edges' as in the SBM, which can also occur frequently in applications given the ``small-world phenomenon''. This says that real graphs often have relatively low diameter (about 6 in the case of Milgram's experiment), which does not take place in purely geometric block models. Therefore, a natural candidate is to superpose an SBM and a GBM to form an hybrid block model (HBM), see for example \cite{powering1}.
The many loops of the GBM can be challenging for some of the algorithms discussed in this monograph, in particular for basic spectral methods and even belief propagation; we discussed in Section \ref{robust} how graph powering provides a more robust alternative in such cases. 

Another variant that circumvents the discussions about non-edges is to consider a {\bf censored block model} (CBM), defined as follows (see \cite{massoulie_label2,abbs}). 
\begin{definition}[Binary symmetric CBM]
Let $G=([n],E)$ be a graph and $\e \in [0,1]$.
Let $X^n=(X_1,\dots,X_n)$ with i.i.d.\ Bernoulli$(1/2)$ components.
Let $Y$ be a random vector of dimension ${n \choose 2}$ taking values in $\{0,1,\star\}$ such that 
\begin{align}
&\pp\{Y_{ij}=1|X_i=X_j,E_{ij}=1\}=\pp\{Y_{ij}=0|X_i\neq X_j,E_{ij}=1\}=\e\\
&\pp\{Y_{ij}=\star|E_{ij}=0\}=1.
\end{align}
\end{definition}
The case of an Erd\H{o}s-R\'enyi graph is discussed in \cite{abbs,abbs-isit,new-vu,florent_CBM,Chen_Huang_Guibas_Graphics,Chen_Goldsmith_ISIT2014} and the case of a grid in \cite{our_grid}. For a random geometric graph on a sphere, it is closely related to the above sphere-GBM.
Inserting the $\star$ symbol simplifies a few aspects compared to SBMs, such as Lemma \ref{indep-l}, which is needed in the weak recovery converse of the SBM. In this sense, the CBM is a more convenient model than the SBM from a mathematical viewpoint, while behaving similarly to the SBM (when $G$ is an Erd\H{o}s-R\'enyi graph of degree $(a+b)/2$ and $\e=b/(a+b)$ for the two community symmetric case).
The CBM can also be viewed as a synchronization model over the binary field, and more general synchronization models have been studied in \cite{wein1,wein2}, with a complete description both at the fundamental and algorithmic level (generalizing in particular the results from Section \ref{wigner}).

Note all previously mentioned models are different forms of latent variable models. 
Focusing on the graphical nature, on can also consider the more general {\bf inhomogenous random graphs} \cite{bollo_inhomo}, which attach to each vertex a label in a set that is not necessarily finite, and where edges are drawn independently from a given kernel conditionally on these labels. This gives in fact a way to model mixed-membership, and is also related to {\bf graphons}, which corresponds to the case where each vertex has a continuous label. 

It may be worth saying a few words about the theory of graphons and its implications for us. Lov\'asz and co-authors introduced graphons \cite{graphon2,graphon3,graphon} in the study of large graphs (also related to Szemer\'edi's Regularity Lemma \cite{regularity}), showing that\footnote{Initially in dense regimes and more recently for sparse regimes \cite{graphon_sparse}.} a convergent sequence of graphs admits a limit object, the graphon, that preserves many local and global properties of the sequence. Graphons can be represented by a measurable function $w : [0,1]^2 \to [0,1]$, which can be viewed as a continuous extension of the connectivity matrix $W$ used throughout this paper. Most relevant to us is that any network model that is invariant under node labelings, such as most models of practical interest, can be described by an edge distribution that is {\it conditionally independent} on hidden node labels, via such a measurable map $w$. This gives a de Finetti's theorem for label-invariant models \cite{hoover,aldous-array,diaconis-graphon}, but does not require the topological theory behind it. Thus the theory of graphons may give a broader meaning to the study of block models, which are precisely building blocks to graphons, but for the sole purpose of studying exchangeable network models, inhomogeneous random graphs give enough degrees of freedom.

%The theory of graphons goes in reality much beyond the actual random graph model that is discussed here. It is a topological theory that allows to define limits of graph sequences, and the random graph model comes into play in the establishment of the limit, but the theory is not necessary to discuss models for machine learning or network problems. However, 

Further, many problems in machine learning and networks are also concerned with interactions of items that go beyond the pairwise setting. For example, citation or metabolic networks rely on interactions among $k$-tuples of vertices. These can be captured by extending previous models to hypergraphs.\footnote{Recent results for community detection in hypergraphs were obtained in \cite{afonso-goemans}.}

\chapter{Concluding remarks and open problems}\label{open}
One may cast the SBM and the previously discussed variants into a comprehensive class of conditional random fields \cite{lafferty2001conditional} or graphical channels \cite{abbetoc}, where edge labels depend on vertex labels. 
%\footnote{A recent paper \cite{rigollet2} has also considered an Ising model with block structure, studying exact recovery and SDPs in this context.} 
\begin{definition}
Let $V=[n]$ and $G=(V,E(G))$ be a hypergraph with $N=|E(G)|$.
Let $\X$ and $\Y$ be two finite sets called, respectively, the input and output alphabets, and $Q(\cdot | \cdot )$ be a channel from $\X^k$ to $\Y$ called the kernel.  
To each vertex in $V$, assign a vertex-variable in $\X$, and to each edge in $E(G)$, assign an edge-variable in $\Y$. 
%, $k=2$, $\X=\{0,1\}$ is a community variable and $\Y=\{0,1\}$ is an edge variable. 
Let $y_I$ denote the edge-variable attached to edge $I$, and $x[I]$ denote the $k$ node-variables adjacent to $I$.
We define a graphical channel with graph $G$ and kernel $Q$ as the channel $P(\cdot |\cdot)$ given by 
\begin{figure}[h]
\begin{center}
\includegraphics[scale=.25]{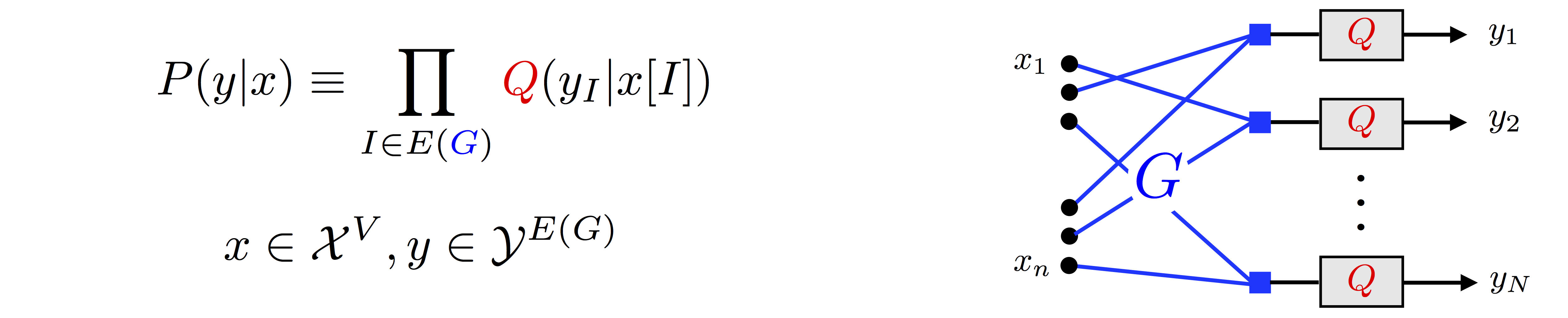}
%\caption{Graphical channel with graph $G$ and kernel $Q$}
\label{graphical_fig2}
\end{center}
\vspace{-0.5cm}
\end{figure}
\end{definition}

%The above departs significantly from a traditionally encoded channel when considering low order edges (e.g., $k=2,3$) and $G$ uniform or complete (the closest would be a special LDGM code \cite{kumar}). As discussed in \cite{colin1focs}, exact recovery in the SBM is verbatim a decoding problem on such a channel with an LDGM code of right-degree 2. 

%Community detection has a strong connection with information theory since $\X$ is typically discrete (as the goal is to obtain `clusters' on the data), which is not common for other applications in machine learning where the real-valued nature of the channel is important.\footnote{Compressed sensing or topic modelling rely instead heavily on real-valued channels.} Graphical channels allow also to capture many extensions of the SBM, such as non-overlapping communities, edge-labeled or non-pairwise interactions. This can be further extended to problems such as topic modelling or ranking, with new notions of recovery.  
%The above model captures a large collection of problems in machine learning, mainly in the unsupervised setting. 
Quantities that are key to understanding how much information can be carried in such graphical channels are: 
\begin{center}how ``rich'' is the observation graph $G$ and \\ how ``noisy'' is the connectivity kernel $Q$. %These are not the usual notions of rates and capacity in information theory, but they are the relevant ones here. These also make the problems novel and interesting. 
\end{center}
This survey quantifies the tradeoffs between these two quantities in the SBM (which corresponds to a discrete $\X$ and a specific graph $G$ and kernel $Q$), in order to recover the input from the output. It shows that depending on the recovery requirements, different phase transitions take place:  For exact recovery, the CH threshold is efficiently achievable for any fixed number of communities. For weak recovery, the KS threshold is efficiently achievable for any fixed number of communities, but it is not necessarily the information-theoretic threshold, leaving a question mark on whether the KS threshold is indeed the fundamental limit for efficient algorithms. We also presented partial results on the optimal tradeoffs between various measures of distortion and the SNR in the partial recovery regime. 

In the quest to achieve these thresholds, novel algorithmic ideas emerged, similarly to the quest to achieve the capacity in channel coding, with sphere-comparison, graph-splitting, linearized BP, nonbacktracking operators and graph powering.  
This program can now be pursued in different directions, refining the models, improving the algorithms and expanding the realm of applications. In particular, similar tradeoffs are expected to take place in other graphical channels, such as ranking, synchronization, topic modelling, collaboration filtering, planted embeddings and more. We list below a series of possible open problems. 
%By understanding various instances of graphical channels, the hope is to build a general theory for the fundamental limits in mining graphical data. 

%Community detection has connections with information theory at various levels. 
%Exact recovery is closely related to the decoding of graph-based codes on memoryless channels, and to $f$-divergences. Weak recovery relates naturally to the broadcasting problem on trees, and partial recovery to information-estimation measures. In \cite{abbetoc}, the SBM is viewed as a special case of a more general family of {\it graphical channels:}

\begin{itemize}
\setlength{\itemsep}{3pt}
\item {\it Exact recovery for sub-linear communities.} The survey gives a comprehensive treatment for exact recovery with linear-size communities, i.e., when the entries of $p$ and its dimension $k$ do not scale with $n$. If $k=o(\log(n))$, most of the developed techniques tend to extend. What happens for larger $k$? In \cite{afonso3,chen-xu}, some of this is captured by looking at coarse regimes of the parameters. It would be interesting to pursue sub-linear communities in the lens of phase transitions and information-computation gaps.

\item {\it Partial recovery.} What is the fundamental tradeoff between the SNR and the distortion (MMSE,  agreement or mutual information) for partial recovery and arbitrary constant degrees? As a preliminary result, one may attempt to show that $I(X;G)/n$ admits a limit in the constant degree regime. This is proved in \cite{abbetoc} for two symmetric disassortative communities, but the assortative case remains open. A recent result from \cite{coja2} further gives the expression for the limit in the disassortative case, but the assortative case remains open. %Establishing the expression for the optimal tradeoff in partial recovery is unknown for $k\ge 3$ (and only known for large enough SNR for $k=2$ \cite{mossel2}).

\item {\it The information-computation gap:} 
\begin{itemize}
\item Related to the last point; can we locate the exact information-theoretic threshold for weak recovery when $k \geq 3$? Recent results and precise conjectures were recently obtained in \cite{lelarge_ks}, for the regime of finite SNR with diverging degrees discussed in Section \ref{wigner}. Arbitrary constant degrees remain open.

\item Can we strengthen the evidence that the KS threshold is the computational threshold? In the general sparse SBM, this corresponds to the following conjecture:

\begin{conjecture}
Let $k \in \mZ_+$, $p\in (0,1)^k$ be a probability distribution, $Q$ be a $k\times k$ symmetric matrix with nonnegative entries. If $\lambda_2^2<\lambda_1$, then there is no polynomial time algorithm that can solve weak recovery (using Definition \ref{all-defs}) in $G$ drawn from $\sbm(n,p,Q/n)$.
\end{conjecture}
%Can we compare/interpolate the computational barriers taking place for community detection in the SBM and for the planted clique problem? 

%\item Can we compare/interpolate the computational barriers taking place for community detection in the SBM and for the planted clique problem? 

%\item What is the optimal fraction of misclassified vertices when $k\ge 3$ under efficient or infomration-theoreic algorithms? 
\end{itemize}

\item \looseness=-1{\it Learning the general sparse SBM.} Under what conditions can we learn the parameters in $\sbm(n,p,Q/n)$ efficiently or information-theoretically? 

\item \looseness=-1{\it Scaling laws:} What is the optimal scaling/exponents of the probability of error for the various recovery requirements? How large need the graph be, i.e., what is the scaling in $n$, so that the probability of error in the discussed results\footnote{Recent work \cite{young} has investigated finite size information-theoretic analysis for weak recovery.} is below a given threshold? 

\item {\it Beyond the SBM:} 
\begin{itemize}

\item How do previous results and open problems generalize to the extensions of SBMs with labels, degree-corrections, overlaps (see \cite{david1}), etc. In the related line of work for graphons \cite{choi,sbm_graphon,borgs_nips}, are there fundamental limits in learning the model or recovering the vertex parameters up to a given distortion? 
The approach of \cite{colin1focs} and sphere-comparison were generalized to the case of overlapping communities in \cite{collaboration_filtering} with applications to collaborative filtering. Can we establish fundamental limits and algorithms achieving the limits for other unsupervised machine learning problems, such as topic modelling, ranking, Gaussian mixture clustering (see \cite{banks3}), low-rank matrix recovery (see \cite{yash_pca} for sparse PCA) or general graphical channels?

\item How robust are the thresholds to model perturbations or adversaries? It was shown in \cite{maka,ankur_SBM} that monotone adversaries can interestingly shift the threshold for weak recovery; what is the threshold for such adversarial models or adversaries having a budget of edges to perturb? What are robust algorithms (see also Section \ref{robust})? What are the exact and weak recovery thresholds in geometric block models (see also previous section)?

\end{itemize}

\item {\it Semi-supervised extensions:} How do the fundamental limits change in a semi-supervised setting,\footnote{Partial results and experiments were obtained for a semi-supervised model \cite{semi}. Another setting with side-information is considered in \cite{clauset} with metadata available at the network vertices. Effects on the exact recovery threshold have also been recently investigated in \cite{asadi1}.} i.e., when some of the vertex labels are revealed, exactly or probabilistically?

\item {\it Dynamical extensions:} In some cases, the network may be dynamical and one may observe different time instances of the network. How does one integrate such dynamics to understand community detection?\footnote{Partial results were recently obtained in \cite{dynamic}.}  

%\item {\it Compression of graphical data:} What is the interplay between community detection and compression of graphical data? Specific results are obtained in \cite{abbe_allerton16,todo asadi1}, but 

\end{itemize}

%\chapter{Appendix}\label{app1}
%todo (Keep or leave)

\section*{acknowledgements}
I am thankful to my collaborators from the main papers discussed in this monograph, in particular, A.~Bandeira, Y.~Deshpande, J.~Fan, A.~Montanari, C.~Sandon, K.~Wang, Y.~Zhong. Special thanks to Colin and Kaizheng for their valuable inputs on this monograph. I would like to also thank the many colleagues with whom I had various conservations on the topic over the past years, in particular, C.~Bordenave, F.~Krzakala, M.~Lelarge, L.~Massouli\'e, C.~Moore, E.~Mossel, Y.~Peres, A.~Sly, V.~Vu, L.~Zdeborov\'a, as well as the colleagues, students and anonymous reviewers who contributed to the significant improvements of the drafts. Finally, I would like to thank the Institute for Advanced Study in Princeton for providing an ideal environment to write part of this monograph.

%BACKMATTER SEE DOCUMENTATION
%\backmatter  % references, restarts sample

%{\small
%\bibliographystyle{amsalpha}
%\bibliography{gen_sbm,mjwain_super}
%}

\newcommand{\etalchar}[1]{$^{#1}$}
\providecommand{\bysame}{\leavevmode\hbox to3em{\hrulefill}\thinspace}
\providecommand{\MR}{\relax\ifhmode\unskip\space\fi MR }
\providecommand{\MRhref}[2]{%
  \href{http://www.ams.org/mathscinet-getitem?mr=#1}{#2}
}
\providecommand{\href}[2]{#2}

\end{document}